\documentclass{amsart}

\usepackage{amsthm,amsmath,amssymb,enumitem,geometry,hyperref}
\usepackage{xcolor, cancel}
\hypersetup{
    colorlinks=true,   
    linkcolor=blue,    
    filecolor=magenta, 
    urlcolor=cyan,     
    citecolor=red,       
}

\makeatletter
\newsavebox{\@brx}
\newcommand{\llangle}[1][]{\savebox{\@brx}{\(\m@th{#1\langle}\)}%
	\mathopen{\copy\@brx\kern-0.5\wd\@brx\usebox{\@brx}}}
\newcommand{\rrangle}[1][]{\savebox{\@brx}{\(\m@th{#1\rangle}\)}%
	\mathclose{\copy\@brx\kern-0.5\wd\@brx\usebox{\@brx}}}
\makeatother

\usepackage{bbm}
\usepackage[utf8]{inputenc}
\usepackage[OT2,T1,T2A]{fontenc}
\usepackage{ucs}

\usepackage{substitutefont}
\substitutefont{OT2}{\rmdefault}{wncyr}
\DeclareRobustCommand{\rus}[1]{
	{\fontencoding{OT2}\selectfont#1}}

\newcommand{\1}{\mathbbm{1}}
\newcommand{\CL}{{\text{\fontencoding{T2A}\selectfont\CYRL}}}

\newcommand{\unmezzo}{\frac{1}{2}}

\newcommand{\eps}{\varepsilon}
\newcommand{\epsi}{\varepsilon}
\renewcommand{\epsilon}{\varepsilon}

\newcommand{\lhs}{left hand side}
\newcommand{\rhs}{right hand side }

\DeclareMathOperator{\deriv}{d}

\DeclareMathOperator{\dt}{dt}

\let\TeXchi\chi
\def\chi{{\setbox0 \hbox{\mathsurround0pt
			$\TeXchi$}\hbox{\raise\dp0 \copy0 }}}

\newcommand{\Q}{{\mathcal Q}}

\newcommand{\derivt}{\frac{\deriv\!{}}{\dt}}

\renewcommand{\u}{\textbf{u}}

\renewcommand{\v}{\textbf{v}}

\newcommand{\grad}{\nabla_x} 
\newcommand{\grady}{\nabla_y} 
\newcommand{\wlap}{\Delta_{\eta_\infty}} 

\newcommand{\R}{\mathbb{R}}

\newcommand{\sprod}[2]{\langle #1 , #2 \rangle} 

\newcommand{\rd}{\mathrm{d}}

\DeclareMathOperator{\dx}{dx}

\DeclareMathOperator{\dz}{dz}
\DeclareMathOperator{\ds}{ds}
\DeclareMathOperator{\dy}{dy}
\DeclareMathOperator{\dd}{d}

\DeclareMathOperator{\dmu}{d\mu}

\newcommand{\ints}[1]{\int_{\R}#1\dy} 

\newcommand{\intsv}[1]{\int_{\R}\int_{\mathbb{V}} #1\, d\mu(v)\dy }
\newcommand{\intv}[1]{\int_{\mathbb{V}} #1\, d\mu(v)}

\renewcommand{\phi}{\varphi}

\newcommand{\finf}{f_\infty}

\newcommand{\etainf}{\eta_{\infty}}

\def\T{\mathsf{T}}
\def\L{\mathsf{L}}
\def\RR{\mathsf{R}}

\def\A{\mathsf{A}}
\def\I{\mathbf{I}}

\newcommand{\calL}{{\mathcal L}}

\DeclareMathOperator{\x}{x}
\newcommand{\xdot}{\dot{\x}}
\newcommand{\bx}{{\bf x}}
\DeclareMathOperator{\sign}{sign}

\newcommand{\s}{\sign(y)}
\newcommand{\sy}{\s}

\newtheorem{thm}{Theorem}[]
\newtheorem*{thm-nn}{Theorem}
\newtheorem{lem}[thm]{Lemma}
\newtheorem{assumption}[thm]{Assumption}
\newtheorem*{lem-nn}{Lemma}
\newtheorem{prop}[thm]{Proposition}
\newtheorem{cor}[thm]{Corollary}
\newtheorem*{cor-nn}{Corollary}
\newtheorem{remark}{Remark}

\numberwithin{equation}{section}

\usepackage{orcidlink}

\begin{document}

\title[Stability analysis for a kinetic bacterial chemotaxis model]{Stability analysis for\\ a kinetic bacterial chemotaxis model}

\author[V.~Calvez]{Vincent Calvez\,\orcidlink{0000-0002-3674-1965}}
\address{CNRS, Univ Brest, UMR 6205, Laboratoire de Mathématiques de Bretagne Atlantique, France}
\email{vincent.calvez@math.cnrs.fr}

\author[G.~Favre]{Gianluca Favre\,\orcidlink{0000-0002-2090-2648}}
\address{Faculty of Mathematics, University of Vienna, Oskar-Morgenstern-Platz 1, 1090, Vienna, Austria}
\email{gianluca.favre@univie.ac.at}

\author[F.~Hoffmann]{Franca Hoffmann\,\orcidlink{0000-0002-1182-5521}}
\address{Department of Computing and Mathematical Sciences, California Institute of Technology, 1200 E California Blvd. MC 305-16, Pasadena, CA 91125, USA}
\email{franca.hoffmann@caltech.edu}

\keywords{Chemotaxis, run-and-tumble model, nonlinear stability, kinetic transport equation, hypocoercivity}

\begin{abstract}
We perform stability analysis of a kinetic bacterial chemotaxis model of bacterial self-organization, assuming that bacteria respond sharply to chemical signals. The resulting discontinuous tumbling kernel represents the key challenge for the stability analysis as it rules out a direct linearization of the nonlinear terms. To address this challenge we fruitfully separate the evolution of the shape of the cellular profile from its global motion. 
We provide a full nonlinear stability theorem in a perturbative setting when chemical degradation can be neglected. With chemical degradation we prove stability of the linearized operator. In both cases we obtain exponential relaxation to equilibrium with an explicit rate using hypocoercivity techniques. To apply a hypocoercivity approach in this setting, we develop two novel and specific approaches: i) the use of the $H^1$ norm instead of the $L^2$ norm, and ii) the treatment of nonlinear terms.  
This work represents an important step forward in bacterial chemotaxis modeling from a kinetic perspective as most results are currently only available for the macroscopic descriptions, which are usually parabolic in nature. Significant difficulty arises due to the lack of regularization of the kinetic transport operator as compared to the parabolic operator in the macroscopic scaling limit.
\end{abstract}

\maketitle



\section{Introduction, model and main results}
Bacteria such as \textit{Escherichia coli} can navigate in heterogeneous environments made of chemical gradients. Moreover, they can produce the chemical positive cues they are sensitive to. This leads to the possible aggregation of the cell population into a stable cluster of bacteria, as reported in the biophysical literature, see {\em e.g.} \cite{mittal_motility_2003}. This cluster results from the interplay between the self-attraction effect of chemotaxis, and the  dispersive effect of cell motion with random reorientations. Additional mechanisms can put these clusters in motion, such as local nutrient consumption, leading to remarkable collective propagation of the cell population, turning the mathematical stationary state into a traveling wave problem \cite{saragosti_directional_2011}. 

This work is devoted to the stability analysis of stationary clusters of bacteria under the effect of positive self-induced chemotaxis. Chemotaxis is subject to an intense mathematical analytical effort (see recent review \cite{Arumugam} and references therein). We shall put the focus at the mesoscopic level of description, using a kinetic-transport equation   as it is common for bacteria which have a persistent motion (run) with intermittent reorientation (tumble). Many results in the literature are dealing with macroscopic derivations of parabolic models from kinetic equations, with increasing level of accuracy, aligned with the accumulation of knowledge of how single bacteria navigate in a heterogeneous environment \cite{HiOth00, MR2120548, MR2250124}. However, quantitative results at the mesoscopic level \textit{per se} are much less available. 

Here, we focus on an underlying kinetic-transport equation  suitable to describe the individual run-and-tumble motion of bacteria. The phase space distribution of cells is denoted by $f(x,v;t):\R^d\times \left(\mathbb{V}\subset\R^d \right)\times[0,\infty) \to \R_{\ge 0}$, where the set of admissible velocities $\mathbb{V}$ is bounded and symmetric with respect to the origin and endowed with a probability measure $\dmu$, renormalized such that $\int_\mathbb{V}\dmu(v)=1$\footnote{Although stated here for general measures $\mu$, we should emphasize that we will restrict the analysis to the two velocity model
		$\dmu(v)=\frac12\delta_{-1}(v)+\frac12\delta_{+1}(v)$, in one dimension of space.
}. The corresponding macroscopic cell density is then given by $\rho(t,x)=\int_{\mathbb{V}}f(x,v;t)\dmu(v)$. The chemoattractant concentration $S(t,x)$ is generated by the cell density. This results in the following model for  cluster formation and stabilization,
\begin{subequations}
	\label{kinetic1}
	\begin{align}
		& f_t + v\cdot\grad f = \int_{\mathbb{V}} \big( T[S](t,x,v',v) f' - T[S](t,x,v,v') f \big)\dmu(v') \,,\label{kinetic1a}\\
		& - \Delta S(t,x) + \alpha S(t,x)  = \rho(t,x) \,,\label{kinetic1b}
	\end{align}
\end{subequations}
where the primed notation stands for the dependence on $v'$ instead of $v$.  The tumbling kernel $T[S](t,x,v',v)$ denotes the turning frequency for cells moving at velocity $v'$ to change direction due to a tumble event resulting in a new velocity $v$. Cells adapt their tumbling frequency locally in space based on the chemoattractant concentration $S(t,x)$. The parameter $\alpha\ge 0$ is the natural decay of the chemical component.
The system \eqref{kinetic1} is equipped with an initial condition $f(x,v;0) = f_0(x,v)$, whose regularity and decay at infinity will be discussed below.


This kinetic model is a particular case of \cite{calvez_chemotactic_2019}. The agreement with experimental data was shown to be satisfactory in \cite{saragosti_directional_2011}. Note that exponential relaxation to equilibrium was established in \cite{calvez_confinement_2015} for the linear problem, where the communication signal $S$ is prescribed, using hypocoercive techniques. This was extended in \cite{mischler_linear_2016} to any dimension of space for a given radially symmetric chemoattractant density $S$, and improved in \cite{evans21,evans24} removing the radial symmetry assumption and considering a wider class of tumbling
kernels including physically more relevant non-uniform kernels. 
In \cite{evans21} the authors further consider a chemoattractant density $S(t,x)$ with a fixed dependence on $\rho(t,x)$ that is more amenable to analysis than \eqref{kinetic1b} and can be considered as an interpolation between the linear setting and the full physically relevant non-linear model \eqref{kinetic1},
resulting in a non-linear model without a stiff non-Lipschitz term in the tumbling rate. This modification removes the need to assume that the chemoattractant density has a single peak. \cite{evans21} is the only other work in the literature we are aware of that provides convergence results for non-linear run and tumble equations of the form \eqref{kinetic1a}. 
Our work goes beyond \cite{evans21}, considers the full physically relevant model \eqref{kinetic1} and uses different techniques.
In particular, our tumbling kernel breaks the differentiability assumptions that are necessary for the convergence results in \cite{evans21} to hold. 
For the sake of simplicity, we will focus here only on the one-dimensional case ($d=1$) and on the two-velocity case.

\subsubsection*{The tumbling kernel.}

From now on, we consider the following choice of the tumbling kernel
\begin{align}\label{eq:tumbling kernel}
	T[S](t,x,v,v')=\sigma K(t,x,v)\,,\qquad
	K(t,x,v) = 1 - \chi \sign(v)\sign(\partial_xS(t,x))\,,
\end{align}
where $\chi\in(0,1)$ is the chemoattractant sensitivity and $\sigma>0$ denotes the relative intensity of the tumbling rate. In other words, under this modelling ansatz, we assume that the probability for a cell to change direction from $v$ to $v'$ only depends on the prior velocity. This kernel was introduced in \cite{saragosti_directional_2011}, and further studied in \cite{calvez_confinement_2015,calvez_chemotactic_2019}. The peculiar form of the $\sign$ function is motivated by biological evidence of efficient signal integration of molecular variations at the single cell level (see discussion in the aforementioned references, but note that more realistic, logarithmic gradient-sensing rules have been derived in the biophysical literature in the past decade, see {\em e.g.} \cite{Perthame2016} and references therein). Moreover, on the side of mathematical considerations, it is compulsory to have a sharp transition at small gradient amplitude, otherwise dispersion will dominate for cells that lie far away from the bulk, resulting in cell leakage hindering the existence of a stationary state. Finally, it is one of the advantages of the choice \eqref{eq:tumbling kernel} that the existence of  stationary states is tractable \cite{calvez_confinement_2015}. They are even explicit in the case of two-velocities considered in this work. Existence of a stationary state is expected also for the logarithmic gradient sensing rule, but this is far beyond the scope of the present work.

The choice of the $\sign$ function, however, comes with a serious price to pay regarding stability analysis, as it rules out a direct linearization of the nonlinear term, obviously. On the other hand, $K$ is a piecewise linear function, which is fully determined by the knowledge of those points where  $\partial_x S$ changes sign. We base our strategy on the latter observation. For our results to hold, we make the following assumption throughout:

\begin{gather}\label{H1}\tag{H1}
	\text{We assume that $S(t,x)$ has a unique critical point}\\
	\text{(which is a global maximum) at all times $t\ge 0$.}\notag
\end{gather}
This assumption enables us to define uniquely the point $\bx(t)$ where $\partial_x S(t,x)$ changes sign: 
\begin{equation}\label{defX}
	\begin{cases}
		\sign(\partial_x S(t,x)) = 1\,, & \quad x<\bx(t)\,,\\
		\sign(\partial_x S(t,x)) = -1\,, & \quad x>\bx(t)\,.\\
	\end{cases}
\end{equation} 
Consequently, the tumbling kernel simplifies to
\begin{align*}
	K(t,x,v) = 1 + \chi \sign(v) \sign(x-\bx(t))\,,
\end{align*}
By changing reference frame $ y= x-\bx(t)$, we can fruitfully separate the shape of the density profile, from the movement of its center $\bx(t)$.  Two of the authors already applied this approach successfully in the macroscopic case \cite{CalHof20}. However, much difficulty arises here due to the lack of regularization of the kinetic transport operator, as compared to the parabolic operator in \cite{CalHof20}. 


\subsection{Main results}

For most turning mechanisms $K$ and velocity measures $\mu$, it is not straight forward to identify a steady state $f_\infty$, and an explicit expression may not be found. This motivates the two velocity model as shown in the following proposition.
\begin{prop}[Existence of steady states]\label{prop:existence-steady-state}
	Under Hypothesis~\eqref{H1}, and in the case of the two velocity model
	\begin{equation*}
		\dmu(v)=\frac12\delta_{-1}(v)+\frac12\delta_{+1}(v)\,,\qquad \mathbb{V}=\{+1,-1\},
	\end{equation*}
	there exists a family of steady states 
 $$\rd f_\infty(x,v)=\frac{M\sigma\chi}{2}e^{-\sigma\chi|x-\bx_\infty|}\rd x\rd\mu(v)$$
with mass $M$ and center $\bx_\infty$ for the system \eqref{kinetic1} with tumbling kernel \eqref{eq:tumbling kernel}. Moreover, the choice $\dmu=\frac12\delta_{-1}+\frac12\delta_{+1}$ is the only velocity measure $\mu$ for which the steady state $f_\infty$ is independent of velocities up to dilatations.
\end{prop}

\begin{proof}
    The expression for $f_\infty(x,v)$ in the case of the two velocity model can be checked by direct substitution. To show the last statement of the proposition, note that at steady state we have $\xdot=0$, and so for $\bx(t)=\bx_\infty$ any steady state $f_\infty(x,v)$ satisfies
	\begin{align*}
		v\partial_x f_\infty &= \sigma \int_{\mathbb{V}}  K(x,v') f_\infty'\dmu(v') - \sigma K(x,v) f_\infty\\
		&= \sigma \left(\int_{\mathbb{V}} f_\infty'\dmu(v')-f_\infty\right)
		+\sigma\chi\sign(x-\bx_\infty)\left(\int_{\mathbb{V}} \sign(v')f_\infty'\dmu(v')-\sign(v)f_\infty\right)\,.
	\end{align*}
	For a steady state $f_\infty$ that is independent of velocities, this reduces to
	\begin{align*}
		v\partial_x f_\infty 
		&= \sigma\chi\sign(x-\bx_\infty)f_\infty\left(\int_{\mathbb{V}} \sign(v')\dmu(v')-\sign(v)\right)
		= - \sigma\chi\sign(x-\bx_\infty)\sign(v)f_\infty\,.
	\end{align*}
	Rearranging and integrating in space, we have
	\begin{align*}
		\log f_\infty = -\sigma\chi\frac{|x-\bx_\infty|}{|v|} +g(v)
	\end{align*}
	for some function $g:\mathbb{V}\to \R$ independent of $x$. From this expression it becomes clear that $f_\infty$ is only independent of the velocity variable $v$ if $|v|$ is constant for all $v\in\mathbb{V}$. In one dimension, this means $v=\pm c$ for some constant $c$, and by rescaling the velocities, without loss of generality, we arrive at $g\equiv 0$, $\mathbb{V}=\{-1,+1\}$ and $\dmu=\frac12\delta_{-1}+\frac12\delta_{+1}$, ensuring that $\int_\mathbb{V}\dmu(v)=1$. In this case, we recover $\rd f_\infty(x,v)=\frac{M\sigma\chi}{2}e^{-\sigma\chi|x-\bx_\infty|}\rd x\rd \mu(v)$.
\end{proof}
For simplicity and without loss of generality (as it corresponds to rescaling the density and time and space variables accordingly), we will from now on consider mass $M=1/\chi$ and relative intensity $\sigma=2$, and we will denote $\lambda:=2\chi+\sqrt{\alpha}$ throughout. In this case, we write the steady state as
\begin{equation*}
    \rd f_\infty(x,v)=\eta_\infty(x-\bx_\infty)\rd x\rd\mu(v)\,,\qquad \eta_\infty(y):= e^{-2\chi|y|} \,.
\end{equation*}
All arguments also hold for general choices of $M$ and $\sigma$. 
We denote by $H^1$ the weighted space of relative energy equipped with the following norm
\begin{equation}\label{eq:defH1}
	\|g\|_{H^1}
	=\left( \iint \left( \left|\frac{g}{\eta_\infty}\right|^2 + \left|\partial_x\left(\frac{g}{\eta_\infty}\right)\right|^2\right)  \eta_\infty \,\rd x\rd \mu(v)\right)^{1/2}
	\,.  
\end{equation}

\begin{thm}[Nonlinear stability: case $\alpha=0$]\label{thm:main1}
	Let $\alpha=0$.
	Under Hypothesis~\eqref{H1}, and for $\dmu(v)=\frac12\delta_{-1}(v)+\frac12\delta_{+1}(v)$, the steady state $\rd f_\infty(x,v)$ is locally nonlinearly stable in the following sense: There exists $\eps_0=\eps_0(\chi)>0$, $C_0>0$ and $\gamma>0$ such that, for all initial data $f_0$ of mass $1/\chi$ satisfying $\iint v  f_0(x,v)\dx\dmu(v)=0$ and
  \begin{equation}\label{ass:initial}
    \left\|\frac{f_0- \eta_\infty(\cdot-\bx(0))}{\eta_\infty (\cdot-\bx(0))}\right\|_\infty
     + 
     \left\|\partial_x\left(\frac{f_0-\eta_\infty(\cdot-\bx(0))}{\eta_\infty (\cdot-\bx(0))}\right)\right\|_\infty
     \le \eps_0\,,
 \end{equation}
	we have for all $t>0$,
	\begin{equation}\label{eq:main1}
		\|f(x,v;t) - \eta_\infty\big(x-\bx(t)\big)\|_{H^1} \leq C_0 e^{-\gamma t}\|f_0 (x,v;t)- \eta_\infty\big(x- \bx(0)\big)\|_{H^1}\, , \qquad \lim_{t\to \infty} \bx(t) = \bx_\infty\,,
	\end{equation}
	where $C_0>1$. Moreover, the limit $\bx_\infty$ exists but has no explicit value, up to our knowledge.
\end{thm}

\begin{remark}
The rate $\gamma$ is quantitative (and could me made explicit in the parameter $\chi$), as it is usual when applying hypocoercivity methods.
\end{remark}

Next, we linearise the operator around the steady state $\rd f_\infty(x,v)=\eta_\infty(x-\bx_\infty)\rd x\dmu(v)$, and obtain a linear stability result for any $\alpha>0$.

\begin{thm}[Linear stability: case  $\alpha>0$]\label{thm:main2}
	Let $\alpha>0$.
	Under Hypothesis~\eqref{H1}, and for $\dmu(v)=\frac12\delta_{-1}(v)+\frac12\delta_{+1}(v)$,
 the steady state $\rd f_\infty(x,v)$ is linearly stable in the following sense: There exists $\eps_\alpha=\eps_\alpha(\chi)>0$, $C_\alpha = C_\alpha(\chi) > 0$ and $\gamma>0$ such that, for all initial data $f_0$ of mass $1/\chi$ satisfying $\iint vf_0e^{-\sqrt{\alpha}|x-\bx(0)|}\dx\dmu(v)=0$ and
	$$\|f_0(x,v) - \eta_\infty\big(x-\bx(0)\big) \|_{H^1}\leq  \eps_\alpha\,,$$
	we have for all $t>0$, up to nonlinear perturbations,
		\begin{equation}\label{eq:main2}
			\|f(x,v;t) - \eta_\infty\big(x-\bx(t)\big)\|_{H^1} \leq C_\alpha e^{-\gamma t}\|f_0(x,v) - \eta_\infty\big(x-\bx(0)\big)\|_{H^1}\, , \qquad \lim_{t\to \infty} \bx(t) = \bx_\infty\,,
		\end{equation}
		where $C_\alpha = C_\alpha(\chi) > 0$.
		Moreover, the limit $\bx_\infty$ exists but has no explicit value, up to our knowledge. 
\end{thm}

\begin{remark}
The assumption $\iint vf_0e^{-\sqrt{\alpha}|x-\bx(0)|}\dx\dmu(v)=0$ can be relaxed. Indeed, if this expression is non-zero initially, Lemma~\ref{lem:cond-initial} guarantees that it decays exponentially in time, and therefore all estimates can be extended to still hold with an exponentially small perturbation; see Remark~\ref{rmk:IC} for more details.
\end{remark}

\subsection{Q\&A Session}

\subsubsection*{Why should we reduce our analysis to the two-velocity model in the one dimensional spatial dimension?} 
\begin{itemize}
    \item The steady state generally depends on the velocity set  $\mathbb{V}$ and this dependence is not explicit in principle, see \cite{calvez_confinement_2015,calvez_chemotactic_2019}, except in this special case $d=1$ and $\mathbb{V} = \{-1,+1\}$ as shown in Proposition~\ref{prop:existence-steady-state}. Achieving results for non-linear settings where the steady state is not known is an interesting and relevant open problem, with recent progress for weakly non-linear run and tumble equations \cite{evans21}.
\item A more serious obstacle arises in a crucial part of the proof. In fact, we work in the moving frame $y = x - \bx(t)$. We run perturbative estimates under the key condition that $\dot \bx$ remains uniformly small, comparably smaller than the speed of individual bacteria, actually. Therefore, allowing for arbitrarily small velocities 
would immediately break down our strategy. 
\end{itemize}

\subsubsection*{Why should we split our result between $\alpha=0$ (nonlinear stability) and $\alpha>0$ (linear stability)?}
\begin{itemize}
\item We realized in a previous work by two of the authors \cite{CalHof20} that the case $\alpha>0$ can induce some additional difficulty due to various exponential weights in the corresponding conservation laws (conservation of mass, and conservation of the center of distribution in the moving frame). Here, we expect the situation to be much more challenging, which is why we proceed with the strategy outlined in the next point.
    \item We provide a full nonlinear stability theorem in a perturbative setting, under the conditional hypothesis \eqref{H1}, in the simpler case $\alpha = 0$. Then, we investigate the case $\alpha > 0$, for which we restrict to a linear stability result to focus on the structure of the linearized operator.  
\end{itemize}

\subsubsection*{How reasonable is Hypothesis~\eqref{H1}?}
\begin{itemize}
\item The condition appearing in Hypothesis \eqref{H1} lies at the heart of analytical results concerning the piecewise constant tumbling kernel \eqref{eq:tumbling kernel}, see \cite{saragosti_mathematical_2010,calvez_confinement_2015,calvez_chemotactic_2019,CalHof20}.  At the level of the macroscopic problem, it was shown to propagate along the flow when it is true at initial time \cite[Proposition 1.1]{CalHof20}. At the level of the mesoscopic problem, the existence of stationary states for which this condition holds true was shown in a fairly general setting \cite{calvez_chemotactic_2019}, far beyond the two-velocity case. However, we are lacking a nice argument as in \cite{CalHof20} in the present setting to propagate the shape of $S$ along the flow. The reason is that we could not derive a tractable equation for the derivative of $\partial_x S(t,x)$, as done in the proof of \cite[Proposition 1.1]{CalHof20}. 
\item The case $\alpha=0$  is an exception, as the condition appearing in Hypothesis \eqref{H1} is always satisfied, simply because $\partial_{xx}S(t,x)= -\rho(t,x)\le 0$ for all $(t,x)$, and so $S$ is concave in space for all times.  
\item The following example is meant to show that the condition is not expected to be propagated in full generality because of inertia in the kinetic transport equation: Let $\rho(0,x)$ be a standard Gaussian. Then $S(0,x)$ has a unique critical point, which is a global maximum (see \cite[Proposition 5.1]{calvez_chemotactic_2019}). Consider that $\rho(0,x)$ is generated by the underlying initial kinetic density
 $$f(0,x,v)=\rho(0,x)\left(\mathrm{1}_{x<0}\,\delta_{v=-1}+ \mathrm{1}_{x>0}\,\delta_{v=1} \right)\,.$$ Then after a small time $t>0$, the spacial density $\rho(t,x)$ is expected to form two peaks, resulting in two peaks for $S(t,x)$, provided $\alpha$ is large enough.
\end{itemize}

\subsubsection*{How do our results compare with current state-of-the-art?}\label{Sect: sota}
\begin{itemize}
    \item This work represents an important step forward in bacterial chemotaxis modelling from a kinetic perspective. The mathematical literature in bacterial chemotaxis is overwhelmingly dominated by questions about blow-up vs. global existence of solutions for the macroscopic descriptions, which are usually parabolic in nature. However, kinetic equations have been shown to be an important tool for capturing collective behaviour at the mesoscopic scale \cite{calvez_chemotactic_2019}. Although the blow-up issue can also be formulated for kinetic chemotaxis equations \cite{Bournaveas}, uniformly bounded, stationary equations are more meaningful for studying bacterial collective at this scale.
    
    \item This work builds on the hypocoercivity technique as presented in \cite{DMS15}, but develops two novel and specific approaches. The first involves the use of the $H^1$ norm instead of the $L^2$ norm. This kind of approach was introduced in a general way for linear problems by Mouhot and Neumann in \cite{Mouhot-Neumann}. The second approach addresses the treatment of nonlinear terms, as introduced by Favre, Pirner, and Schmeiser in \cite{favre_nonlinhypo}. Furthermore, these two aspects are intertwined, as it is in the nonlinear term that higher-order terms emerge.  Up to our knowledge, this is the first result about $H^1$ nonlinear hypocoercivity.
    
    \item Here, we take advantage that the nonlinearity occurs only through a scalar component $\bx(t)$, which enables to split the problem into subparts with different features (PDE \& ODE). Of note, this type of structure is not isolated, an important example for instance is the dynamics of neuron populations whose individual firing rate depends on the global activity of the network, through the (scalar) total intensity of firing neurons \cite{Caceres}; another important class of models is mutation-selection models in quantitative genetics where individuals are in competition through the consumption of a single resource \cite{Barles}. Yet, despite many simplifications, this should not be considered a toy problem, as many difficulties arise in the coupling due to the lack of regularity of the hyperbolic operator.
\end{itemize}

\subsection{Strategy of the proof}\label{sec:context}

We borrow some key ideas from the macroscopic case  \cite{CalHof20}. In view of Hypothesis~\eqref{H1} and \eqref{defX}, it is natural to reformulate \eqref{kinetic1} in the moving frame $y = x - \bx(t)$. For $$(\bar f(y,v;t), \bar S(t,y)):=(f(x,v;t),S(t,x))\,,$$
and with corresponding steady state $(\bar f_\infty(y,v), \bar S_\infty(y))$.
After performing the change of variable   $y= x-\bx(t)$, we derive energy estimate of order 0 and order 1 for the density $f$. We apply the operator splitting {\em à la} Dolbeault-Mouhot-Schmeiser \cite{DMS15}, since it is well suited for kinetic equations. Namely, we can write (in abstract form), 
\begin{equation}
	\label{hypo eq formulation-intro}
	\partial_t W_y + \T W_y = \L W_y + \RR(W)\,, 
\end{equation}
where $W$ collects some proxies of the cell density (after appropriate reformulation), and $W_y$ denotes the derivative along the spatial coordinate. 
Here, $\T$ and $\L$ encode linear transport and collision terms respectively, as introduced in \cite{calvez_confinement_2015}, whereas $\RR$ collects nonlinear and nonlocal remainder terms depending on $W$, $W_y$ and $W_{yy}$.
By investigating the asymptotic behaviour of the semigroup generated by $\L-\T$, our goal is to quantify its stability or, to be precise, to determine the rate of convergence of $W$ towards zero as $t\to\infty$. Since the operator $\L-\T$ acts linearly on $W_y$, there is no restriction to study fluctuations around a global equilibrium, and we can apply the general methodology of hypoercoercivity for linear kinetic equations conserving mass as outlined in \cite{DMS15}. The important difference is that here, we are working in $H^1$ as opposed to $L^2$ as it is done in the classical setting. The reason for choosing $H^1$ is because we need to control the dynamics of $\bx(t)$, which contains pointwise values of $f$ at $x=\bx(t)$ which cannot be evaluated by energy estimates of order 0. 
More precisely,
we have to control separately $f(\bx(t),+1;t)$ and $f(\bx(t),-1;t)$ and their derivatives. 
The difference of $f(\bx(t),+1;t)$ and $f(\bx(t),-1;t)$ can be readily controlled by microscopic dissipation, but to control each separately we need to use the macroscopic dissipation. This requires the careful combination of microscopic and macroscopic estimates. Achieving such a control is one of the main contributions in this work, allowing the application of the hypocoercivity approach to this nonlinear model.

The paper is organized as follows: In Section \ref{Sect: 2}, we carefully describe the modeling set-up, separating the evolution of the shape of the solution from the evolution of the chemoattractant peak $\bx(t)$ and writing the two-velocity model in the perturbative setting close to the steady state; finally, we discuss conservation laws and notation used throughout the paper.
Section \ref{Sect: 3} provides an explicit expression for the $H^1$-entropy dissipation which is the starting point of our hypocoercivity analysis; the technical proof is postponed to Appendix~\ref{appendix: technical proof}. Section~\ref{Sect: 4} is dedicated to several functional inequalities which is the ground work for the later hypocoercivity estimates as well as for the control of $\xdot(t)$ in Section~\ref{sec:xdot}.
Section \ref{Sect: 5} is focused on hypocoercivity, a key mathematical tool to achieve convergence to equilibrium. A Lyapunov (\rus{Lyapunov}) functional is introduced modifying the entropy of the system, and controling its dissipation requires careful estimation of the linear and nonlinear contributions in combination with the functional inequalities from Section~\ref{Sect: 4}.
Subsequently, in Section \ref{Sect: 6}, we apply the previous machinery to analyse the nonlinear case with the parameter $\alpha = 0$, for which we need a detailed examination of entropy decay and the control of nonlocal terms. Since in this case the weights of the involved functional spaces are the same,
we are able to provide a nonlinear stability analysis, concluding with the proof of Theorem~\ref{thm:main1} in Section~\ref{sec:proof-main1}. Section \ref{Sect: 7} treats the linear case for $\alpha > 0$.  
To prove Theorem~\ref{thm:main2}, we split the argument in two cases: 1) we show the convergence of the shape profile centered at $\xdot$ applying hypocoercivity techniques to a different modified entropy that includes an $\alpha$-dependent nonlocal term, and 2) we show convergence of the chemoattractant peak using the first step.
Technical details of some of the proofs from across all sections are postponed to Appendices~\ref{appendix: technical proof} and \ref{sec:char-method}.


\section{Preliminaries}\label{Sect: 2}

\subsection{Set-up}

\subsubsection*{Moving frame.}
We consider the change of reference frame
\[y = x - \bx(t)\]
For $(\bar f(y,v;t), \bar S(t,y)):=(f(x,v;t),S(t,x))$, and writing $(f,S)$ instead of $(\bar f, \bar S)$ for brevity, equation \eqref{kinetic1} simplifies to 
\begin{subequations}
	\label{kinetic2}
	\begin{align}
		& 	f_t + (v-\xdot) \partial_y f=\sigma \int_{\mathbb{V}} \big( K(y,v') f' - K(y,v) f \big)\dmu(v') \label{kinetic2a}\\
		&\hspace{2.33cm}
		= \sigma \int_{\mathbb{V}}  K(y,v') f'\dmu(v') - \sigma K(y,v) f \,,
		\notag\\
		& - \partial_{yy} S(t,y) + \alpha S(t,y)  = \rho(t,y)\,,\label{kinetic2b}\\
		& K(y,v) = 1 + \chi \sign(y)\sign(v)\label{kinetic2c}\,.
	\end{align}
\end{subequations}
Note that the evolution of the kinetic cell density $f$ is still coupled to the chemoattractant $S$ via the evolution of the chemoattractant peak $\bx(t)$.

\subsubsection*{Dynamics of the chemoattractant peak.}

The dynamics of $\bx(t)$ inherit from the condition $\partial_y S(t,0) = 0$ by using the explicit representation formula for $S$ in terms of $\rho$. Indeed, we can write $S = G_\alpha*\rho$ for $\alpha>0$, where $G_\alpha$ is the fundamental solution of \eqref{kinetic2b}: 
\begin{equation}\label{eq:S conv}
	S(t,y)
	= \frac{1}{2\sqrt{\alpha}}
	\int_{-\infty}^\infty  e^{- \sqrt{\alpha} |y-z|}  \rho(t,z) \, \dz\, .
\end{equation}
\begin{lem}\label{lem:xdot}
	The dynamics of the chemoattractant peak $\xdot=\xdot[f]$ are given by
	\begin{equation*}
		\xdot(t) = \frac{\int_{\mathbb{V}}v f (0,v;t)d\mu(v) - \frac{\sqrt{\alpha}}{2}\int_{\R\times\mathbb{V}} vfe^{-\sqrt{\alpha}|y|}\dy d\mu(v)}{\rho(0;t)- \frac{\sqrt{\alpha}}{2}\int_{\R}\rho e^{-\sqrt{\alpha}|y|}\dy} \,.
	\end{equation*}
\end{lem}
\begin{proof}
	To derive an equation for $\xdot$, we integrate \eqref{kinetic2a} over velocities to obtain the evolution of the mass
	\begin{equation}
		\label{eq:massevolution}
		\rho_t + \grady\cdot\int_{\mathbb{V}}{(v - \xdot) f}\dmu(v) =0\,.
	\end{equation}
	Using the explicit expression for $S$, we can rewrite the condition $\partial_y S(t,0)=0$ as
	\begin{equation*}
		0 = \ints{\sign(y) e^{-\sqrt{\alpha}|y|} \rho(t,y)}\,.
	\end{equation*}
 Differentiating in time, we get
	\begin{align*}
		0 &= \unmezzo \ints{\sign(y) e^{-\sqrt{\alpha}|y|} \partial_t\rho(t,y)}\\
		&= -\unmezzo \ints{\sign(y) e^{-\sqrt{\alpha}|y|} \intv{(\xdot-v)\partial_y f(y,v,t)}}\\
		&=\intv{(v-\xdot)f(0,v,t)} - \unmezzo\sqrt{\alpha} \intsv{e^{-\sqrt{\alpha}|y|}(v-\xdot)f(y,v,t)}\\
		&= - \xdot \rho(0,t) + \xdot \unmezzo \sqrt{\alpha} \ints{e^{-\sqrt{\alpha}|y|}\rho(y,t)}\\
		&\quad+\intv{v f(0,y,t)} - \unmezzo \sqrt{\alpha} \intsv{e^{-\sqrt{\alpha}|y|}vf(y,v,t)}\,.
	\end{align*}
	This completes the proof.
\end{proof}
For a detailed derivation of the chemoattractant peak dynamics of the corresponding macroscopic model, see \cite[Lemma 2.1]{CalHof20}.

\subsection{Two velocity model}
We reduce $\mathbb{V}$ to $\{-1,1\}$ by choosing
\begin{equation*}
	\dmu(v)=\frac12\delta_{-1}(v)+\frac12\delta_{+1}(v)\,.
\end{equation*}
Together with the choice of tumbling kernel~\eqref{kinetic2c}, the kinetic cell evolution~\eqref{kinetic2a} simplifies to
\begin{subequations}
	\label{kinetic4}
	\begin{align}
		\partial_t f_+ + (1-\xdot) \partial_y f_+ &=  (1-\chi \sign(y)) f_- - (1+\chi \sy) f_+ \,,\\
		\partial_t f_- - (1+\xdot) \partial_y f_- &=  \frac{\sigma}{2}(1+\chi \sign(y)) f_+ - \frac{\sigma}{2}(1-\chi \sy) f_-\,,
	\end{align}
\end{subequations}
for $f_{\pm}(y,t):=f(y,\pm 1;t)$.  
Next, we renormalize the two-velocity model by using instead 
$$
\phi(y;t) = \frac{f(y,1;t)}{f_\infty(y,1)}\,,\qquad
\psi(y;t) = \frac{f(y,-1;t)}{f_\infty(y,-1)}\,,
$$
and reducing model \eqref{kinetic2} to 
\begin{subequations}
	\label{kinetic5}
	\begin{align}
		\partial_t \phi + (1-\xdot) \partial_y \phi &= \phi (\xdot-1) \partial_y\ln\finf|_{v=+1} + (1-\chi \sign(y)) \psi - (1+\chi \sy) \phi \,,\\
		\partial_t \psi - (1+\xdot) \partial_y \psi &= \psi (\xdot+1) \partial_y\ln\finf|_{v=-1} + (1+\chi \sign(y)) \phi - (1-\chi \sy) \psi\,.
	\end{align}
\end{subequations}
Observe that in this case, the steady state of \eqref{kinetic4} can be determined explicitly: for initial datum $f_0$ of mass $1/\chi$, it is given by 
$$
f_\infty(y,+1)=f_\infty(y,-1)= \eta_\infty\,,\quad
\eta_\infty:=e^{-2\chi|y|}\,,\quad \xdot=0\,,
$$
with $\int_\R \eta_\infty(y)\dy =1/\chi$.
Substituting the explicit expression for the steady state, the normalized two velocity model \eqref{kinetic5} reduces to
\begin{subequations}
	\label{kinetic6}
	\begin{align}
		\partial_t \phi + (1-\xdot) \partial_y\phi &= -2\xdot\chi\sy \phi - (1-\chi \sign(y)) (\phi - \psi)\,, \\
		\partial_t \psi - (1+\xdot) \partial_y \psi &= -2\xdot\chi\sy \psi + (1+\chi \sign(y)) (\phi - \psi)\,,
	\end{align}
\end{subequations}
with
\begin{equation*}
	\label{xdot}
	\xdot = \frac{\phi(0;t) - \psi(0;t) - \frac{\sqrt{\alpha}}{2} \int_{\R^d} (\phi-\psi) e^{-(2\chi+\sqrt{\alpha})|y|} \dy   }{\phi(0;t)+\phi(0;t) - \frac{\sqrt{\alpha}}{2} \int_{\R^d} (\phi+\psi)e^{-(2\chi+\sqrt{\alpha})|y|} \dy } \,.
\end{equation*}
Model~\eqref{kinetic6} can be written in more compact form as
\begin{subequations}
	\label{kinetic7}
	\begin{align}
		\partial_t \phi + \partial_y\phi &= \xdot\eta_\infty^{-1}\partial_y\left(\phi\eta_\infty\right) - (1-\chi \sign(y)) (\phi - \psi)\,, \\
		\partial_t \psi - \partial_y \psi &= \xdot\eta_\infty^{-1}\partial_y\left(\psi\eta_\infty\right) + (1+\chi \sign(y)) (\phi - \psi)\,.
	\end{align}
\end{subequations}

\subsection{Perturbative Setting}

As we are concerned with convergence to equilibrium when we are close to a steady state, we consider a linear perturbation around the steady state,
\[\phi(y;t) = 1 + u(y;t)\quad,\quad \psi(y;t) = 1 + v(y;t)\,. \]
The evolution of $\xdot$ in terms of the perturbations $u,v$ is then given by
\begin{equation}
	\label{xdot pert}
	\xdot(t) = \frac{u(0;t) - v(0;t) - \frac{\sqrt{\alpha}}{2} \int_{\R^d} (u-v) e^{-(2\chi+\sqrt{\alpha})|y|} \dy   }{\frac{4\chi}{2\chi+\sqrt{\alpha}}+u(0;t)+v(0;t) - \frac{\sqrt{\alpha}}{2} \int_{\R^d} (u+v)e^{-(2\chi+\sqrt{\alpha})|y|} \dy } \,.
\end{equation}
Model~\eqref{kinetic6} then becomes
\begin{subequations}
	\label{kinetic10}
	\begin{align}
		\label{kinetic10a}
		\partial_t u + (1-\xdot) \partial_y u &= -2\xdot\chi\sy u - (1-\chi \sign(y)) (u - v) -2\xdot\chi\sy\,, \\
		\label{kinetic10b}
		\partial_t v - (1+\xdot) \partial_y v &= -2\xdot\chi\sy v + (1+\chi \sign(y)) (u - v) -2\xdot\chi\sy\,.
	\end{align}
\end{subequations}
Analogous to \eqref{kinetic7}, we can rewrite equations \eqref{kinetic10} as follows
\begin{subequations}
	\label{kinetic11}
	\begin{align}
		\label{kinetic11aa}
		\partial_t u + \partial_y u &= \xdot \eta_\infty^{-1} \partial_y \big( (u+1 ) \eta_\infty \big) - (1 - \chi \sy) (u-v)\,,\\
		\label{kinetic11bb}    
		\partial_t v - \partial_y v &= \xdot \eta_\infty^{-1} \partial_y \big( (v+1 ) \eta_\infty \big) + (1 + \chi \sy) (u-v) \,.
	\end{align}
\end{subequations}
Differentiating system~\eqref{kinetic11}, we obtain for $u_y$ and $v_y$,
\begin{subequations}
	\label{kinetic13}
	\begin{align}
		\label{kinetic11a}
		\partial_t u_y + \partial_y u_y &= \xdot \partial_y\left(\eta_\infty^{-1} \partial_y \big( (u+1 ) \eta_\infty \big)\right) - (1 - \chi \sy) (u_y-v_y) +2\chi\delta_0(u-v)\,,\\
		\label{kinetic11b}    
		\partial_t v_y - \partial_y v_y &= \xdot \partial_y\left(\eta_\infty^{-1} \partial_y \big( (v+1 ) \eta_\infty \big)\right) + (1 + \chi \sy) (u_y-v_y) +2\chi\delta_0(u-v)\,.
	\end{align}
\end{subequations}

\subsection{Conservation laws}\label{sec:csl}
For a scalar $a\in\R^+$, we define the scalar product $$\sprod{f}{g}_a = \int fg \, e^{-a|y|} \dy\,,$$ with corresponding weighted average $$\langle f \rangle_a := \frac{a}{2} \sprod{f}{1}_a\,.$$
With this notation, we have $\langle 1 \rangle_a=1$ for any $a\in\R^+$. The kinetic model \eqref{kinetic2} has two conservation laws: (1) conservation of mass, and (2) invariance by translation.
Fixing the mass to $1/\chi>0$, these write in the moving frame
\begin{align*}
	\int_\R\int_\mathbb{V}f(y,v;t)\dmu(v)\dy=\frac{1}{\chi}\,,\qquad 
	\int_\R\int_\mathbb{V}\partial_y f(y,v;t) e^{-\sqrt{\alpha}|y|}\,\dmu(v)\dy=0\,.
\end{align*}
We can then derive the corresponding conservation laws for the normalizations $\phi$, $\psi$ from \eqref{kinetic4}.
\begin{lem}\label{lem:CSL}
 Denote $\lambda=2\chi+\sqrt{\alpha}$. Then the solution $(\phi,\psi)$ to \eqref{kinetic4} satisfies the two conservation laws
	\begin{align*}
		\langle\phi+\psi\rangle_{2\chi} = 2 \,,\qquad
		\langle\phi_y+\psi_y\rangle_{\lambda}=0\,. 
	\end{align*}
\end{lem}
\begin{proof}
	For the first conservation law, this result follows by direct substitution. For the second conservation law, we calculate directly
	\begin{align*}
		0&= \int_\R\int_\mathbb{V}\partial_y f(y,v;t) e^{-\sqrt{\alpha}|y|}\,\dmu(v)\dy
		= \frac12\int_\R \left(\partial_yf_++\partial_yf_-\right) e^{-\sqrt{\alpha}|y|}\,\dy\\
		&= \frac12\int_\R \left(\partial_y(\phi\eta_\infty)+\partial_y(\psi\eta_\infty)\right) e^{-\sqrt{\alpha}|y|}\,\dy\\
		&= \frac12\int_\R \left(\phi_y+\psi_y\right) e^{-\lambda|y|}\,\dy
		-
		\chi\int_\R \left(\phi+\psi\right) \sign(y) e^{-\lambda|y|}\,\dy\\
		&= \frac12\int_\R \left(\phi_y+\psi_y\right) e^{-\lambda|y|}\,\dy
		+\frac{\chi}{\lambda}\int_\R \left(\phi+\psi\right) \partial_y e^{-\lambda|y|}\,\dy\\
		&=\left( \frac12-\frac{\chi}{\lambda}\right)\int_\R \left(\phi_y+\psi_y\right) e^{-\lambda|y|}\,\dy\,,
	\end{align*}
	and the result follows since $\frac12-\frac{\chi}{\lambda}= \frac{\sqrt{\alpha}}{2\lambda}\neq 0$.
\end{proof}

\begin{remark}
	We remark that the second conservation law is equivalent to the statement of a centering frame $\partial_yS(0;t)=0$. To see this, observe from the proof of Lemma~\ref{lem:xdot} that $\partial_yS(0;t)=0$ corresponds to 
	\begin{align*}
		0 &= \ints{\sign(y) e^{-\sqrt{\alpha}|y|} \rho(y;t)}
		=  \int_\R\int_\mathbb{V} \sign(y)f(y,v;t) e^{-\sqrt{\alpha}|y|}\,\dmu(v)\dy\\
		&= -\frac{1}{\sqrt{\alpha}} \int_\R\int_\mathbb{V} f(y,v;t) \partial_ye^{-\sqrt{\alpha}|y|}\,\dmu(v)\dy
		=\frac{1}{\sqrt{\alpha}}\int_\R\int_\mathbb{V}\partial_y f(y,v;t) e^{-\sqrt{\alpha}|y|}\,\dmu(v)\dy\,.
	\end{align*}
\end{remark}
Finally, the following formulation of the two conservation laws for the perturbed densities $u,v$ solving \eqref{kinetic10} is a direct consequence of Lemma~\ref{lem:CSL}.
\begin{cor}\label{cor:CSL-pert}
	Denote $\lambda=2\chi+\sqrt{\alpha}$. Then the solution $(\phi,\psi)$ to \eqref{kinetic4} satisfies the two conservation laws
	\begin{align*}
		\langle u+ v\rangle_{2\chi} = 0 \,,\qquad
		\langle u_y+ v_y\rangle_{\lambda}=0\,. 
	\end{align*}
\end{cor}
Notably, there is a discrepancy between the weights,  here $e^{-2\chi|y|}$ and $e^{-(2\chi + \sqrt{\alpha})|y|}$. It was shown in \cite{CalHof20}, that the choice of weight $\eta_\infty=e^{-2\chi  |y|}$ for the entropy estimates, which is more restrictive in terms of the class of functions considered, comes without restrictions on $\alpha$. It is also the natural choice, given it corresponds to the steady state of the system.

\begin{lem}\label{lem:cond-initial}
	Given the initial condition $\langle u - v \rangle_{\lambda}(0)$ then the solutions $u$ and $v$ to \eqref{kinetic10} satisfy
	\begin{equation}
		\label{eq: third cons law}
		\langle u - v \rangle_{\lambda}(t) = e^{-2t} \langle u - v \rangle_{\lambda}(0) \,.
	\end{equation}
\end{lem}
\begin{proof}
    From \eqref{kinetic10}, we deduce
    \begin{align*}
        \partial_t(u-v) + \partial_y(u+v) &= \xdot \eta_\infty^{-1} \partial_y \big( (u-v) \eta_\infty \big) - 2(u-v)\,,    
    \end{align*}
    and hence
    \begin{align*}
       \frac{\dd}{\dt} \langle u - v \rangle_{\lambda}
       = - \langle u_y + v_y \rangle_{\lambda}-2\langle u - v \rangle_{\lambda}=-2 \langle u - v \rangle_{\lambda}
    \end{align*}
    thanks to the second conservation law above.
\end{proof}
This third condition \eqref{eq: third cons law} is fundamental while considering the Poincar\'e inequality. In Section \ref{P inequality} in the case $\alpha=0$ we show that the Poincar\'e inequality in the $\exp{(-2\chi|y|)}$ weighted space holds up to a term that coincides exactly with $\langle u - v \rangle_{2\chi}$. This condition allows us to avoid this term entirely by choosing a suitable initial condition satisfying $\langle u - v \rangle_{2\chi}(0)=0$. It follows directly from the above lemma that $\langle u - v \rangle_{2\chi}(t)$ vanishes for all times if $\langle u - v \rangle_{2\chi}(0)=0$.

\subsection{Notation}

Denoting by $\Pi$ the projection on velocity averages, we  write
\begin{align}\label{def:Pi}
	W := 
	\begin{bmatrix}
		u\\v
	\end{bmatrix}\,,\qquad
	\Pi W     := \frac{1}{2}
	\begin{bmatrix}
		u+v\\u+v
	\end{bmatrix} =  \frac{1}{2}(u+v)\1\quad
 \text{ with }
 \1:=\begin{bmatrix}
		1\\1
	\end{bmatrix}\,.
\end{align}    
Further, let us introduce the following weighted inner product and corresponding norm
\begin{align*}
	\langle W_1,W_2\rangle := \langle u_1,u_2\rangle_{2\chi} + \langle v_1,v_2\rangle_{2\chi}
	= \int \left( u_1u_2+v_1v_2\right)\eta_\infty\dy\,,\quad
	\|W\|:=\langle W, W \rangle^{1/2}\,.
\end{align*}
Then
$$
\|(\I-\Pi)W\|^2=\frac12\int |u-v|^2\eta_\infty\dy\,,\qquad
\|\Pi W\|^2=\frac12\int |u+v|^2\eta_\infty\dy\,,
$$
and so
$$
\|W\|^2 = \|(\I-\Pi)W\|^2+\|\Pi W\|^2 = \int\left(|u|^2+ |v|^2\right)\eta_\infty\dy\,.
$$
We would like to estimate the decay of the $H^1$-norm $\|W_y\|^2$. Since the evolution of the derivatives $u_y$, $v_y$ contains singular terms, we first regularize equation~\eqref{kinetic13} before differentiating. We choose a suitable regularization kernel $K\in C^\infty(\R;\R_{\ge 0})$ satisfying 
\begin{align*}
	K(-y)=K(y)\,,\qquad \int_\R K(y)\,\dy=1\,,
\end{align*}
and define $u^\eps=K^\eps\ast u$, $v^\eps=K^\eps \ast v$, where $K^\eps(y):=\frac{1}{\eps}K\left(\frac{y}{\eps}\right)$ for any $\eps>0$. Further, we write for any function $f:\R\to\R$,
\begin{equation*}
	\llangle f \rrangle :=\lim_{\eps\to 0} \iint K(y)K(z)f(\eps(y-z))\,\dy\dz\,,
\end{equation*}
if the limit exists. Note that if $f$ is continuous at the origin, then $\llangle f\rrangle = f(0)$. If $f$ has right and left limits at the origin, then $\llangle f \rrangle= (f(0^+)+f(0^-))/2$. This notation allows to express the $H^1$-entropy dissipation in a nice compact form.
\section{Entropy Dissipation
}\label{Sect: 3}
In this short section we express explicitly the $H^1$-entropy dissipation.
\begin{prop}[$H^1$-entropy dissipation]\label{prop:H1dot}
	If $(u,v)$ solve equation~\eqref{kinetic11}, then
	\begin{align*}
		\frac12\frac{\dd}{\dt} \|W_y\|^2
		&=  - \int |u_y-v_y|^2\eta_\infty\dy
		+ 2\chi \left[u(0)-v(0)-2\xdot\right]\llangle u_y+v_y \rrangle 
		\\
		&\quad 
		- \frac{\xdot}{2} \int \partial_y \left(|u_y|^2+|v_y|^2\right)\eta_\infty\dy
		-4\chi\xdot\left[u(0)\llangle u_y\rrangle+ v(0)\llangle v_y\rrangle \right]\,.
	\end{align*}
\end{prop}
In the above expression for the dissipation, the first line corresponds to linear terms in equation~\eqref{kinetic10}, whereas the second line contains the nonlinear contributions.

The proof of Proposition~\ref{prop:H1dot} is rather technical, and we postpone it to Appendix~\ref{appendix: technical proof}. We give here a formal version. The main difficulty arises from the fact that $u_y(0)$ and $v_y(0)$ may not be well-defined as they are not necessarily continuous at zero. For the formal proof below we assume that $u_y$ and $v_y$ are continuous as the origin, an assumption that cannot be expected to hold for the dynamics~\eqref{kinetic11}. In Appendix \ref{appendix: technical proof} we rigorously derive the $H^1$-entropy dissipation using instead $\llangle u_y \rrangle$, $\llangle v_y \rrangle$ and $\llangle u_y + v_y \rrangle$ which exist and are well-defined.
\begin{proof}[Formal proof of Proposition~\ref{prop:H1dot}]
	We write the $L^2_{2\chi}$-norm of $W_y$ and we compute its time derivative
	\begin{align*}
		\unmezzo \derivt \sprod{W_y}{W_y}_{2\chi} = \sprod{\partial_t u_y}{u_y}_{2\chi} + \sprod{\partial_t v_y}{v_y}_{2\chi}\,. 
	\end{align*}
	From equation~\eqref{kinetic11a} the first term can be written as  
	\begin{align*}
		\sprod{\partial_t u_y}{u_y}_{2\chi} &= \sprod{- u_{yy} + \xdot u_{yy} -2\xdot\chi \sy u_y - 4\xdot \delta_0 (u+1) + 2\chi\delta_0(u-v)}{u_y}_{2\chi}\\
 &\quad -\sprod{(1-\chi\sy)u_y}{u_y}_{2\chi} +\sprod{(1-\chi\sy)v_y}{u_y}_{2\chi} 
  \,.
	\end{align*}
	Computing explicitly the scalar product $\sprod{\cdot}{\cdot}_{2\chi}$ we get
	\begin{align*}
		\sprod{\partial_t u_y}{u_y}_{2\chi} &= 2 \chi (u(0;t) - v(0;t))\,\partial_y u(0;t) - \|u_y\|^2 - 4 \xdot \chi \partial_y u(0;t) (u(0;t) + 1)\\
  &\quad- \xdot \chi \sprod{\sy u_y}{u_y}_{2\chi}
  - \chi \sprod{\sy v_y}{u_y}_{2\chi} + \sprod{u_y}{v_y}_{2\chi}\,.
	\end{align*}
	We do the same for $\sprod{\partial_t v_y} {v_y}_{2\chi}$ and then we sum the two results,
	\begin{align*}
		\sprod{\partial_t u_y}{u_y}_{2\chi} + \sprod{\partial_t v_y} {v_y}_{2\chi}
		&= 2 \chi (u(0;t) - v(0;t)) (u_y(0;t) + v_y(0;t)) - \int |u_y-v_y|^2\eta_\infty\dy\\
		&\quad - \frac{\xdot}{2} \int \partial_y (|u_y|^2 + |v_y|^2)\eta_\infty\dy\\
		&\quad - 4\xdot \chi (u(0;t) + 1) \partial_y u(0;t) - 4\xdot \chi (v(0;t) + 1) \partial_y v(0;t)
	\end{align*}
 This argument is only formal as $u_y,v_y$ may not be continuous at zero. The fully rigorous expression is instead
	\begin{align*}
		\sprod{\partial_t u_y}{u_y}_{2\chi} + \sprod{\partial_t v_y} {v_y}_{2\chi}
		&= - \int |u_y-v_y|^2\eta_\infty\dy - \frac{\xdot}{2} \int \partial_y (|u_y|^2 + |v_y|^2)\eta_\infty\dy \\
		&\quad + 2 \chi (u(0;t) - v(0;t)) \llangle u_y + v_y \rrangle \\
		&\quad - 4 \chi \xdot ( u(0;t) \llangle u_y \rrangle + v(0;t)  \llangle v_y \rrangle ) - 4 \chi \xdot \llangle u_y + v_y \rrangle\,.
	\end{align*}
\end{proof}

\section{Poincar\'e and interpolation inequalities}\label{P inequality}\label{Sect: 4}

In order to estimate $\xdot$, we make use of the following inequalities.
\begin{lem}[Interpolation inequalities]\label{lem:interpol}
	Fix $a\ge b>0$. For any function $f\in L^1_a(\R)$ such that $f'\in L^2_b(\R)$, we have
	\begin{equation}\label{interpol}
		\left | f(0) -  \langle f \rangle_a \right|^2 \leq  \left(\frac1{2a-b}\right) \frac12 \int |f'(y)|^2e^{-b|y|}\,\dy
	\end{equation}
	and 
	\begin{equation}\label{interpol-average}
		\langle f\rangle_a \le \frac{2}{\sqrt{2(2a-b)}} \left(\int |f(y)|^2e^{-b|y|}\,\dy\right)^{1/2}\,.
	\end{equation}
\end{lem}
\begin{proof}
	The first inequality is simply an interpolation result, and was shown in \cite[Lemma 3.5]{CalHof20}. The second estimate follows in the same way as a direct consequence of H\"older's inequality,
	\begin{align*}
		\langle f\rangle_a &= \frac{a}{2}\int f(y)e^{-(a-b)|y|}e^{-b|y|}\dy\\
		&\le \frac{a}{2}\left(\int |f(y)|^2e^{-b|y|}\dy\right)^{1/2
		}\left(\int e^{-2(a-b)|y|}e^{-b|y|}\dy\right)^{1/2}\\
		&= \frac{2}{\sqrt{2(2a-b)}} \left(\int |f(y)|^2e^{-b|y|}\,\dy\right)^{1/2}\,.
	\end{align*}
\end{proof}

\begin{prop}[Poincar\'e Inequalities]\label{prop:Poincare}
	For any $u,v\in L^1\left(e^{-2\chi|y|}\right)$ such that $u_y,v_y \in L^2\left(e^{-2\chi|y|}\right)$, we have
	\begin{subequations}
		\begin{align}\label{Poincare-1}
			\|\Pi W\|&\le \frac{1}{\chi} \|\Pi W_y\|\,,\\
			\label{Poincare-2}
			\|(\I-\Pi)W\|^2&\le \frac{1}{\chi^2}\|(\I-\Pi)W_y\|^2 +\frac{1}{2\chi}\langle u-v\rangle_{2\chi}^2\,,\\
			\label{Poincare-3}
			\|W\|^2 &\le\frac{1}{\chi^2}\|W_y\|^2 +\frac{1}{2\chi}\langle u-v\rangle_{2\chi}^2\,.
		\end{align}
	\end{subequations}
\end{prop}

\begin{proof}
	These estimates follow from the classical Poincar\'e inequality with exponential weight:
	For any $w\in L^1\left(e^{-2\chi|y|}\right)$ such that $w' \in L^2\left(e^{-2\chi|y|}\right)$,
	\begin{equation}\label{eq:poincare-classic}
		\int_\R |w(y) - \langle w \rangle_{2\chi}|^2 e^{-2\chi |y|}\,\dy \leq \frac 1{\chi^2} \int_\R | w'(y)|^2 e^{-2\chi|y|}\,\dy\, .
	\end{equation}
	Moreover, the constant $\chi^{-2}$ is optimal.
	Choosing $w=u+v$ and using the first conservation law in Corollary~\ref{cor:CSL-pert} immediately gives the first estimate \eqref{Poincare-1}. Similarly, with $w=u-v$, the Poincar\'e inequality \eqref{eq:poincare-classic} translates into estimate \eqref{Poincare-2}. Adding \eqref{Poincare-1} and \eqref{Poincare-2} immediately yields \eqref{Poincare-3}.
\end{proof}

\begin{cor}\label{cor:Poincare-improved-est}
Assume $\alpha\ge 0$.
For any $u,v\in L^1\left(e^{-\lambda|y|}\right)$ such that $u_y,v_y \in L^2\left(e^{-2\chi|y|}\right)$ that solve \eqref{kinetic10}, the following improved Poincar\'e inequality holds
	\begin{align*}
		\|(\I-\Pi)W\|^2&\le \frac{1}{\chi^2}\|(\I-\Pi)W_y\|^2 +\frac{1}{\chi}\langle u-v\rangle_{\lambda}\left(\langle u-v\rangle_{2\chi}-\frac12\langle u-v\rangle_{\lambda} \right)\,.
	\end{align*}
\end{cor}
\begin{proof}
Using the improved Poincar\'e inequality in \cite[Proposition 3.1]{CalHof20} it immediately follows that we can improve \eqref{Poincare-2} to the above estimate.
\end{proof}

\begin{cor}\label{cor: P ineq without extra term}
Assume $\alpha=0$.
For any $u,v\in L^1\left(e^{-2\chi|y|}\right)$ such that $u_y,v_y \in L^2\left(e^{-2\chi|y|}\right)$ that solve \eqref{kinetic10} with initial condition $\langle u - v \rangle_{2\chi}|_{t=0}= 0$, the following classical Poincar\'e inequality holds
	\begin{subequations}\label{eq: P ineq without extra term}
		\begin{align}\label{Poincare-1 a}
			\|\Pi W\|^2&\le \frac{1}{\chi^2} \|\Pi W_y\|^2\,,\\
			\label{Poincare-2 a}
			\|(\I-\Pi)W\|^2&\le \frac{1}{\chi^2}\|(\I-\Pi)W_y\|^2    \,,\\
			\label{Poincare-3 a}
			\|W\|^2 &\le\frac{1}{\chi^2}\|W_y\|^2 \,.
		\end{align}
	\end{subequations}
\end{cor}
\begin{proof}
This follows directly from Proposition \ref{prop:Poincare} and Lemma \ref{lem:cond-initial}. 
\end{proof}

\subsection{\texorpdfstring{$\xdot$ estimate}{Peak evolution estimate}}\label{sec:xdot}

Recall the expression for $\xdot$ as stated in \eqref{xdot pert}; we can simplify it to
\begin{subequations}
	\label{xdot2}
	\begin{align}
		\xdot &= \frac{\lambda\left(u(0) - v(0)\right) - \sqrt{\alpha} \langle u-v\rangle_\lambda}{c(u,v)}\,,\label{xdot2a}\\
		c(u,v)&:=4\chi+\lambda(u(0)+v(0)) - \sqrt{\alpha} \langle u+v\rangle_\lambda \,.\label{xdot2b}
	\end{align}
\end{subequations}
Comparing with the expression for $\xdot$ in the macroscopic setting \cite{CalHof20}, we see that $ \langle u-v\rangle_\lambda$ is a truly kinetic contribution which does not appear in the macroscopic limit.
\begin{prop}\label{prop:xdotbound}
	For parameters $\chi>0$, $\alpha\ge0$ and solutions $(u,v)$ to equation~\eqref{kinetic11} satisfying 
	\begin{equation*}
		\mu:=\frac{4\chi\sqrt{2(\chi+\sqrt{\alpha})}}{\lambda} > \|\Pi W_y\|+2\chi\|\Pi W\| \,,
	\end{equation*}
	$\xdot$ is controlled by
	\begin{align}\label{eq:xdot-bound}
		|\xdot(t)|\le \frac{\|(\I-\Pi)W_y\|+2\chi\|(\I-\Pi)W\|}{\mu-\|\Pi W_y\|-2\chi\|\Pi W\|}\,.
	\end{align}
\end{prop}
\begin{proof}
	Let us start by providing an upper bound for the numerator in the expression for $\xdot$. We have
	\begin{align*}
		\left|\lambda\left(u(0) - v(0)\right) - \sqrt{\alpha} \langle u-v\rangle_\lambda\right|
		&= \left|\lambda\left[u(0) - v(0) -\langle u-v\rangle_\lambda\right] +2\chi \langle u-v\rangle_\lambda\right|\\
		&\le \lambda\left|u(0) - v(0) -\langle u-v\rangle_\lambda\right| +2\chi \left|\langle u-v\rangle_\lambda\right|\,.
	\end{align*}
	From Lemma~\ref{lem:interpol}, we estimate
	\begin{align*}
		\left|u(0) - v(0) -\langle u-v\rangle_\lambda\right|^2 \le \frac{1}{2(\chi+\sqrt{\alpha})} \frac12 \int |u_y-v_y|^2\eta_\infty\dy
		= \frac{1}{2(\chi+\sqrt{\alpha})} \|(\I-\Pi)W_y\|^2\,,
	\end{align*}
	and
	\begin{align*}
		\left|\langle u-v\rangle_\lambda\right|\le \frac{\lambda}{2\sqrt{\chi+\sqrt{\alpha}}} \left(\int |u-v|^2\eta_\infty\dy\right)^{1/2}
		=\frac{\lambda}{\sqrt{2(\chi+\sqrt{\alpha})} } \|(\I-\Pi)W\|\,.
	\end{align*}
	Hence, the numerator in \eqref{xdot2a} is bounded by
	\begin{align*}
		\left|  \lambda\left(u(0) - v(0)\right) - \sqrt{\alpha} \langle u-v\rangle_\lambda\right|
		\le \frac{\lambda}{\sqrt{2(\chi+\sqrt{\alpha})} } \left(\|(\I-\Pi)W_y\|+2\chi \|(\I-\Pi)W\|\right)\,.
	\end{align*}
	Similarly, we can derive a bound from below for the denominator in \eqref{xdot2b} by estimating
	\begin{align*}
		\left|u(0) + v(0) -\langle u+v\rangle_\lambda\right|^2 \le \frac{1}{2(\chi+\sqrt{\alpha})} \frac12 \int |u_y+v_y|^2\eta_\infty\dy
		= \frac{1}{2(\chi+\sqrt{\alpha})} \|\Pi W_y\|^2\,,
	\end{align*}
	and
	\begin{align*}
		\left|\langle u+v\rangle_\lambda\right|\le \frac{\lambda}{2\sqrt{\chi+\sqrt{\alpha}}} \left(\int |u+v|^2\eta_\infty\dy\right)^{1/2}
		=\frac{\lambda}{\sqrt{2(\chi+\sqrt{\alpha})} } \|\Pi W\|\,.
	\end{align*}
	Putting the above estimates together, we obtain
	\begin{align}\label{c-bound}
		c(u,v)&=4\chi+\lambda\left[u(0)+v(0)-\langle u+v\rangle_\lambda\right] +2\chi \langle u+v\rangle_\lambda\notag\\
		&\ge 4\chi-\frac{\lambda}{\sqrt{2(\chi+\sqrt{\alpha})} }\|\Pi W_y\| -\frac{2\chi\lambda}{\sqrt{2(\chi+\sqrt{\alpha})} } \|\Pi W\| \,.
	\end{align}
	Substituting in the expression for $\xdot$ completes the proof.
\end{proof}

\begin{remark}
	Using Proposition~\ref{prop:Poincare}, the bound on $\xdot$ can be further simplified to
	\begin{align*}
		|\xdot(t)|\le \frac{3\|(\I-\Pi)W_y\|+\sqrt{2\chi}|\langle u-v\rangle_{2\chi}|}{\mu-3\|\Pi W_y\|}\,.
	\end{align*}
\end{remark}


\section{Hypocoercivity}\label{Sect: 5}

This section is devoted to bounding the dissipation of a suitable modified entropy which is a well-chosen $\delta$-perturbation of the $H^1$-entropy. In particular, the goal is to write an upper bound on the dissipation of this modified entropy that depends on the $H^1$-distance of the solution from the local and the global equilibria.

In Section~\ref{sec:set-up} we introduce the notions of microscopic and macroscopic coercivity as well as a suitable operator splitting that allows to derive the $H^1$-entropy dissipation in a compact form. Next, we introduce the exact expression for the modified entropy in Section~\ref{sec:modified entropy}. The remainder of Section~\ref{Sect: 5} is then dedicated to control the dissipation of the additional terms in this modified entropy; Section~\ref{sec: control nonlin rem} deals with the linear terms, whereas Section~\ref{sec:remainder} handles the nonlinear contributions. The nonlinear terms are mainly controlled through the regularizing properties of the operator introduced to modify the entropy. This operator morally acts as the inverse of the space derivative, allowing to control the higher order derivatives appearing in the nonlinear operator. 
We then collect all these estimates for hypocoercive control of the modified entropy dissipation in Section~\ref{sec: hypo control}. The bound we obtain has coefficients depending on the nonlocal term $\xdot$, on the entropy modification parameter $\delta$, on $|u(0;t) - v(0;t)|$, and on the parameters of the model.


\subsection{Set-up}\label{sec:set-up}
In order to show exponential convergence to equilibrium, we write our system with a more general notation consistent with the one adopted by Dolbeault, Mouhot, and Schmeiser in \cite{DMS15}. Following the operator splitting in \cite{favre_nonlinhypo}, model \eqref{kinetic11} can be expressed as
\begin{equation}
	\label{hypo eq formulation}
	\partial_t W_y + \T W_y = \L W_y + \RR(W)\,, 
\end{equation}
where the collision operator $\L$ and the transport operator $\T$ are defined as
\begin{equation*}
	\L W := 
	\begin{bmatrix}
		v-u\\u-v
	\end{bmatrix}
	= -2(\I-\Pi)W
	,\qquad 
	\T W := 
	\begin{bmatrix}
		u_y - \chi \s (u-v) \\
		-v_y - \chi \s (u-v)
	\end{bmatrix}\,.
\end{equation*}
Then $\L$ and $\T$ are symmetric and skew-symmetric respectively,
$$
\langle \L W_1,W_2\rangle=\langle W_1,\L W_2\rangle\,,\qquad
\langle \T W_1,W_2\rangle=-\langle W_1,\T W_2\rangle\,.
$$
With $\RR(W)$ we refer to the remainder terms
\begin{align*}
	\RR(W) &:= 2\chi \delta_0 (u-v)  \1
	+ \xdot  \partial_y\bigg(\eta_\infty^{-1}\partial_y\Big(\eta_\infty\big(W + \1\big)\Big)\bigg)\\
	&= 2\chi \delta_0 (u-v)  \1
	+ \xdot  \bigg(W_{yy} -2\chi\sign(y)W_y-4\chi\delta_0 \Big(W + \1\Big)\bigg)\,.
\end{align*}
Note that $\RR$ acts also on $W$ and $W_{yy}$, not just on $W_y$. 
Constant vectors $(1,1)$ are in the intersection of the kernels $\ker\T$ and $\ker\L$.
The operator $\Pi$ as defined in \eqref{def:Pi} represents the orthogonal projection on the set of local equilibria:
$$
\L\Pi=0\,.
$$
We notice that the entropy dissipation is coercive with respect to the distance to the set of local equilibria. It follows that the entropy does not contain the whole information needed to show the convergence to the global equilibrium. In order to fill the information gap we should introduce the $H^1$-equivalent Lyapunov (\rus{Lyapunov}) functional $\CL[W_y]$, which will be defined in Section \ref{sec:modified entropy}. 

In order to apply the general methodology of hypocoercivity in our context, we work on derivatives $W_y$ and verify the corresponding microscopic and macroscopic coercivity assumptions \cite{Mouhot-Neumann,Villani-hypo2009}:
\begin{lem}\label{lem:micro-macro-coercivity}
	Let $u_y+v_y\in L^1_\lambda$ and $u_{yy}+v_{yy}\in L^2_{2\chi}$.
	\begin{enumerate}
		\item Microscopic Coercivity: $\L^*=\L$ and
		\begin{equation}\label{micro}
			\langle \L W_y,W_y\rangle_{2\chi} = -2\|(\I-\Pi)W_y\|^2\,.
		\end{equation}
		\item Macroscopic Coercivity: $\T^*=-\T$ and
		\begin{equation}\label{macro}
			\|\T\Pi W_y\|\ge \chi \|\Pi W_y\|\,.
		\end{equation}
		\item Diffusive macroscopic limit: 
		\begin{equation}
			\label{eq: macro diff limit}
			\Pi \T\Pi=0\,.
		\end{equation}
	\end{enumerate}
\end{lem}
In particular, microscopic coercivity means that the collision operator $\L$ is dissipative in the sense that $\langle \L W_y,W_y\rangle_{2\chi}\le 0$. Macroscopic coercivity means that the transport operator $\T$ is coercive when restricted to $\ker\L$. It corresponds to a spectral gap-like inequality for the operator obtained when taking an appropriate macroscopic diffusion limit \cite{DMS15,calvez_confinement_2015}.  
\begin{proof}[Proof of Lemma \ref{lem:micro-macro-coercivity}]
	Following the strategy in \cite{DMS15}, the first statement is a direct consequence of the identity
	\begin{align*}
		\langle \L W_y, W_y\rangle_{2\chi} 
		&= \langle \L W_y, (\I-\Pi)W_y+\Pi W_y\rangle_{2\chi} 
		= -2\|(\I-\Pi)W_y\|^2 +\langle W_y, \L\Pi W_y\rangle_{2\chi} 
	\end{align*}
	and $\L\Pi=0$,
	whereas the second statement can be rewritten as
	\begin{align*}
		\frac{4}{(2\chi)^2} \frac12\int\left|u_{yy}+v_{yy}\right|^2\eta_\infty\dy
		\ge \frac{1}{2}\int \left|u_{y}+v_{y}\right|^2\eta_\infty\dy\,.
	\end{align*}
	Thanks to the second conservation law $\langle u_y+v_y\rangle_{\lambda}=0$, we can directly apply the improved Poincar\'e inequality as shown in \cite[Proposition 3.1]{CalHof20} to obtain macroscopic coercivity as stated above.
	
	The diffusive macroscopic limit property holds by direct investigation.
\end{proof}

Thanks to the fact that $\T^*=-\T$ implies $\langle \T W_y,W_y\rangle=0$, we have
\begin{equation}\label{eq: diss of L2 entropy}
	\frac12\frac{\dd}{\dt}\|W_y\|^2 = \langle \L W_y,W_y\rangle + \langle \RR(W) , W_y\rangle\,.
\end{equation}
If $\langle \L W_y,W_y\rangle$ was coercive with respect to the norm $\|\cdot\|$, and the remainder term $\langle W_y,\RR(W)\rangle$ could be controlled with a small enough multiple of $\|W_y\|^2$, then exponential decay to zero as $t\to\infty$ would follow. The microscopic coercivity assumption \eqref{micro} states that such a coercivity for $\L$ does not hold fully for $W_y$, but only with respect to a bound on $(\I-\Pi)W_y$ away from the local equilibria. In other words, the restriction of $\L$ to $(\ker\L)^\perp$ is coercive, and as soon as the evolution has driven the solution to the set of local equilibria $\ker\L$, the action of $\L$ does not contribute anymore, and a different approach is needed to show that the solution reaches the global equilibrium. 

\subsubsection*{Control of the nonlinear and nonlocal terms in the entropy dissipation}
From Proposition \ref{prop:H1dot} we know that
\begin{align}
	\label{eq: diss nonlin entropy}
	\sprod{W_y}{\RR(W)} &= 2\chi \left[u(0)-v(0)-2\xdot\right]\llangle u_y+v_y \rrangle \\
	\nonumber
	&\quad + \frac{\xdot}{2} \int \partial_y \left(|u_y|^2+|v_y|^2\right)\eta_\infty\dy
	-4\chi\xdot\left[u(0)\llangle u_y\rrangle+ v(0)\llangle v_y\rrangle \right]\,.
\end{align}
We will control this term separately in the two different cases: a linearized version of \eqref{kinetic10} with $\alpha \ge 0$, and the full nonlinear model with $\alpha = 0$.

\subsection{Modified Entropy}\label{sec:modified entropy}

The key tool for this aim is the \textit{modified entropy} which is given by the following Lyapunov (\rus{Lyapunov}) functional \cite{DMS15,Herau}
\begin{align*}
	\CL[W_y]:=\frac12 \|W_y\|^2+\delta\langle \A W_y,W_y\rangle_{2\chi}\,,\qquad
	\text{ with } \A:=\left(\I+(\T\Pi)^*(\T\Pi)\right)^{-1}(\T\Pi)^*\,,
\end{align*}
for a well-chosen constant $\delta>0$ to be fixed later.
The idea is to define the modification in such a way that $\CL[W_y]$ is equivalent to $\|W_y\|^2$. It is worth noticing that the operator $\A$ is bounded, furthermore we have 
\begin{lem}
	With the assumptions 1 and 3 of Lemma \ref{lem:micro-macro-coercivity}, and given $f\in L^2_{2\chi}$, then the operator $\A$ is bounded by $1/2$ in $L^2_{2\chi}$, furthermore $$\|\A f\| \le \unmezzo \|(\I - \Pi)f\| \,.$$
\end{lem}
\begin{proof}
	For completeness, we recall the proof given in Lemma 1 of \cite{DMS15}.
	
	Let introduce $h = \A f$, using the definition of $\A$ and the properties of $\T$ and $\Pi$ we can write
	\begin{equation*}
		h = (\T \Pi)^* f - (\T\Pi)^*(\T \Pi) h\,.
	\end{equation*}
	Testing by $h$ the latter equation we immediately get
	\begin{align*}
		\|h\|^2 + \| \T \Pi h\|^2 &= \sprod{f}{\T\Pi h}  = \sprod{(\I - \Pi)f}{\T\Pi h}\\
		&\le \frac{1}{4} \|(\I - \Pi)f\|^2 + \|\T\Pi h\|^2\,,
	\end{align*}
	which concludes the proof.
\end{proof}

\begin{lem}
	Given $W$ solution to equation \eqref{hypo eq formulation}, then the following chain of inequalities holds
	\begin{equation}\label{eq: equivalence entropy L2-norm}
		\frac{1-\delta}{2} \|W_y\|^2 \le \CL[W_y] \le \frac{1+\delta}{2} \|W_y\|^2\,.
	\end{equation}
\end{lem}
\noindent
A general proof has been given in Lemma 1 of \cite{DMS15}, it concerns the boundedness of the $\A$ operator.

If the Lyapunov (\rus{Lyapunov}) functional is coercive with respect to the $L^2$-norm of $W_y$ then it is possible to obtain Gr\"onwall's inequality
$$
\frac{\dd}{\dt} \CL[W_y](t) \le - \lambda' \|W_y\|^2 \le - \lambda \CL[W_y]
$$
for some computable convergence rate $\lambda$. Such a strategy is referred to as \emph{hypocoercivity} \cite{DMS15,Herau,Villani-hypo2009}. We start with a preliminary result on the dissipation of the perturbation term in the modified entropy.
\begin{lem}\label{lem:delta-dissipation}
For $W$ solving \eqref{hypo eq formulation}, it holds
    \begin{equation*}
      \frac{\dd}{\dt}   \langle \A W_y,W_y\rangle_{2\chi}
      = - D(t) +R(t)\,,
    \end{equation*}
    where $D(t)$ and $R(t)$ are given by
\begin{subequations}\label{eq:remainders}
   \begin{align}
	D(t)&= \langle \A\T\Pi W_y,W_y\rangle_{2\chi}
	+\langle \A\T(\I-\Pi)W_y,W_y\rangle_{2\chi}
	-\langle \T\A W_y, W_y\rangle_{2\chi}
	-\langle \A\L W_y, W_y\rangle_{2\chi}\,,
 \label{eq:remainderD}\\
	R(t)&=\langle \A\RR(W),W_y\rangle_{2\chi} + \langle\A W_y,\RR(W)\rangle_{2\chi}\,.
  \label{eq:remainderR}
\end{align} 
\end{subequations}
\end{lem}
\begin{proof}
This result follows by direct computation, noting that $\T^*=-\T$ and $\L^*\A=\L\A=\L\Pi\A=0$. 
\end{proof}
Here, we used the splitting as proposed in \cite{favre_nonlinhypo}.
We can think of $D(t)$ and $R(t)$ as a linear and nonlinear contributions, respectively. Indeed, the expression for $D(t)$ in \eqref{eq:remainderD} only contains linear operators acting on $W_y$. And we will see that the linear contributions in the first term $\langle \A\RR(W),W_y\rangle_{2\chi}$ in \eqref{eq:remainderR} disappear, and the linear contributions in the second part $\langle\A W_y,\RR(W)\rangle_{2\chi}$ can all be expressed in terms of $\|(\I -\Pi)W_y\|^2$ which will be absorbed into the decay term $-\|(\I -\Pi)W_y\|^2$ in the expression for $\frac{\dd}{\dt} \CL[W_y](t)$ above.

\subsection{\texorpdfstring{Control of the linear remainder $D(t)$}{Control of the linear remainder D(t)}}\label{sec: control nonlin rem}
From the macroscopic coercivity \eqref{macro} and writing $\mathcal{L}$ as $(\T\Pi)^*(\T\Pi)$ we obtain that
\begin{align*}
	\langle \mathcal{L} W_y,W_y\rangle_{2\chi} \ge \chi\, \|\Pi W_y\|^2 \,.
\end{align*}
Writing $\A\T\Pi = \left(\I+\mathcal{L}\right)^{-1}\mathcal{L}$ we can consider the spectral gap of $\frac{\calL}{1+\calL}$ (see \cite{Favre_BGK+T} Lemma 5). Then we notice that $\A=\Pi\A$, concluding that 
\begin{equation}\label{eq: ATP}
	\langle \A\T\Pi W_y,W_y\rangle_{2\chi} = \langle \frac{\calL}{\I+\calL} W_y,\Pi W_y\rangle_{2\chi} \ge \frac{\chi}{1+\chi}  \|\Pi W_y\|^2\,.
\end{equation}

We observe that $\calL=-\Pi \T^2 \Pi$ is the weighted Laplacian operator in the $L^2_{2\chi}$-space, i.e.
\begin{equation*}
	\calL z = - \eta_\infty^{-1} \partial_y\big( \eta_\infty \partial_y(\Pi z) \big) =  - \Delta_\eta(\Pi z)\,,\qquad \text{ where}\,\quad
	\Delta_\eta u: =  \eta_\infty^{-1} \partial_y\big( \eta_\infty \partial_y(u) \big)\,.
\end{equation*}
Here, the operator $\Delta_\eta$ acts component-wise, $\Delta_\eta(z_1,z_2)=(\Delta_\eta z_1, \Delta_\eta z_2)$, and $\calL^*=\calL$ since
\begin{equation}
	\sprod{\calL z}{\xi}_{2\chi} 
	=\frac12\int\partial_y(z_1+z_2)\partial_y(\xi_1+\xi_2)\eta_\infty\dy
	= \sprod{z}{\calL \xi}_{2\chi}  \quad \text{ for all } z, \xi \in H^1_{2\chi}\,.
 \label{eq:L*=L}
\end{equation}

\begin{lem}\label{lem:hypo-lin-bounds}
	For all $W\in H^1_{2\chi}$,
	\begin{gather}
		\label{eq: A}
		\|\A W_y\| \le \frac12 \|(\I-\Pi)W_y\|\,,\\
		\label{eq: TA}
		\|\T\A W_y\| \le \|(\I-\Pi)W_y\|\,,\\
		\label{eq: AT(I-P)}
		\A\T(\I-\Pi) = 0\,,\\
		\label{eq: AL}
		\|\A \L W_y\|  \le \|(\I-\Pi)W_y\|\,.
	\end{gather}
\end{lem}
\begin{proof}
	The first two inequalities have been shown in \cite[Lemma 1]{DMS15}. The third identity is a consequence of
	\begin{align*}
		\T(\I-\Pi)W_y=\frac12\Delta_\eta(u-v)\1 \,,
	\end{align*}
	and the fact that for any scalar function $z$, we have $\Pi\T\begin{bmatrix} z\\z \end{bmatrix}
	=\Pi \begin{bmatrix} z'\\-z'\end{bmatrix}=0$.
	Finally, for the last inequality, we recall that 
	\begin{equation*}
		\|\L W_y\|^2 = 2 \|(\I - \Pi)W_y\|^2
	\end{equation*}
	and using the boundedness of $\A$ we conclude the proof.
\end{proof}

\begin{cor}\label{cor:hypo-lin-bounds}
The linear remainder term can be controlled by
	\begin{align*}
		-D(t) \le &- \frac{\chi}{1+\chi} \|\Pi W_y\|^2 + \|(\I - \Pi)W_y\|^2 + \|(\I-\Pi)W_y\| \|\Pi W_y\|  \,.
	\end{align*}
\end{cor}
\begin{proof}
We can control the linear term $D(t)$ by applying the estimates in Lemma~\ref{lem:hypo-lin-bounds} to the expression for $D(t)$ in Lemma~\ref{lem:delta-dissipation} and using the properties of the operators $\L, \T, \Pi$ shown before. Indeed, the first term can be bounded as in \eqref{eq: ATP}. The second term vanishes using \eqref{eq: AT(I-P)}. Since $\A = \Pi\A$, $\Pi\T\Pi=0$ from Lemma~\ref{lem:micro-macro-coercivity} and $(\I-\Pi)^*=\I-\Pi$, the third term can be bounded as follows:
\begin{align*}
\langle \T\A W_y, W_y\rangle_{2\chi}
&= \langle \T\Pi\A W_y, W_y\rangle_{2\chi}
= \langle (\I-\Pi)\T\Pi\A W_y, W_y\rangle_{2\chi}
=\langle \T\Pi\A W_y, (\I-\Pi)W_y\rangle_{2\chi}\\
&\le \|\T\Pi\A W_y\|\,\|(\I-\Pi)W_y\|
=\|\T\A W_y\|\,\|(\I-\Pi)W_y\|
\le \|(\I-\Pi)W_y\|^2\,,
\end{align*}
where we used \eqref{eq: TA} in the last inequality. For the last term in \eqref{eq:remainderD}, we use $\A = \Pi\A$ together with \eqref{eq: AL} to conclude
\begin{align*}
 - \langle \A \L W_y, W_y\rangle_{2\chi}   
 &=  - \langle \Pi\A \L W_y, W_y\rangle_{2\chi}  
  =  - \langle \A \L W_y, \Pi W_y\rangle_{2\chi}  \\
  &\le \|\A \L W_y\|\,\|\Pi W_y\|
  \le \|(\I-\Pi) W_y\|\,\|\Pi W_y\|\,.
\end{align*}
This concludes the proof.
\end{proof}

\subsection{\texorpdfstring{Control of the nonlinear remainder $R(t)$}{Control of the nonlinear remainder R(t)}}\label{sec:remainder}

There are no known results in the literature that can directly be applied to control the nonlinear remainder $R(t)$, it depends on the specific structure of the $\RR(W)$ term. This subsection is devoted to bound the nonlinear remainder in terms of the norm of $W_y$. The way to get the wished inequalities is to use recursively the regularizing property of the operator $\calL$, appearing in $\A = -\big( \I + \mathcal{L}\big)^{-1} \Pi \T $. The operator $\calL$ acts as a Laplacian. Roughly speaking, the operator $\A$ can be seen as acting as the inverse of the space derivative. Morally, we gain one derivative order each time that the operator $\A$ is applied.

\begin{lem}\label{lem:hypo-nonlin-bounds-1}
	If $W$ solves \eqref{hypo eq formulation}, then
	\begin{gather}
		\label{eq: AR}
		\sprod{\A \RR(W)}{W_y}_{2\chi} \le |\xdot|\,\|W_y\|\, \Big( \|(\I - \Pi)W_y\| + 4 \chi^2 \|(\I - \Pi)W\| \Big)\,.
	\end{gather}
\end{lem}
\begin{proof}
	First we write $\A$ explicitly, then we use the selfadjointness of $(\I + \calL)^{-1}$ to get
	\begin{equation*}
		\sprod{\A \RR(W)}{W_y}_{2\chi} = - \sprod{\Pi \T \RR(W)}{(\I + \calL)^{-1} W_y}_{2\chi} \,.
	\end{equation*}
	According to the definition of $\T$ and $\Pi$ we know that 
	\begin{equation*}
		\Pi \T \RR(W) = 2 \chi \,\Pi \T \big[\delta_0 (u-v)\big] \1
		+ \xdot \, \Pi\T\Bigg[ \partial_y \bigg( \eta_\infty^{-1} \partial_y \Big( \eta_\infty \big( W + \1\big)
		\Big) \bigg) \Bigg]\,.
	\end{equation*}
	The first and the last terms belong to the set of local equilibria, this implies that they vanish as a direct consequence of the property $\Pi \T \Pi = 0$. This is why we can think of the first contribution to $R(t)$ as a truly nonlinear term. Defining the operator $\Q$ as
 \begin{equation}\label{def:Q}
      \Q[W] = \begin{bmatrix}
		q_1 \\ q_2 \end{bmatrix} = \partial_y \big( \eta_\infty^{-1} \partial_y (\eta_\infty W)\big)\,,
 \end{equation}
  the structure of the previous term reduces to
	\begin{equation*}
		\Pi \T \RR(W) = \frac{\xdot}{2} \eta_\infty^{-1} \partial_y \big( \eta_\infty (q_1 - q_2) \big)\1 \,.
	\end{equation*}
	Given $\zeta = (\zeta_1,\zeta_2)$ solution of 
	\begin{equation}
		\label{eq zeta}
		\zeta + \calL \zeta = \Pi\T\RR(W), \quad \text{equivalent to}\quad (\I  + \mathcal{L})^{-1} \Pi\T\RR(W) = \zeta \,,
	\end{equation}
	we notice that $\zeta_1 = \zeta_2$ as $\calL$ acts component-wise. Let define $\xi = (\xi_1, \xi_2)$ as the auxiliary function that solves 
	\begin{equation}
		\label{eq xi_zeta}
		\xi + \calL \xi = \zeta , \quad \text{equivalent to}\quad (\I  + \mathcal{L})^{-1} \zeta = \xi \,,
	\end{equation}
	it follows by construction that also for $\xi$ it holds $\xi_1 = \xi_2$. From the latter property we deduce that we can consider the scalar problem instead. This follows from the fact that also $\calL$ acts component-wise on the first and second variable. For this reason, to lighten the notation, without loss of generality we reduce the computation to the scalar version.

	\subsubsection*{Boundedness of $\calL \xi$} Let us test \eqref{eq xi_zeta} with $\calL \xi$, we get
	\begin{equation*}
		\sprod{\xi}{\calL \xi}_{2\chi} + \|\calL \xi\|^2 = \sprod{\zeta}{\calL \xi}_{2\chi}\,.
	\end{equation*}
	Then, thanks to the positiveness of the first term on the \lhs, and using the Cauchy-Schwarz inequality we obtain the wished estimate
	\begin{equation*}
		\|\calL \xi\| \le \|\zeta\| \,.
	\end{equation*}
	
	\subsubsection*{Estimate on $\zeta$}
	Testing the \lhs~of equation \eqref{eq zeta} by $\zeta$, recalling the property \eqref{eq:L*=L} of the weighted Laplacian, and using the definition of $\xi$ we get
	\begin{equation}
		\label{eq: zeta 1}
		\sprod{\zeta}{\xi}_{2\chi} + \sprod{\calL \zeta}{\xi}_{2\chi} = \sprod{\zeta}{(\I + \calL)\xi}_{2\chi} = \|\zeta\|^2 \,.
	\end{equation}
	
	In order to estimate the \rhs of equation \eqref{eq zeta} we use the definition of $\Q$,
	\begin{align*}
		\sprod{\Pi\T\RR(W)}{\xi}_{2\chi} &= 
		\frac{\xdot}{2} \ints{ \partial_y \big( \eta_\infty (q_1-q_2) \big) \xi \cdot \1  }
		= - \frac{\xdot}{2} \ints{ \eta_\infty (q_1-q_2)  \partial_y \xi \cdot \1 }\\
		&= - \frac{\xdot}{2} \ints{ \eta_\infty \partial_y \Big(\eta_\infty^{-1} \partial_y \big(\eta_\infty (u-v)\big)  \Big)  \partial_y \xi \cdot \1 }\\
		&= \frac{\xdot}{2} \ints{\partial_y (\eta_\infty \partial_y \xi) \eta_\infty^{-1} \partial_y\big(\eta_\infty(u-v)\big)  \cdot \1 }\\
		&= - \frac{\xdot}{2} \ints{ \partial_y\big(\eta_\infty(u-v)\big) \calL\xi \cdot \1  }
		\le \frac{|\xdot|}{2} \|\calL \xi\|\, \|\eta_\infty^{-1} \partial_y\big(\eta_\infty(u-v)\big)\1 \| \\
		&\le \frac{|\xdot|}{2} \|\zeta\|\, \left\|\eta_\infty^{-1} \partial_y\big(\eta_\infty(u-v)\big) \1 \right\| \,.
	\end{align*}
	It follows from \eqref{eq: zeta 1} and \eqref{eq zeta} that $\|\zeta\|^2 =\sprod{\Pi\T\RR(W)}{\xi}_{2\chi}$, and so
	\begin{equation}
		\label{estimate on zeta}
		\|\zeta\| \le \frac{|\xdot|}{2}  \left\|\eta_\infty^{-1} \partial_y\big(\eta_\infty(u-v)\big) \1 \right\| \,.
	\end{equation}
	
	\subsubsection*{Estimate on the \rhs of \eqref{estimate on zeta}}
	Using the definition of the operator $\Pi$ we have that
	\begin{align*}
		 \left\|\eta_\infty^{-1} \partial_y\big(\eta_\infty(u-v)\big) \1 \right\|^2 &= 2\ints{\eta_\infty \Big(\eta_\infty^{-1} \partial_y \big( \eta_\infty (u-v) \big) \Big)^2}\\
		&= 2\ints{\eta_\infty \big( u_y - v_y - 2\chi \s (u-v) \big)^2}\\
		&\le 2 \ints{\eta_\infty |u_y - v_y|^2} + 8 \chi^2 \ints{\eta_\infty |u-v|^2}\\
		&= 2 \|(\I - \Pi)W_y\|^2 + 8 \chi^2 \|(\I - \Pi)W\|^2 \,.
	\end{align*}
	
	Putting together the previous estimates we conclude with the wished result
	\begin{align*}
		\sprod{\A\RR(W)}{W_y}_{2\chi} &=- \sprod{\zeta}{W_y}_{2\chi}
		\le \|\zeta\| \|W_y\| \\
		&\le |\xdot| \|W_y\| \Big( \|(\I - \Pi)W_y\| + 4\chi^2 \|(\I-\Pi)W\| \Big) \,.
	\end{align*}
\end{proof}

\begin{lem}\label{lem:hypo-nonlin-bounds-2}
	Let $W$ be a solution of \eqref{hypo eq formulation}. Then 
	\begin{align}\label{eq: RW AW}
		\sprod{\RR(W)}{\A W_y}_{2\chi} &\le  \frac{c_0}{2} |u(0;t) - v(0;t) | \|(\I-\Pi)W_y\|\notag \\ &\quad + |\xdot|\left( c_1 \|(\I-\Pi)W_y\| + c_2 \|W_y\|^2 +  c_3 \|\Pi W\|^2\right)
	\end{align}
	with $c_1 = 4 \sqrt{\chi}\, (1+2\chi)$, $c_2 = (1+4\chi^2)$, $c_3 = 4 \chi^2$ and 
 $$
 c_0:=
\begin{cases}
4\sqrt{\chi} & \text{ if } \chi^2\le 1/2\,,\\
8\chi^2\sqrt{\chi}/\sqrt{(2\chi-1)(2\chi+1)} & \text{ if } \chi^2 > 1/2\,.
\end{cases}
 $$
\end{lem}

\begin{proof}
	Similarly as for the previous lemma, we define $\zeta$ as the (now scalar) solution of
	\begin{equation}
		\label{eq zetaa}
		\zeta + \calL \zeta = \wlap (u-v)\,, \quad \text{equivalent to}\quad \zeta = (1 + \mathcal{L})^{-1}\wlap(u-v)\,,
	\end{equation}
	then $$\A W_y = - \frac{\zeta}{2} \1\,.$$ 
	We split $\RR(W)$ in three terms $\RR(W) = \RR_1 + \RR_2 +\RR_3$, where
	\begin{equation}\label{R-splitting}
		\RR_1(W):= 2\chi \delta_0 (u-v) \1\,,\quad
     \RR_2(W):=\xdot \partial_y \bigg( \etainf^{-1} \partial_y \big( \etainf W \big) \bigg)\,,\quad
     \RR_3(W):=\xdot  \partial_y \bigg( \etainf^{-1} \partial_y \big( \etainf \big) \bigg)\1\,.
	\end{equation}

 \subsubsection*{Preliminary results for $\zeta$}
 The weighted interpolation introduced in \cite{CalHof20} and defined in \eqref{interpol} gives the following estimate
	\begin{equation*}
		|\zeta(0) - \langle \zeta \rangle_{2\chi}|^2 \le \frac{1}{4\chi}\, \|\zeta_y\|^2 \,.
	\end{equation*}
	At this point we look for a bound on $\zeta_y$ testing \eqref{eq zetaa} with $\zeta$,
	\begin{align}
		\nonumber
		\|\zeta\|^2 + \|\zeta_y\|^2 
  &=\int \zeta ((1+\calL)\zeta) \etainf
  = \int \zeta \Delta_\eta(u-v) \etainf \notag\\
  &\le \bigg( \int{|\zeta_y|^2 \etainf} \bigg)^{\frac{1}{2}}  \bigg( \int{|u_y-v_y|^2 \etainf} \bigg)^{\frac{1}{2}}
		= \|\zeta_y\| \sqrt{2}\|(\I - \Pi)W_y\|\\
		\label{eq zeta zeta' estim}
		&\le \frac{\epsi}{2} \|\zeta_y\|^2 + \frac{1}{\epsi}\|(\I - \Pi)W_y\|^2 \,,
	\end{align}
	which implies
	\begin{equation*}
		\|\zeta_y\|^2 \le \frac{c_\epsi}{\epsi} \|(\I - \Pi)W_y\|^2 - c_\epsi \|\zeta\|^2 
	\end{equation*}
	with $c_\epsi = 2/(2-\epsi)$ and $\epsi\in(0,2)$, i.e.~$c_\epsi \in (1,\infty)$. As an immediate consequence, we have
 \begin{align}
     \|\zeta\|^2 + \|\zeta_y\|^2 
     &\le\frac{\epsi}{2} \left( \frac{c_\epsi}{\epsi} \|(\I - \Pi)W_y\|^2 - c_\epsi \|\zeta\|^2 \right)+ \frac{1}{\epsi}\|(\I - \Pi)W_y\|^2\notag\\
     &\le\left(\frac{c_\epsi}{2} +\frac{1}{\epsi}\right)\|(\I - \Pi)W_y\|^2
     = \frac{2}{(2-\epsi)\epsi}\|(\I - \Pi)W_y\|^2\label{eq zeta sum}
 \end{align}
 which is optimized by choosing $\epsi=1$.
 Next, in order to control $\zeta(0)$, we use Jensen's inequality to deduce $ \langle \zeta \rangle_{2\chi}^2 \le \chi \|\zeta\|^2$, and by combining with the above bounds we estimate
 \begin{align*}
|\zeta(0)|^2 &\le \left(|\zeta(0) - \langle \zeta \rangle_{2\chi}| + |\langle \zeta \rangle_{2\chi}|\right)^2
\le  2|\zeta(0) - \langle \zeta \rangle_{2\chi}|^2 + 2 \langle \zeta \rangle_{2\chi}^2\\
&\le \frac{1}{2\chi} \|\zeta_y\|^2 + 2\chi  \|\zeta\|^2
\le \frac{c_\epsi}{2\chi\epsi} \|(\I - \Pi)W_y\|^2 + \left(2\chi - \frac{c_\epsi}{2\chi}\right) \|\zeta\|^2 \,.
 \end{align*}
We choose $\epsi=1$ if $\chi^2\le 1/2$ and $\epsi=2-1/(2\chi^2)$ if $\chi^2>1/2$. With this choice, the last term in the inequality above is less or equal to zero. We conclude that
\begin{equation}\label{eq bdd zeta0}
|\zeta(0)|^2 
\le \left(\frac{c_0}{4\chi}\right)^2 \|(\I - \Pi)W_y\|^2\,,
\end{equation}
where
\begin{equation*}
    c_0:=4\chi\left(\frac{c_\epsi}{2\chi\epsi}\right)^{1/2}=
\begin{cases}
4\sqrt{\chi} & \text{ if } \chi^2\le 1/2\,,\\
8\chi^2\sqrt{\chi}/\sqrt{(2\chi-1)(2\chi+1)} & \text{ if } \chi^2 > 1/2\,.
\end{cases}
\end{equation*}

	\subsubsection*{Bound on $\sprod{\RR_1}{\zeta \1}_{2\chi}$}
	Using the estimate \eqref{eq bdd zeta0} from the preceding paragraph, we immediately obtain the bound
	\begin{align}
		\sprod{\RR_1}{\zeta\1}_{2\chi} &= 4 \chi \zeta(0;t) (u(0;t) - v(0;t))\notag\\
		&\le c_0 |u(0;t) - v(0;t)| \, \|(\I - \Pi)W_y\|\,.\label{eq bdd first nonlin term}
	\end{align}

\subsubsection*{Bound on $\sprod{\RR_2}{\zeta\1}$}
 	Integrating by parts $\sprod{\RR_2}{\zeta\1}_{2\chi}$ we get
	\begin{align*}
 \sprod{\RR_2}{\zeta\1}_{2\chi}&=
		-\xdot \ints{\partial_y(\zeta \etainf) \partial_y\big(\etainf(u+v)\big)\etainf^{-1}} \\
  &\le |\xdot| \bigg[ \ints{|\partial_y (\zeta \etainf)|^2 \etainf^{-1}} \bigg]^{\unmezzo} \bigg[ \ints{|\partial_y \big((u+v) \etainf\big)|^2 \etainf^{-1}} \bigg]^{\unmezzo} \\
		&= 2|\xdot| \bigg[ \ints{\big(|\zeta_y|^2 + 4\chi^2 |\zeta|^2 \big) \etainf} \bigg]^{\unmezzo} \bigg[ \ints{\big(|u_y+v_y|^2 + 4\chi^2 |u+v|^2\big) \etainf} \bigg]^{\unmezzo}\\
		&\le |\xdot| \Big(\|\zeta_y\|^2 + 4 \chi^2 \|\zeta\|^2 + \|u_y+v_y\|^2 + 4 \chi^2 \|u+v\|^2\Big)\\
		&= |\xdot| \Big(\|\zeta_y\|^2 + 4 \chi^2 \|\zeta\|^2 + 2\|\Pi W_y\|^2 + 8 \chi^2 \|\Pi W\|^2\Big)\,.
	\end{align*}
	
	\subsubsection*{Bound on $\sprod{\RR_3}{\zeta\1}$}
	Similarly, integrating by parts the last term yields
	\begin{align*}
		\sprod{\RR_3}{\zeta\I}_{2\chi} &= - 2 \xdot \ints{ \partial_y(\zeta \etainf )\bigg( \etainf^{-1} \partial_y \big( \etainf \big) \bigg)}\\
		&=  4 \chi \xdot \ints{ \partial_y(\zeta \etainf ) \s } \\
		&\le 4 \chi |\xdot| \ints{ \big(|\zeta_y| + 2 \chi \,|\zeta|\big) \etainf }\\
		&\le 4 |\xdot| \sqrt{\chi}\, (\|\zeta_y\| + 2\chi \|\zeta\|)\,.
	\end{align*}

	\subsubsection*{Combining estimates}
 It follows from \eqref{eq zeta sum} that for $\epsi=1$ we have 
 $$\|\zeta\| + \|\zeta_y\| \le 2 \|(\I - \Pi)W_y\|\,.$$ Combining all the previous bounds, we then obtain
	\begin{align*}
		\sprod{\RR(W)}{\zeta\1}_{2\chi} &\le c_0|u(0;t) - v(0;t) |\|(\I-\Pi)W_y\| + 4 |\xdot| \sqrt{\chi}\, (\|\zeta_y\| + 2\chi \|\zeta\|) \\ &\quad +|\xdot| \Big(\|\zeta_y\|^2 + 4 \chi^2 \|\zeta\|^2 + 2\|\Pi W_y\|^2 + 8 \chi^2 \|\Pi W\|^2\Big)\\
		&\le c_0|u(0;t) - v(0;t) |\|(\I-\Pi)W_y\| + 8 |\xdot| \sqrt{\chi}\, (1+2\chi) \|(\I-\Pi)W_y\| \\
		&\quad + 2|\xdot| \bigg( (1+4\chi^2) \|(\I-\Pi)W_y\|^2 + \|\Pi W_y\|^2 + 4 \chi^2 \|\Pi W\|^2\bigg) \\
		&\le c_0|u(0;t) - v(0;t) |\|(\I-\Pi)W_y\| \\
		&\quad + 2|\xdot| \bigg(4 \sqrt{\chi}\, (1+2\chi) \|(\I-\Pi)W_y\| + (1+4\chi^2) \|W_y\|^2 + 4 \chi^2 \|\Pi W\|^2\bigg)\,.
	\end{align*}
 This concludes the proof, recalling $\A W_y =  - \frac{\zeta}{2} \1$.
\end{proof}

The above estimates allow us finally to obtain  control of the nonlinear remainder term.
\begin{cor}\label{cor:hypo-nonlin-bounds}
If $W$ solves \eqref{hypo eq formulation}, then there exists constants $c_0,c_1,c_2,c_3>0$ only depending on $\chi$ such that
    \begin{align*}
        R(t)
        &\le
         |\xdot|\,\|W_y\|\, \Big( \|(\I - \Pi)W_y\| + 4 \chi^2 \|(\I - \Pi)W\| \Big)\\
         &\quad+\frac{c_0}{2} |u(0;t) - v(0;t) | \|(\I-\Pi)W_y\|\notag \\ &\quad + |\xdot|\left( c_1 \|(\I-\Pi)W_y\| + c_2 \|W_y\|^2 +  c_3 \|\Pi W\|^2\right)\,.
    \end{align*}
\end{cor}
\begin{proof}
    This bound follows directly by combining the results in Lemma~\ref{lem:hypo-nonlin-bounds-1} and Lemma~\ref{lem:hypo-nonlin-bounds-2}.
\end{proof}

\subsection{Hypocoercive control}\label{sec: hypo control}
At this point we have all the estimates needed to control the time derivative of the Lyapunov (\rus{Lyapunov}) functional $\CL[W](t)$ in terms of norms of $W_y$, $\Pi W_y$, $(\I - \Pi)W_y$.

\begin{prop}\label{prop: gronwall_hypo_general_nonlocal}
Let $W$ be a solution of equation \eqref{hypo eq formulation} satisfying $\langle u - v \rangle_{2\chi}|_{t=0}= 0$. Then there exists constants $c_0, c_1, c_2, c_3>0$ only depending on $\chi$ such that for any $\delta>0$,
	\begin{align*}
	\frac{\dd}{\dt} \CL[W_y](t)
 &\le- 2 \|(\I -\Pi)W_y\|^2  - \delta\frac{\chi}{1+\chi} \|\Pi W_y\|^2
		\\
		&\quad 
		+ |\xdot| \chi \|W_y\|^2 
  + \delta\|(\I - \Pi)W_y\|^2 + \delta\|(\I-\Pi)W_y\| \|\Pi W_y\|\\
  &\quad +\delta |\xdot| (1+4\chi)\,\|W_y\|\, \|(\I - \Pi)W_y\| \\
    &\quad + \delta |\xdot|\left( c_1 \|(\I-\Pi)W_y\| + c_2 \|W_y\|^2 +  c_3 \|\Pi W\|^2\right)\\
  &\quad+ 2\chi \left|u(0)-v(0)-2\xdot\right|\,|\llangle u_y+v_y \rrangle |
       +\delta \frac{c_0}{2} |u(0;t) - v(0;t) | \|(\I-\Pi)W_y\|\\
       &\quad+|\xdot|4\chi\left|u(0)\llangle u_y\rrangle+ v(0)\llangle v_y\rrangle \right|\,.
\end{align*}
\end{prop}
\begin{proof}
Applying the estimates from Corollary~\ref{cor:hypo-lin-bounds} and Corollary~\ref{cor:hypo-nonlin-bounds} to the terms in
Lemma~\ref{lem:delta-dissipation} and combining with  Proposition~\ref{prop:H1dot}, we obtain
\begin{align*}
	\frac{\dd}{\dt} \CL[W_y](t)&= \frac12\frac{\dd}{\dt} \|W_y\|^2  -\delta D(t) +\delta R(t)\\
 &\le- 2 \|(\I -\Pi)W_y\|^2
		+ 2\chi \left[u(0)-v(0)-2\xdot\right]\llangle u_y+v_y \rrangle 
		\\
		&\quad 
		- \frac{\xdot}{2} \int \partial_y \left(|u_y|^2+|v_y|^2\right)\eta_\infty\dy
		-4\chi\xdot\left[u(0)\llangle u_y\rrangle+ v(0)\llangle v_y\rrangle \right]\\
  &\quad
  - \delta\frac{\chi}{1+\chi} \|\Pi W_y\|^2 + \delta\|(\I - \Pi)W_y\|^2 + \delta\|(\I-\Pi)W_y\| \|\Pi W_y\|\\
  &\quad +\delta |\xdot|\,\|W_y\|\, \Big( \|(\I - \Pi)W_y\| + 4 \chi^2 \|(\I - \Pi)W\| \Big)\\
         &\quad+\delta \frac{c_0}{2} |u(0;t) - v(0;t) | \|(\I-\Pi)W_y\|\notag \\ &\quad + \delta |\xdot|\left( c_1 \|(\I-\Pi)W_y\| + c_2 \|W_y\|^2 +  c_3 \|\Pi W\|^2\right)\,.
\end{align*}
Next, we use the Poincar\'e inequality \eqref{Poincare-2 a} together with the assumption $\langle u - v \rangle_{2\chi}|_{t=0}= 0$ as well as integration by parts. This allows to estimate the energy decay by
\begin{align*}
	\frac{\dd}{\dt} \CL[W_y](t)
 &\le- 2 \|(\I -\Pi)W_y\|^2
		+ 2\chi \left|u(0)-v(0)-2\xdot\right|\,|\llangle u_y+v_y \rrangle |
		\\
		&\quad 
		+ |\xdot| \chi \|W_y\|^2 
		+|\xdot|4\chi\left|u(0)\llangle u_y\rrangle+ v(0)\llangle v_y\rrangle \right|\\
  &\quad
  - \delta\frac{\chi}{1+\chi} \|\Pi W_y\|^2 + \delta\|(\I - \Pi)W_y\|^2 + \delta\|(\I-\Pi)W_y\| \|\Pi W_y\|\\
  &\quad +\delta |\xdot| (1+4\chi)\,\|W_y\|\, \|(\I - \Pi)W_y\| \\
         &\quad+\delta \frac{c_0}{2} |u(0;t) - v(0;t) | \|(\I-\Pi)W_y\|\notag \\ &\quad + \delta |\xdot|\left( c_1 \|(\I-\Pi)W_y\| + c_2 \|W_y\|^2 +  c_3 \|\Pi W\|^2\right)\,,
\end{align*}
which completes the proof.
\end{proof}
The two terms $- 2 \|(\I -\Pi)W_y\|^2$ and $- \delta\frac{\chi}{1+\chi} \|\Pi W_y\|^2$ are the good contributions as they give us the desired decay, and allow controlling all the other (bad) terms with carefully chosen parameter $\delta$. Indeed, they morally give the full control on the $L^2$ norm of $W_y$:
\begin{equation*}
	- 2\|(\I - \Pi)W_y\|^2 - \delta \frac{\chi}{1+\chi} \|\Pi W_y\|^2 \le -2\alpha_1 \|(\I - \Pi)W_y\|^2 - \alpha_2\|W_y\|^2 \,,
\end{equation*}
where $\alpha_1\in[0,1)$ can be chosen conveniently to control the other terms in Proposition~\ref{prop: gronwall_hypo_general_nonlocal}, and $\alpha_2:=\min\left\{2(1-\alpha_1), \delta\frac{\chi}{1+\chi}\right\}$.
The goal is to choose the first term on the right-hand side of order $1$, and to combine only a small fraction $(1-\alpha_1)$ of it with the second term to control the norm of $W_y$.


\section{\texorpdfstring{Nonlinear case with $\alpha=0$}{Nonlinear case with a=0}}\label{Sect: 6}

In this section, we consider the case $\alpha=0$ which results in $\lambda=2\chi$. This means all the weights of our functional spaces are the same and in this simplified setting we are able to provide a nonlinear stability analysis.
Note that the choice of $\alpha$ only enters the dynamics through the movement of the chemoattractant peak. For $\alpha=0$, the expression for $\xdot$ in \eqref{xdot pert} simplifies to
\begin{equation}\label{def:xdot-alpha0}
	\xdot(t) = \frac{u(0;t) - v(0;t)}{2+u(0;t)+v(0;t) } \,,
\end{equation}
and hence
\begin{align}
	\label{eq: u0 -v0 - 2xdot}
	u(0;t)-v(0;t)-2\xdot
	&=\frac{|u(0;t)|^2-|v(0;t)|^2}{2+u(0;t)+v(0;t)}\,.
\end{align}
Choosing $\langle u-v\rangle_{2\chi}(0) = 0$ and applying the third condition \eqref{eq: third cons law} together with the interpolation inequality \eqref{interpol} we have that 
	 \begin{equation}
		     \label{eq: control on u_0 - v_0}
		     |u(0;t) - v(0;t)| 
       = \left|u(0;t) - v(0;t)  - \langle u-v\rangle_{2\chi}(t)\right|\le \frac{1}{\sqrt{2\chi}}\|(\I - \Pi)W_y\| \,.
		 \end{equation}

In what follows, we derive the decay of the modified entropy in the case $\alpha=0$, starting from the bound on the entropy dissipation derived in Section~\ref{sec: hypo control}. In order to prove Theorem~\ref{thm:main1} we proceed in two steps: (1) we show in Section~\ref{sec:alpha0-entropydecay} that entropy decay holds as long as a certain smallness condition is satisfied at all times, and (2) we use this result in Section~\ref{sec:proof-main1} to obtain that the decay ensures for the smallness condition to be propagated from initial time to all later times.

\subsection{Control of nonlocal terms}\label{sec:alpha0-nonlocal}

Subsequently, we will make frequent use of the following assumption.

\begin{assumption}\label{ass:xdot-bound}
	There exists $T>0$ and $p\in(0,1)$ such that
	\begin{align*}
		\left|\xdot(t)\right|\le p \quad \text{ for all } t\in [0,T]\,,
	\end{align*}
	where $\xdot$ is defined by \eqref{def:xdot-alpha0}.
\end{assumption}

We begin by deriving explicit expressions for $\llangle u_y\rrangle$ and $\llangle v_y\rrangle$. 

\begin{prop}\label{prop:badterms}
	Let Assumption~\ref{ass:xdot-bound} hold for some $T>0$ and $p\in (0,1)$.
	Denote $$z_\pm(s,t):=\pm s+x(t)-x(t-s) \quad \text{ for }\quad  t>0\,, s\in[0,t]\,.$$ 
 Then $z_+(s,t)>0$ and $z_-(s,t)<0$. Moreover,
for any solution $(u,v)$ of \eqref{kinetic10}, we have for $t\in [0.T]$
	\begin{align*}
		\llangle  u_y \rrangle (t)
		& =\eta_\infty(z_-(t,t)) \left[u_y(z_-(t,t);0)+2\chi u(z_-(t,t);0)\right]\\
		&\quad+\int_0^t e^{-s} \eta_\infty(z_-(s,t)) \left[v_y(z_-(s,t);t-s)+2\chi v(z_-(s,t);t-s)\right]\ds\\
		&\quad+\chi\int_0^t
		e^{-s} \eta_\infty(z_-(s,t)) \left[(u_y+v_y)(z_-(s,t);t-s)+2\chi(u+v)(z_-(s,t);t-s)\right]\ds\\
		&\quad+4\chi^2\int_0^t\xdot(t-s)e^{-s} \eta_\infty(z_-(s,t))\ds
		-\frac{\chi}{1-\xdot(t)}\left[u(0;t)+v(0;t)+2\xdot(t)\right] 
	\end{align*}
	and
	\begin{align*}
		\llangle  v_y \rrangle(t) 
		& =\eta_\infty(z_+(t,t)) \left[v_y(z_+(t,t);0)-2\chi v(z_+(t,t);0)\right]\\
		&\quad+\int_0^t e^{-s}\eta_\infty(z_+(s,t)) \left[u_y(z_+(s,t);t-s)-2\chi u(z_+(s,t);t-s)\right]\ds\\
		&\quad+\chi\int_0^t
		e^{-s} \eta_\infty(z_+(s,t)) \left[(u_y+v_y)(z_+(s,t);t-s)-2\chi(u+v)(z_+(s,t);t-s)\right]\ds\\
		&\quad+4\chi^2\int_0^t\xdot(t-s)e^{-s} \eta_\infty(z_+(s,t))\ds
		+\frac{\chi}{1+\xdot(t)}\left[u(0;t)+v(0;t)-2\xdot(t)\right] \,.  
	\end{align*}
\end{prop}

\begin{proof}
    The proof of this result uses the method of characteristics, and we postpone it to Appendix~\ref{sec:char-method} to not distract from the main argument here.
\end{proof}

With these preliminary computations at hand, we are now ready to derive an upper bound for the higher-order nonlocal terms.
\begin{prop}\label{prop:badterms-est}
	Let Assumption~\ref{ass:xdot-bound} hold with $T>0$ and $p\le 1/(4\chi)$.
	If there is a constant $c>0$ such that the initial data satisfies
	\begin{align*}
		\|u(\cdot;0)\|_\infty + \|v(\cdot;0)\|_\infty+\|u_y(\cdot;0)\|_\infty+\|v_y(\cdot;0)\|_\infty\le c\,,
	\end{align*}
	then
 \begin{align*}
		&\int_0^{T} \left(|\llangle u_y(t)\rrangle|+ |\llangle v_y(t)\rrangle| \right)\dt  \le   \int_0^{T} g(t,p) \dt\\
  &\quad+ \frac{\chi}{1-p}\int_0^{T}|u(0;t)+v(0;t)+2\xdot| \dt 
  + \frac{\chi}{1-p}\int_0^{T}|u(0;t)+v(0;t)-2\xdot| \dt \\
		&\quad + 2\chi g_0(p)\int_0^{T} \left(  \|\Pi W_y\| + 2\chi  \|\Pi W\| \right) \dt
+  g_0(p) \int_0^{T} \left( \|W_y\|  + 2\chi \|W\| \right) \dt\,.
	\end{align*}
	with $r:=1+2\chi(1-p)$ and
	\begin{align}\label{def:g}
 g(t,p)&:=2(1+2\chi)ce^{-rt} + \frac{8p\chi^2}{r} (1-e^{-rt}) \,,\quad
		g_0(p):= \sqrt{\frac{4(1+p)}{r(1-p)^2}}\,.
	\end{align}
\end{prop}

\begin{proof}[Proof of Proposition~\ref{prop:badterms-est}]
	Using the explicit expression of $\llangle  u_y \rrangle$ provided by Proposition~\ref{prop:badterms}, we estimate for $t\in [0,T]$
	\begin{align*}
		\left| \llangle  u_y \rrangle \right|
		& \le (1+2\chi) ce^{-t}\eta_\infty(z_-(t,t))
		+\int_0^t e^{-s} \eta_\infty(z_-(s,t)) \left|v_y(z_-(s,t);t-s)+2\chi v(z_-(s,t);t-s)\right|\ds\\
		&\quad+\chi\int_0^t
		e^{-s} \eta_\infty(z_-(s,t)) \left|(u_y+v_y)(z_-(s,t);t-s)+2\chi(u+v)(z_-(s,t);t-s)\right|\ds\\
		&\quad+4p\chi^2\int_0^t e^{-s} \eta_\infty(z_-(s,t))\ds
		+\frac{\chi}{1-p}\left|u(0;t)+v(0;t)+2\xdot(t)\right|\,. 
	\end{align*}
	Observe that $|z_-(s,t)|=s-x(t)+x(t-s)\ge(1-p)s$ for all $s\in[0,t]$, and so with $r:=1+2\chi(1-p)$, we have
	\begin{align*}
		\left| \llangle  u_y \rrangle \right|
		& \le( 1+2\chi) ce^{-rt} 
		+\frac{4p\chi^2(1-e^{-rt})}{r}
		+\frac{\chi}{1-p}\left|u(0;t)+v(0;t)+2\xdot(t)\right|\\
		&\quad +\int_0^t e^{-rs} \left|v_y(z_-(s,t);t-s)+2\chi v(z_-(s,t);t-s)\right|\ds\\
		&\quad+\chi\int_0^t
		e^{-rs} \left|(u_y+v_y)(z_-(s,t);t-s)+2\chi(u+v)(z_-(s,t);t-s)\right|\ds\,.
	\end{align*}
	Integrating in time between $0$ and ${T}$, and using Lemma \ref{lem:badterms-aux1} from Appendix~\ref{sec:time-to-space} we get
	\begin{align*}
		&\int_0^{T} |\llangle u_y \rrangle| \dt\\
  & \quad\le  \int_0^{T} \left( 
  (1+2\chi)ce^{-rt} + \frac{4p\chi^2}{r}(1-e^{-rt})
  +
  \frac{\chi}{1-p} |u(0;t) + v(0;t) + 2\xdot| \right)\dt\\
		&\qquad+ \frac{1}{1-p} \int_0^{T} \ints{ e^{-\frac{r}{1+p}|y|}\left( |v_y| + \chi |u_y + v_y| + 2\chi^2 |u+v| + 2\chi |v| \right) } \dt\,.
	\end{align*}
	The same procedure can be done for the $v_y$ term. Summing the two results we get
	\begin{align*}
		&\int_0^{T} |\llangle u_y \rrangle| + |\llangle v_y \rrangle| \dt   \le   \int_0^{T} g(t,p) \dt\\
  &\quad+ \frac{\chi}{1-p}\int_0^{T}|u(0;t)+v(0;t)+2\xdot| \dt 
  + \frac{\chi}{1-p}\int_0^{T}|u(0;t)+v(0;t)-2\xdot| \dt \\
		&\quad+ \frac{2\chi}{1-p} \int_0^{T} \ints{ e^{-\frac{r}{1+p}|y|}\left( |u_y + v_y| + 2\chi |u+v| \right) } \dt \\
		&\quad+ \frac{1}{1-p} \int_0^{T} \ints{ e^{-\frac{r}{1+p}|y|}\left( |v_y| + |u_y| + 2\chi (|v| + |u|) \right) } \dt\,.
	\end{align*}
	Using Cauchy-Schwarz and the condition $4\chi p \le 1$, we have $\frac{r}{1+p} \ge 2\chi$, and so for any function $z(y)$ with sufficient integrability,
	\begin{align*}
		\ints{ e^{-\frac{r}{1+p}|y|} z(y) } 
  &\le \sqrt{\ints{ e^{-\frac{r}{1+p}|y|}}} \sqrt{\ints{ e^{-\frac{r}{1+p}|y|} z(y)^2}}\\
  &= \sqrt{\frac{2(1+p)}{r}} \sqrt{\ints{ e^{-\frac{r}{1+p}|y|} z(y)^2}}
\le \sqrt{\frac{2(1+p)}{r}} \|z\|\,,
	\end{align*}
 where  $\|z\| := \left(\int_\R |z(y)|^2\eta_\infty\dy\right)^{1/2}$.
Applying this idea to our last estimate, we obtain
\begin{align*}
		&\int_0^{T} \left(|\llangle u_y(t)\rrangle|+ |\llangle v_y(t)\rrangle| \right)\dt  \le   \int_0^{T} g(t,p) \dt\\
  &\quad+ \frac{\chi}{1-p}\int_0^{T}|u(0;t)+v(0;t)+2\xdot| \dt 
  + \frac{\chi}{1-p}\int_0^{T}|u(0;t)+v(0;t)-2\xdot| \dt \\
		&\quad + 2\chi \sqrt{\frac{2(1+p)}{r(1-p)^2}} \int_0^{T} \left(  \|u_y + v_y\| + 2\chi  \|u+v\| \right) \dt\\
  &\quad +  \sqrt{\frac{2(1+p)}{r(1-p)^2}} \int_0^{T} \left( \|u_y\| + \|v_y\|  + 2\chi \left(\|u\| + \|v\|\right)  \right) \dt\,.
	\end{align*}
 Notice that $\|u+v\|=\sqrt{2}\|\Pi W\|$ and $\|u\|+\|v\| \le \sqrt{2}\left(\|u\|^2+\|v\|^2\right)^{1/2}=\sqrt{2}\|W\|$ etc. This concludes the proof.
\end{proof}

It remains to control the lower-order non-local terms $\llangle u\rrangle$ and $\llangle v \rrangle$, which can be done directly with the bounds already established. 
\begin{lem}\label{lem:aux_est}
If $W$ solves \eqref{kinetic10}, then
\begin{align*}
    |\llangle u\rrangle|^2 +  |\llangle v\rrangle|^2&\le \frac{1}{2\chi} \|W_y\|^2 + 2\chi\|W\|^2\,,\\
    |\llangle u + v\rrangle|^2&\le \frac{1}{\chi} \|\Pi W_y\|^2 + 4\chi\|\Pi W\|^2\,.
\end{align*}
\end{lem}

\begin{proof}
Since $u$ is continuous at $y=0$, we compute directly using \eqref{interpol} and Jensen's inequality,
\begin{align*}
|\llangle u\rrangle|^2
&\le 2|u(0;t)-\langle u \rangle_{2\chi}|^2 + 2\langle u \rangle_{2\chi}^2
\le \frac{1}{2\chi} \|u_y\|^2 + 2\chi\|u\|^2\,,
\end{align*}
and similarly for $v$. Summing the expressions for $|\llangle u\rrangle|^2$ and $|\llangle v\rrangle|^2$, we obtain the first result. The second inequality follows by the same argument replacing $u$ by $u+v$ above and using $\|u+v\|^2 = 2\|\Pi W\|^2$ etc.
\end{proof}

\begin{cor}\label{cor:nonlocal-control}
Let $\langle u-v\rangle_{2\chi}(0) = 0$.
	Let Assumption~\ref{ass:xdot-bound} hold with $T>0$ and $p\le 1/(4\chi)$.
	If there is a constant $c>0$ such that the initial data satisfies
	\begin{align*}
		\|u(\cdot;0)\|_\infty + \|v(\cdot;0)\|_\infty+\|u_y(\cdot;0)\|_\infty+\|v_y(\cdot;0)\|_\infty\le c\,,
	\end{align*}
	then for any $\delta\in (0,1)$ there exist constants $\lambda_1, \lambda_2>0$ only depending on $p$, $\chi$ and $\delta$ such that
	\begin{equation*}
		 \int_0^{T} \big(|\llangle u_y(t)\rrangle| + |\llangle v_y(t) \rrangle|\big)\, \dt
   \le \int_0^{T} \left(\lambda_1 + \lambda_2\, \sqrt{\CL[W_y]}\right) \, \dt \,.
	\end{equation*}
\end{cor}
\begin{proof}
	The proof essentially reduces to an application of Proposition \ref{prop:badterms-est} and the Poincar\'e inequality. More precisely, together with the control on the lower-order nonlocal terms in Lemma~\ref{lem:aux_est} we obtain
  \begin{align*}
		&\int_0^{T} \left(|\llangle u_y(t)\rrangle|+ |\llangle v_y(t)\rrangle| \right)\dt  \le   \int_0^{T} g(t,p) \dt
    + \frac{2\chi}{1-p}\int_0^{T}|\llangle u(t)+v(t)\rrangle|+2|\xdot| \dt \\
		&\qquad + 2\chi g_0(p)\int_0^{T} \left(  \|\Pi W_y\| + 2\chi  \|\Pi W\| \right) \dt
+  g_0(p) \int_0^{T} \left( \|W_y\|  + 2\chi \|W\| \right) \dt\\
&\quad\le   \int_0^{T} \left( \frac{4p\chi}{1-p} + g(t,p)\right) \dt
    + \left(\frac{2\sqrt{\chi}}{1-p}+2\chi g_0(p)
    \right)\int_0^{T}
    \left(\|\Pi W_y\| + 2\chi\|\Pi W\|\right)\dt \\
		&\qquad 
+  g_0(p) \int_0^{T} \left( \|W_y\|  + 2\chi \|W\| \right) \dt\\
&\quad\le   \int_0^{T} \left( \frac{4p\chi}{1-p} + g(t,p)\right) \dt
    + 3\left(\frac{2\sqrt{\chi}}{1-p}+(1+2\chi) g_0(p)
    \right)\int_0^{T}
    \|W_y\|\dt \,.
    \end{align*}
 where we used $\langle u-v\rangle_{2\chi}(0) = 0$ together with \eqref{Poincare-3 a}. This concludes the proof
 with constants
 \begin{align*}
     \lambda_1(p,\chi)&:= 2(1+2\chi)c + \frac{4p\chi}{1-p} + 8p\chi^2\,,\\
     \lambda_2(p,\chi,\delta)&:= \frac{6\sqrt{2\chi}}{(\sqrt{1-\delta})(1-p)}\left(1+\left(\frac{1+2\chi}{2\chi}\right)\sqrt{1+4\chi}\right)\,.
 \end{align*}
 since $g(t,p)\le 2(1+2\chi)c+8p\chi^2$ for all $t\ge 0$, and using the norm equivalence \eqref{eq: equivalence entropy L2-norm}. 
 \end{proof}
\subsection{Entropy decay}\label{sec:alpha0-entropydecay}
Our strategy of proof to obtain entropy decay uses a 'security cylinder': we show (1) that entropy decay holds as long as a certain smallness condition is satisfied at all times, and (2) that the decay ensures for the smallness condition to be propagated from initial time to all later times. The following proposition is the key estimate to establish the first step. We will then combine (1) and (2) in the next section for the proof of Theorem~\ref{thm:main1}.

\begin{prop}\label{prop:entropy-decay}
Let $W$ be a solution of equation \eqref{hypo eq formulation} satisfying initially $\langle u-v\rangle_{2\chi}(0) = 0$. Let Assumption~\ref{ass:xdot-bound} hold with $p<\chi/(16(1+\chi))$ small enough and with $T>0$ small enough such 
\begin{align}\label{ass:c}
    \|\Pi W_y(t)\|+2\chi\|\Pi W(t)\|\le \sqrt{2\chi}\quad\text{ for all } t\in [0,T]\,. 
\end{align}
Then there exists constants $\mu_0=\mu_0(p,\chi)>0$, $\delta=\delta(p,\chi)\in(0,1)$ and a function $\mu(t)>0$ of the nonlocal terms $|\llangle u_y(t)\rrangle| + |\llangle v_y(t)\rrangle|$ 
such that
\begin{align*}
	\derivt \CL[W_y]  \le& - \mu_0 \CL[W_y]  + \mu(t)  \CL[W_y]
 \quad\text{ for all } t\in [0,T]\,. 
\end{align*}
\end{prop}

\begin{proof}
Starting from the expression in Proposition~\ref{prop: gronwall_hypo_general_nonlocal} combined with equation \eqref{eq: u0 -v0 - 2xdot}, substituting the expression \eqref{def:xdot-alpha0} for $\xdot$ for the $c_1$ term, applying Young's inequality with parameters $\beta_1, \beta_2>0$, using the Poincar\'e inequality~\eqref{Poincare-1 a} and bound \eqref{eq: control on u_0 - v_0}, we get
\begin{align} 
	\derivt \CL[W_y] \le &- \left[ 2 - \delta \bigg( 1+ \frac{c_0}{2\sqrt{2\chi}} + \frac{\sqrt{2\chi} c_1}{|c(u,v)|}+ \frac{\beta_1+\beta_2}{2} \bigg) \right] \, \|(\I - \Pi)W_y\|^2
 \notag\\
	&-  \delta \bigg( \frac{\chi}{1+\chi} - \frac{1}{2\beta_1} - |\xdot| \frac{c_3}{\chi^2} \bigg)\, \|\Pi W_y\|^2\notag\\
	\label{eq: lyapunov estimate}
	&+ |\xdot|\bigg( \chi+\frac{\delta(1+4\chi)}{2\beta_2} + \delta c_2\bigg)\, \|W_y\|^2
 \\
	&+2\chi \frac{\left|u(0;t)^2-v(0;t)^2\right|}{\left|2+u(0;t)+v(0;t)\right|} \left|\llangle u_y+v_y \rrangle\right|
 \notag\\
	& +4\chi|\xdot|\left|u(0;t)\llangle u_y\rrangle + v(0;t)\llangle v_y\rrangle \right|\notag
\end{align}
with $c(u,v)=4\chi+2\chi(u(0;t)+v(0;t))$ as defined in \eqref{xdot2b}. 
Thanks to Lemma~\ref{lem:aux_est} and the Poincar\'e inequality \eqref{Poincare-3 a}, we have 
\begin{align*}
    \left|u(0;t)^2-v(0;t)^2\right| \le u(0;t)^2+v(0;t)^2 
    \le  \frac{1}{2\chi} \|W_y\|^2 + 2\chi\|W\|^2
    \le\frac{5}{2\chi} \|W_y\|^2\,.
\end{align*}
As a consequence and using the expression for $\xdot$ in \eqref{def:xdot-alpha0}, we have
\begin{align*}
   &|\xdot| \left|u(0;t)\llangle u_y\rrangle + v(0;t)\llangle v_y\rrangle \right|\\
   &\quad \le \frac{2\chi|u(0;t) - v(0;t)|}{|c(u,v)| }\left(|u(0;t)| + |v(0;t)|\right)\left(|\llangle u_y\rrangle| + |\llangle v_y\rrangle| \right)\\
   &\quad \le 2\chi\frac{ u(0;t)^2+ v(0;t)^2}{|c(u,v)| }\left(|\llangle u_y\rrangle| + |\llangle v_y\rrangle| \right)
 \le \frac{ 5}{|c(u,v)| }\left(|\llangle u_y\rrangle| + |\llangle v_y\rrangle| \right) \|W_y\|^2\,.
\end{align*}
Combining these estimates, we can control the last two lines in \eqref{eq: lyapunov estimate} by
\begin{align*}
&2\chi \frac{\left|u(0;t)^2-v(0;t)^2\right|}{\left|2+u(0;t)+v(0;t)\right|} \left|\llangle u_y+v_y \rrangle\right|
+4\chi|\xdot|\left|u(0;t)\llangle u_y\rrangle + v(0;t)\llangle v_y\rrangle \right|\\
&\quad \le
\frac{ 30\chi}{|c(u,v)| }\left(|\llangle u_y\rrangle| + |\llangle v_y\rrangle| \right) \|W_y\|^2\,.
\end{align*}
Next, we make use of the lower bound on $c(u,v)$ derived in \eqref{c-bound}, together with Assumption~\eqref{ass:c}, to obtain
\begin{align*}
    |c(u,v)|\ge 4\chi - \sqrt{2\chi}\left(\|\Pi W_y\| +2\chi\|\Pi W\|\right)\ge 2\chi\,.
\end{align*}
Using Assumption~\ref{ass:xdot-bound}, we obtain the energy estimate
\begin{align*}
	\derivt \CL[W_y] \le& - \left[ 2 - \delta \bigg( 1+ \frac{c_0}{2\sqrt{2\chi}} + \frac{c_1}{\sqrt{2\chi}} + \frac{\beta_1+\beta_2}{2} \bigg) \right] \, \|(\I - \Pi)W_y\|^2 \\
	&-  \delta \bigg( \frac{\chi}{1+\chi} - \frac{1}{2\beta_1} - p \frac{c_3}{\chi^2} \bigg)\, \|\Pi W_y\|^2\\
	&+ p\bigg( \chi+\frac{\delta(1+4\chi)}{2\beta_2} + \delta c_2\bigg)\, \|W_y\|^2
	+15\,\left(|\llangle u_y\rrangle| + |\llangle v_y\rrangle| \right) \|W_y\|^2\,.
\end{align*}
Recall the definition of constants $c_0, c_1, c_2$ and $c_3$ in Lemma~\ref{lem:hypo-nonlin-bounds-2}. We choose $\beta_1=\beta_2:=(1+\chi)/\chi>0$ and $p<\chi/(8(1+\chi))$ to guarantee $\frac{\chi}{1+\chi} - \frac{1}{2\beta_1} - p \frac{c_3}{\chi^2} >0$. Next, we choose 
\begin{align}\label{def:delta}
    \delta=\delta(p,\chi):= 2\left[\bigg( 1+ \frac{c_0}{2\sqrt{2\chi}} + \frac{c_1}{\sqrt{2\chi}}+ \frac{\beta_1+\beta_2}{2} \bigg) +\frac{\chi}{2(1+\chi)}-p \frac{c_3}{\chi^2}\right]^{-1}
\end{align}
Then $\delta< 2/\left( 1+ \frac{c_0}{2\sqrt{2\chi}} + \frac{c_1}{\sqrt{2\chi}} + \frac{\beta_1+\beta_2}{2}\right) <1$. Finally, we define
\begin{align*}
    \eta(p,\chi):=2 - \delta \bigg( 1+ \frac{c_0}{2\sqrt{2\chi}} + \frac{c_1}{\sqrt{2\chi}} + \frac{\beta_1+\beta_2}{2} \bigg)
   =
     \delta \bigg( \frac{\chi}{1+\chi} - \frac{1}{2\beta_1} - p \frac{c_3}{\chi^2} \bigg)
\end{align*}
and choose $p$ small enough such that $\eta>p\bigg( \chi+\frac{\delta(1+4\chi)}{2\beta_2} + \delta c_2\bigg)$. We obtain 
\begin{align*}
	\derivt \CL[W_y] \le& - \left[\eta -p\bigg( \chi+\frac{\delta(1+4\chi)}{2\beta_2} + \delta c_2\bigg) - 15\,\left(|\llangle u_y\rrangle| + |\llangle v_y\rrangle| \right) \right] \|W_y\|^2\\
 \le& - \frac{2}{(1-\delta)}\left[\eta -p\bigg( \chi+\frac{\delta(1+4\chi)}{2\beta_2} + \delta c_2\bigg) - 15\,\left(|\llangle u_y\rrangle| + |\llangle v_y\rrangle| \right) \right] \CL[W_y]
\end{align*}
thanks to the norm equivalence \eqref{eq: equivalence entropy L2-norm}.
This concludes the proof with
\begin{align*}
    &\mu_0:=\frac{2}{(1-\delta)}\left[\eta -p\bigg( \chi+\frac{\delta(1+4\chi)}{2\beta_2} + \delta c_2\bigg)\right]\,, \\
    &\mu(t):=\frac{30}{(1-\delta)}\,\left(|\llangle u_y(t)\rrangle| + |\llangle v_y(t)\rrangle| \right)\,.
\end{align*}
\end{proof}

\subsection{Proof of Theorem~\ref{thm:main1}}\label{sec:proof-main1}

\begin{prop}\label{prop:main1-part1}
Let $W$ be a solution of equation \eqref{hypo eq formulation} satisfying initially $\langle u-v\rangle_{2\chi}(0) = 0$, and
	\begin{align*}
		\|u(\cdot;0)\|_\infty + \|v(\cdot;0)\|_\infty+\|u_y(\cdot;0)\|_\infty+\|v_y(\cdot;0)\|_\infty\le c
	\end{align*}
 for some small enough constant $0<c\le 2\chi/(1+2\chi)$.
 Let Assumption~\ref{ass:xdot-bound} hold with some $\tilde T>0$ and $p\le \min\left\{ \frac{1}{4\chi}\,;\, \frac{\chi}{8(1+\chi)} \right\}$ small enough. Let $T^*$ denote the first time when $ \|\Pi W_y(\cdot)\|+2\chi\|\Pi W(\cdot)\|= \sqrt{2\chi}$ and define $T:=\min\{ T^*,\tilde T\}$.
Then there exists a rate $\gamma=\gamma(p,c,\chi)>0$ and a constant $C=C(\CL[W_y(0)],p, c,\chi)>1$ such that
\begin{align*}
	 \|W_y(t)\| \le  C  \|W_y(0)\| \exp\left(-\gamma t\right) \quad\text{ for all } t\in [0,T]\,. 
\end{align*}
\end{prop}

\begin{proof}
Denote $\gamma:=\left(\mu_0-(\frac{30}{1-\delta}) \lambda_1\right)/2$ with $\mu_0$, $\delta$ and $\lambda_1$ as defined in the proofs of Corollary~\ref{cor:nonlocal-control} and Proposition~\ref{prop:entropy-decay}.
By choosing $p$ and $c$ small enough, we can guarantee that $\gamma>0$.
 Note that $c< 2\chi/(1+2\chi)$ guarantees that \eqref{ass:c} is satisfied initially with strict inequality. By definition of $T$, we have $\|\Pi W_y(t)\|+2\chi\|\Pi W(t)\|\le \sqrt{2\chi}$ for all $t\in [0,T]$ and Assumption~\ref{ass:xdot-bound} holds for all $t\in [0,T]$. 
 Integrating the decay estimate in Proposition~\ref{prop:entropy-decay} and writing $\CL[W_y](t)=\CL(t)$ for short, we have for all $t\in [0,T]$ that it holds
  \begin{align*}
  &\derivt\left[\CL\left(t\right) \exp\left(\mu_0 t-\int_0^{t} \mu(s)\,\ds\right)\right] \\
  &\qquad= \exp\left(\mu_0t-\int_0^{t} \mu(s)\,\ds\right) \left[\dot{\CL}\left(t\right) + \mu_0\CL\left(t\right) - \mu\left(t\right)\CL\left(t\right)\right]\le 0\,.
\end{align*}
Then the control of the nonlocal terms from Corollary~\ref{cor:nonlocal-control} yields
  \begin{align*}
  &\derivt\left[\CL\left(t\right) \exp\left((\mu_0-\lambda_1')t-\lambda_2'\int_0^{t} \sqrt{\CL(s)}\,\ds\right)\right] \le 0 \quad \text{ for all } t\in [0,T]
\end{align*}
with $\lambda_i':=(\frac{30}{1-\delta}) \lambda_i$ for $i=1,2$ and for $\delta$ as defined in \eqref{def:delta} in the proof of Proposition~\ref{prop:entropy-decay}. Denoting $q(t):= \int_0^{t} \sqrt{\CL(s)}\,\ds$, we conclude
  \begin{align}\label{eq:entropy-bound-1}
  &\CL\left(t\right) \exp\left((\mu_0-\lambda_1')t\right)\le 
  \CL(0) \exp\left(\lambda_2'q(t)\right) \quad \text{ for all } t\in [0,T]\,.
\end{align}
We claim that $q(t)$ is uniformly bounded above for all $t\in [0,T]$. Indeed, it follows from \eqref{eq:entropy-bound-1} that 
  \begin{align*}
  &\CL\left(t\right) \le 
  \left[\CL(0) \sup_{0\le s \le T}\exp\left(\lambda_2'q(s)\right)\right] \exp\left(-(\mu_0-\lambda_1')t\right) \quad \text{ for all } t\in [0,T]\,.
\end{align*}
Taking square roots on both sides and integrating, we obtain
  \begin{align*}
  q(t)&=\int_0^{t} \sqrt{\CL(s)}\,\ds \le 
  \frac{2\sqrt{\CL(0)}}{(\mu_0-\lambda_1')} \sup_{0\le s \le t}\exp\left(\frac{\lambda_2'}{2}q(s)\right)\left[ 1- \exp\left(-\frac{(\mu_0-\lambda_1')t}{2}\right)\right]\\
  &\le \frac{\sqrt{\CL(0)}}{\gamma} \exp\left(\frac{\lambda_2'}{2}q(t)\right)\left[ 1- \exp\left(-\gamma t\right)\right]
  \quad \text{ for all } t\in [0,T]\,,
\end{align*}
where we used the fact that $q(\cdot)$ is an increasing function and the definition of $\gamma:=\left(\mu_0- \lambda_1'\right)/2$. Denote
\begin{align*}
   0\le R(t):=\frac{\sqrt{\CL(0)}}{\gamma} \left[ 1- \exp\left(-\gamma t\right)\right]
   \le\frac{\sqrt{\CL(0)}}{\gamma}\,.
\end{align*}
Then
\begin{align}\label{eq:q-cond}
    q(t)&\le R(t)\exp\left(\frac{\lambda_2'}{2}q(t)\right) 
    \le \frac{\sqrt{\CL(0)}}{\gamma}\exp\left(\frac{\lambda_2'}{2}q(t)\right) \,.
\end{align}
Notice that the equation $x=\frac{\sqrt{\CL(0)}}{\gamma}\exp\left(\frac{\lambda_2'}{2}x\right)$ has exactly two roots in $x>0$, say $0<q_1<q_2$. In other words, condition \eqref{eq:q-cond} means that either $q(t)\le q_1$, or $q(t)\ge q_2$. Further, $q(\cdot)$ is a continuous function due to its integral definition. As a result, if $q(0)<q_1$, then $q(t)\le q_1$ for all $t\in [0,T]$. Returning to expression \eqref{eq:entropy-bound-1} with this uniform bound on $q(\cdot)$, we arrive at
  \begin{align*}
  &\CL\left(t\right) \le 
  \CL(0) \exp\left(\lambda_2'q_1\right)\exp\left(-2\gamma t\right)\,,\quad \text{ for all } t\in [0,T]\,.
\end{align*} 
By norm equivalence \eqref{eq: equivalence entropy L2-norm},
\begin{align*}
	\left(\frac{1-\delta}{2}\right) \|W_y(t)\|^2 &\le \CL[t] 
 \le \CL(0) \exp\left(\lambda_2'q_1\right)\exp\left(-2\gamma t\right)\\
 &\le \left(\frac{1+\delta}{2}\right) \exp\left(\lambda_2'q_1\right) \|W_y(0)\|^2\exp\left(-2\gamma t\right)  \,,
\end{align*}
and so we conclude
\begin{align*}
	 \|W_y(t)\| \le  C  \|W_y(0)\| \exp\left(-\gamma t\right)
   \quad \text{ for all } t\in [0,T]
\end{align*}
with hypocoercivity constant
\begin{equation}\label{eq:constantC}
C:=\sqrt{\frac{1+\delta}{1-\delta}} \exp\left(\frac{\lambda_2'q_1}{2}\right)>1
\end{equation}
only depending on $\CL(0)$, $c$, $p$ and $\chi$.
\end{proof}

\begin{cor}[Convergence of the centre]\label{cor:main1-part2}
Let $W$ be a solution of equation \eqref{hypo eq formulation} satisfying initially $\langle u-v\rangle_{2\chi}(0) = 0$, and
	\begin{align*}
		\|u(\cdot;0)\|_\infty + \|v(\cdot;0)\|_\infty+\|u_y(\cdot;0)\|_\infty+\|v_y(\cdot;0)\|_\infty\le c
	\end{align*}
 for some small enough constant $0<c\le \min\{ 2\chi/(1+2\chi)\,;\, \chi/(3C)\}$ with $C$ as defined in \eqref{eq:constantC}.
 Let Assumption~\ref{ass:xdot-bound} hold with some $\tilde T>0$ and $p\le \min\left\{ \frac{1}{4\chi}\,;\, \frac{\chi}{8(1+\chi)} \right\}$ small enough.
Then
\begin{align*}
	 |\xdot(t)| \le \frac{3C c}{\chi} \,e^{-\gamma t} \quad \text{ for all } t\in [0,T]
\end{align*}
with $T>0$ and $\gamma>0$ as defined in Proposition~\ref{prop:main1-part1}.
\end{cor}

\begin{proof}
 Thanks to the decay estimate from Proposition~\ref{prop:main1-part1} and Poincar\'e inequality \eqref{Poincare-1 a}, we deduce that for any $t\in [0,T]$,
\begin{align*}
    \|\Pi W_y(t)\|+2\chi\|\Pi W(t)\|
    &\le 3\|\Pi W_y(t)\|\le 3 \| W_y(t)\|\\
    &\le 3C \|W_y(0)\|e^{-\gamma t}<3 C\|W_y(0)\| \le3\sqrt{\frac{2 }{\chi} }Cc
\end{align*}
and so choosing $c<\chi/(3C)$ small enough, we can guarantee that $\sqrt{8\chi}-\|\Pi W_y(t)\|-2\chi\|\Pi W(t)\| \ge \sqrt{2\chi}$. We recall the dynamics \eqref{def:xdot-alpha0} of the centre $\xdot$. Using the above, the bound from Proposition~\ref{prop:xdotbound} on $\xdot$ simplifies to
	\begin{align*}
|\xdot(t)|&\le \frac{\|(\I-\Pi)W_y(t)\|+2\chi\|(\I-\Pi)W(t)\|}{\sqrt{8\chi}-\|\Pi W_y(t)\|-2\chi\|\Pi W(t)\|}
\le \frac{3\|(\I-\Pi)W_y(t)\|}{\sqrt{2\chi}}
\le \frac{3\|W_y(t)\|}{\sqrt{2\chi}}\\
&\le \frac{3C \|W_y(0)\|}{\sqrt{2\chi}} \,e^{-\gamma t}
\le \frac{3C c}{\chi} \,e^{-\gamma t}
	\end{align*}
 for any $t\in [0,T]$. We conclude that $\xdot$ is integrable on $[0,t]$.
\end{proof}

\begin{proof}[Proof of Theorem~\ref{thm:main1}]
    Let $\bar f(t)$ be the reformulation in the moving frame $y=x-\x(t)$ solving \eqref{kinetic2} with initial condition satisfying $\iint v\bar f_0\dy\dmu(v)=0$ and \eqref{ass:initial}. The latter is equivalent to 
	\begin{align*}
		\|u(\cdot;0)\|_\infty + \|v(\cdot;0)\|_\infty+\|u_y(\cdot;0)\|_\infty+\|v_y(\cdot;0)\|_\infty\le c
	\end{align*}
 with $c:=\eps_0$. 
Take $\eps_0=\eps_0(\chi)>0$ small enough such that  $$
0<\eps_0 < \min\left\{ \frac{2\chi}{(1+2\chi)}\,;\, \frac{\chi}{3C}\,;\, \frac{1}{12C}\,;\, \frac{\chi^2}{24 C (1+\chi)}\right\}
$$
with $C$ as defined in \eqref{eq:constantC}. Then Assumption~\ref{ass:xdot-bound} is satisfied  at $t=0$ with strict inequality with  $p\le \min\left\{ \frac{1}{4\chi}\,;\, \frac{\chi}{8(1+\chi)} \right\}$ thanks to Corollary~\ref{cor:main1-part2}. Denote 
$$\tilde T:=\inf\left\{s>0 \,:\, |\xdot(s)|=\min\left\{ \frac{1}{4\chi}\,;\, \frac{\chi}{8(1+\chi)} \right\}\right\}\,.$$
It follows from Proposition~\ref{prop:main1-part1} that
\begin{align*}
	 \|W_y(t)\| \le  C  \|W_y(0)\| \exp\left(-\gamma t\right) \quad\text{ for all } t\in [0,T]\,,
\end{align*}
with  $T:=\min\{ T^*,\tilde T\}$, where $T^*$ denotes the first time when $ \|\Pi W_y(\cdot)\|+2\chi\|\Pi W(\cdot)\|= \sqrt{2\chi}$. If $T^*<\tilde T$, then $ \|\Pi W_y(T^*)\|+2\chi\|\Pi W(T^*)\| \le C\|W_y(0)\|e^{-\gamma T^*}<\sqrt{2\chi}$ following the proof of Corollary~\ref{cor:main1-part2}, which contradicts the definition of $T^*$. Hence $T=\tilde T$. Assume $\tilde T<\infty$. Then from Corollary~\ref{cor:main1-part2}, 
\begin{equation*}
    \min\left\{ \frac{1}{4\chi}\,;\, \frac{\chi}{8(1+\chi)} \right\} =|\xdot(\tilde T)| =|\xdot(T)| \le \frac{3C \eps_0}{\chi} \,e^{-\gamma T} \le \frac{3C \eps_0}{\chi}.
\end{equation*}
This is a contradiction with the upper bound on $\eps_0$. Hence $\tilde T=+\infty$.\\
We conclude applying Proposition~\ref{prop:main1-part1} with $\tilde T=+\infty$ and noting that
\begin{align*}
\|\bar f(t) - \bar f_\infty\|^2_{H^1} &= \|W(t)\|^2 + \|W_y(t)\|^2
\le \left(1+\frac{1}{\chi^2}\right) \|W_y(t)\|^2 \\
&\le \left(1+\frac{1}{\chi^2}\right) C^2 \|W_y(0)\|^2 e^{-2\gamma t} 
\le C_0 \|\bar f_0 - \bar f_\infty\|^2_{H^1} e^{-2\gamma t}
\end{align*}
for all $t\ge 0$
with $C_0:=\left(1+\frac{1}{\chi^2}\right) C^2$ and where $\bar f(t)$ solves \eqref{kinetic2}. This concludes the first part of the theorem statement. For the second part regarding the convergence of the center, the exponential decay of $|\xdot|$ for all $t\ge 0$ guaranteed by Corollary~\ref{cor:main1-part2} immediately implies that $\xdot$ is integrable between $0$ and $+\infty$, which concludes the existence of a limit $\x_\infty$.
\end{proof}

\section{\texorpdfstring{Linear case with $\alpha>0$}{Linear case with a>0}}\label{Sect: 7}

Recall the operator splitting from equation~\eqref{hypo eq formulation},
\begin{equation*}
	\partial_t W_y + \T W_y = \L W_y + \RR(W)\,. 
\end{equation*}
The only nonlinear terms in this equation appear in the operator $\RR(W)$ and correspond precisely to 
$$
\xdot \Q[W]:=  \xdot\partial_y\left(\eta_\infty^{-1}\partial_y\left(\eta_\infty W\right)\right)
$$
with the operator $\Q$ as in \eqref{def:Q}. For the statement in Theorem~\ref{thm:main2}, we discard the nonlinear terms in the evolution, and consider the solution $W$ of the following linearized equation instead
\begin{equation}\label{eq:lin}
   \partial_t W_y + \T W_y = \L W_y + 2\chi\delta_0(u-v-2\xdot_{lin})\1\,,
\end{equation}
where $\xdot_{lin}$ is the linearization of $\xdot$,
\begin{equation}
	\label{xdot lin}
	\xdot_{lin}(t) = \frac{\lambda}{4\chi}\left(u(0;t) - v(0;t) - \frac{\sqrt{\alpha}}{2} \int_{\R^d} (u-v) e^{-\lambda|y|} \dy   \right)\,.
\end{equation}
We obtain the following expression for the entropy dissipation for solutions $W$ of the linearized system \eqref{eq:lin}.
\begin{prop}\label{prop:H1dot-lin}
	If $W$ solves equation~\eqref{eq:lin}, then
	\begin{align*}
		\frac12\frac{\dd}{\dt} \|W_y\|^2
		&=  - 2\|(\I-\Pi)W_y\|^2
		+ 2\chi \left[u(0)-v(0)-2\xdot_{lin}\right]\llangle u_y+v_y \rrangle \,.
	\end{align*}
\end{prop}
\begin{proof}
   The argument is analogous to Proposition~\ref{prop:H1dot}. 
\end{proof}
The only dependence on $\alpha>0$ in the linearized system \eqref{eq:lin} is via $\xdot_{lin}$. We start with an auxiliary lemma rewriting this term.
\begin{lem}\label{lem:lin-estimate}
	For any $t>0$ it holds
	\begin{align*}
		&2\chi \left[u(0;t)-v(0;t)-2\xdot_{lin}(t)\right]\llangle u_y+v_y \rrangle(t)\\
		&=
		\frac{\sqrt{\alpha}}{2} \dfrac{\dd}{\dt}  \llangle u-v\rrangle^2(t) 
		+2\sqrt{\alpha} \llangle u-v\rrangle^2(t) 
		-\sqrt{\alpha}e^{-2t}\langle u-v \rangle_\lambda|_{t=0} \left[\dfrac{\dd}{\dt}\llangle u - v\rrangle(t)+2\llangle u-v\rrangle(t)\right]\,.  
	\end{align*}
\end{lem}
\begin{proof}
The linearized equation~\eqref{eq:lin} for $W_y$ can be derived from the following equation on $W$, obtained by linearizing \eqref{kinetic11}:
\begin{subequations}
	\label{L+Q}
	\begin{align}
		\label{lina}
		\partial_t u + \partial_y u &=- \xdot_{lin} 2\chi\sy  - (1 - \chi \sy) (u-v)\,,\\
		\label{linb}    
		\partial_t v - \partial_y v &=- \xdot_{lin}
  2\chi\sy  + (1 + \chi \sy) (u-v) \,.
	\end{align}
\end{subequations}
	From \eqref{L+Q} it follows that
	\begin{equation}\label{eq:u-v:lin}
		\partial_t(u-v)=-(u_y+v_y)-2(u-v)\,.
	\end{equation}
	Therefore, using the definition of $\llangle\cdot\rrangle$ directly, we can calculate explicitly
	\begin{align*}
		\frac{\dd}{\dt}\llangle u-v\rrangle
		&= \lim_{\eps\to0} \iint K(y)K(z)\partial_t(u-v)(\eps(y-z))\dy\dz\\
		&= - \lim_{\eps\to0} \iint K(y)K(z)\left[(u_y+v_y)+2(u-v)\right](\eps(y-z))\dy\dz\\
		&= - \llangle u_y+v_y\rrangle-2\llangle u-v\rrangle\,.
	\end{align*}
	We conclude that 
	\begin{equation*}
		\llangle u_y + v_y\rrangle  = -\dfrac{\dd}{\dt} \llangle u - v\rrangle-2\llangle u-v\rrangle\,.
	\end{equation*}
	Using the explicit expression in \eqref{xdot lin}, we can simplify the last term in the entropy dissipation as follows:
	\begin{align*}
		& 2\chi \left[u(0)-v(0)-2\xdot_{lin}\right]\llangle u_y+v_y \rrangle\\
		&\qquad=
		-\sqrt{\alpha}\left[ u(0)-v(0)-\langle u-v \rangle_\lambda  \right]  \llangle u_y+v_y \rrangle  \\
		&\qquad=
		-\sqrt{\alpha}\left[ \llangle u-v\rrangle-\langle u-v \rangle_\lambda \right]  \left[-\dfrac{\dd}{\dt} \llangle u - v\rrangle-2\llangle u-v\rrangle \right]  \\    
		&\qquad=
		\frac{\sqrt{\alpha}}{2} \dfrac{\dd}{\dt}  \llangle u-v\rrangle^2 
		+2\sqrt{\alpha} \llangle u-v\rrangle^2 
		-\sqrt{\alpha}\langle u-v \rangle_\lambda  \left[\dfrac{\dd}{\dt}\llangle u - v\rrangle+2\llangle u-v\rrangle\right] \,.  
	\end{align*}
	Using~\eqref{eq:u-v:lin} together with the second conservation law $\langle u_y+ v_y\rangle_{\lambda}=0$  from Corollary~\ref{cor:CSL-pert}
 we immediately have
	\begin{align*}
		\frac{\dd}{\dt}\langle u-v\rangle_\lambda = -2\langle u-v\rangle_\lambda \,,
	\end{align*}
	and so $\langle u-v\rangle_\lambda$ decays exponentially (also see Lemma~\ref{lem:cond-initial}). 
	Hence,
	\begin{align*}
		& 2\chi \left[u(0)-v(0)-2\xdot_{lin}\right]\llangle u_y+v_y \rrangle\\
		&\quad=
		\frac{\sqrt{\alpha}}{2} \dfrac{\dd}{\dt}  \llangle u-v\rrangle^2 
		+2\sqrt{\alpha} \llangle u-v\rrangle^2 
		-\sqrt{\alpha}e^{-2t}\langle u-v \rangle_\lambda|_{t=0} \left[\dfrac{\dd}{\dt}\llangle u - v\rrangle+2\llangle u-v\rrangle\right]\,.  
	\end{align*}
\end{proof}

As a result, we are able to obtain energy decay for the $H^1$-entropy with correction term.

\begin{prop}[Coercivity of the local equilibria]\label{prop:coercivity-lin}
	Assume  $\langle u-v\rangle_\lambda|_{t=0}=0$.
	Then the following inequality holds
	\begin{align*}
		&\frac12\frac{\dd}{\dt} \left(\|W_y\|^2- \sqrt{\alpha} \llangle u-v\rrangle^2 \right)
		\le  - \frac{2\chi+\sqrt{\alpha}}{\chi + \sqrt{\alpha}}\|(\I-\Pi)W_y\|^2\,.
	\end{align*}
\end{prop}

\begin{remark}\label{rmk:IC}
	We can remove the assumption $\langle u-v\rangle_\lambda|_{t=0}=0$ at the price of obtaining additional terms in the estimates on the right-hand side of Proposition~\ref{prop:coercivity-lin} of order $O(e^{-2t})$. This follows from the fact that  $\langle u-v\rangle_\lambda(t) \le e^{-2t} \langle u-v\rangle_\lambda(0)$ up to nonlinear terms. In this case however a uniform bound on $\llangle u_y+v_y\rrangle$ is needed.
\end{remark}

\begin{proof}
	Using the assumption  $\langle u-v\rangle_\lambda|_{t=0}=0$, Lemma~\ref{lem:cond-initial} and Lemma~\ref{lem:interpol}, we have 
	\begin{align}\label{eq: <<u-v>>^2}
		\llangle u-v\rrangle^2 
  = \left|u(0;t)-v(0;t)-\langle u-v\rangle_\lambda\right|^2
		&\le\frac{1}{2\chi + 2\sqrt{\alpha}} \|(\I-\Pi)W_y\|^2\,.
	\end{align}
	Combining estimate \eqref{eq: <<u-v>>^2} with Proposition~\ref{prop:H1dot-lin} and Lemma~\ref{lem:lin-estimate} we conclude the proof:
	\begin{align*}
		\frac12\frac{\dd}{\dt} \left(\|W_y\|^2- \sqrt{\alpha} \llangle u-v\rrangle^2 \right)
		&=  - 2\|(\I-\Pi)W_y\|^2
		+2\sqrt{\alpha} \llangle u-v\rrangle^2 \\
		&\le  - \left(\frac{2\chi+\sqrt{\alpha}}{\chi + \sqrt{\alpha}} \right)\|(\I-\Pi)W_y\|^2\,.
	\end{align*}
\end{proof}

Now we are ready to conclude for the proof of Theorem~\ref{thm:main2}. We split the argument in two cases: Step 1) convergence of the shape, and Step 2) convergence of the chemoattractant peak.

\begin{proof}[Proof of Theorem~\ref{thm:main2} - Step 1: convergence of the shape]
Let $\bar f(t)$ be the reformulation in the moving frame $y=x-\x(t)$ solving the linearized version of \eqref{kinetic2} with initial condition satisfying $\iint v\bar f_0 e^{-\sqrt{\alpha}|y|}\dy\dmu(v)=0$ for $\alpha>0$. Consider instead perturbations $(u,v)$ solving the linearized equation \eqref{L+Q} with initial condition satisfying $\langle u-v\rangle_\lambda|_{t=0}=0$, which implies that $W_y$ solves \eqref{eq:lin}.
	To handle the nonlocal terms involving $\alpha$, we consider here the following \rus{Lyapunov} functional:
	\begin{align}\label{def: lyapunv functional alpha positive}
		\CL_\alpha[W_y] &:=
  \CL[W_y]- \frac{\sqrt{\alpha}}{2} \llangle u-v\rrangle^2
  =\frac{1}{2}\|W_y\|^2- \frac{\sqrt{\alpha}}{2} \llangle u-v\rrangle^2 + \delta \sprod{\A W_y}{W_y}_{2\chi} \,.
	\end{align}
	Recalling the operator splitting \eqref{hypo eq formulation}, the nonlinear terms in the evolution of $W_y$ all appear in the remainder $\RR$. More precisely, we use the notation from \eqref{R-splitting} and denote $\RR_{lin}=\RR_1+\RR_3$ the linear part of $\RR$, where
 	\begin{equation*}
		\RR_1(W):= 2\chi \delta_0 (u-v) \1\,,\quad
     \RR_3(W):=-\xdot_{lin} 2\chi\sy  \1\,.
	\end{equation*}
By direct inspection, we have
	\begin{align*}
	    \sprod{\A \RR_{lin}(W)}{W_y}_{2\chi}=0
	\end{align*}
following the arguments in the proof of Lemma~\ref{lem:hypo-nonlin-bounds-1} for the linear terms only. Similarly, to estimate $\sprod{ \RR_{lin}(W)}{\A W_y}_{2\chi}$ we follow the argument in Lemma~\ref{lem:hypo-nonlin-bounds-2} to obtain
\begin{align*}
		\sprod{\RR_{lin}(W)}{\A W_y}_{2\chi} 
		&\le \frac{c_0}{2}|u(0;t) - v(0;t) |\|(\I-\Pi)W_y\| \\
  &\quad+ 4 |\xdot_{lin}| \sqrt{\chi}\, (1+2\chi) \|(\I-\Pi)W_y\| \,.
\end{align*}
Using Proposition~\ref{prop:coercivity-lin}, and Lemma~\ref{lem:delta-dissipation} together with Corollary~\ref{cor:hypo-lin-bounds} as well as the above estimates on the linear part of the remainder $\RR_{lin}$, we have
\begin{align*}
 \derivt \CL_\alpha[W_y]
 &= \frac12\frac{\dd}{\dt} \left(\|W_y\|^2- \sqrt{\alpha} \llangle u-v\rrangle^2 \right)
 + \delta \frac{\dd}{\dt}\sprod{\A W_y}{W_y}_{2\chi}\\
 &\le - \frac{2\chi+\sqrt{\alpha}}{\chi + \sqrt{\alpha}}\|(\I-\Pi)W_y\|^2 - \delta D(t) + \delta \sprod{\RR_{lin}(W)}{\A W_y}_{2\chi} \\
&\le - \frac{2\chi+\sqrt{\alpha}}{\chi + \sqrt{\alpha}}\|(\I-\Pi)W_y\|^2\\
&\quad - \delta \frac{\chi}{1+\chi} \|\Pi W_y\|^2 + \delta\|(\I - \Pi)W_y\|^2 + \delta\|(\I-\Pi)W_y\| \|\Pi W_y\| \\
&\quad + \frac{\delta c_0}{2}|u(0;t) - v(0;t) |\|(\I-\Pi)W_y\|+ 4 \delta|\xdot_{lin}| \sqrt{\chi}\, (1+2\chi) \|(\I-\Pi)W_y\|\,.
\end{align*}
Further, since $\langle u-v\rangle_\lambda|_{t=0}=0$, equation~\eqref{xdot lin} reduces to
\begin{equation*}
	|\xdot_{lin}(t)| = \frac{\lambda}{4\chi}\left|u(0;t) - v(0;t)\right| \,.
\end{equation*}
Using Young's inequality and the estimate on $|u(0;t) - v(0;t)|$ in \eqref{eq: <<u-v>>^2} we get
\begin{align*}
 \derivt \CL_\alpha[W_y]
&\le - \left(\frac{2\chi+\sqrt{\alpha}}{\chi + \sqrt{\alpha}}-\delta c_0'\right)\|(\I-\Pi)W_y\|^2
- \delta \frac{\chi}{2(1+\chi)} \|\Pi W_y\|^2\\
&= - \left(c_1'-\delta c_0'\right)\|(\I-\Pi)W_y\|^2
- \delta c_2' \|\Pi W_y\|^2\,,
\end{align*}
where
\begin{gather*}
    c_0'(\chi,\alpha):= 1+\frac{(1+\chi)}{2\chi}+ \frac{c_0}{2\sqrt{2\chi+\sqrt{\alpha}}} + \frac{\lambda(1+2\chi)}{\sqrt{2\chi}\sqrt{\chi+\sqrt{\alpha}}}\,,\\
    c_1'(\chi,\alpha):=\frac{2\chi+\sqrt{\alpha}}{\chi + \sqrt{\alpha}}\,,\quad c_2'(\chi):=\frac{\chi}{2(1+\chi)}\,.
\end{gather*}
Choosing $\delta=\delta(\chi,\alpha):=\min\left\{ \frac{1}{2}\,,\, \frac{c_1'}{c_0'+c_2'}\right\}$,
we obtain
\begin{align*}
 \derivt \CL_\alpha[W_y]
&\le - \delta c_2' \|W_y\|^2 \,,
\end{align*}
and using the norm equivalence \eqref{eq: equivalence entropy L2-norm} and the positiveness of $\llangle u-v \rrangle^2$ we get the wished result
\begin{align*}
    \derivt \CL_\alpha[W_y] &\le - \frac{2 \delta c_2'}{1 + \delta} \left( \CL[W_y] - \frac{\sqrt{\alpha}}{2} \llangle u-v\rrangle^2 \right) 
    = - 2\gamma \CL_\alpha[W_y] 
\end{align*}
with $\gamma:=\frac{\delta c_2'}{1+\delta}$.
By Gr\"onwall's inequality, we conclude for exponential decay in entropy,
	\begin{align*}
		\CL_\alpha[W_y](t) \le e^{- 2\gamma t} \CL_\alpha[W_y](0)\,.
	\end{align*}
Finally, we observe that the norm equivalence \eqref{eq: equivalence entropy L2-norm} carries over to the perturbed entropy $\CL_\alpha$ thanks to the estimate \eqref{eq: <<u-v>>^2} and the bound $\|(\I-\Pi) W_y\|^2\le \|W_y\|^2$:
 	\begin{equation}\label{eq: equivalence entropy L2-norm alpha positive}
    \left(\frac{1-\delta}{2} 
  - \frac{\sqrt{\alpha}}{4(\chi+\sqrt{\alpha})}\right)
  \|W_y\|^2 \le \CL_\alpha[W_y] \le \frac{1+\delta}{2} \|W_y\|^2\,,
	\end{equation}
where the l.h.s. is strictly positive thanks to $\delta\le \frac{1}{2}$.
As a direct consequence, we obtain the desired decay estimate
	\begin{align}\label{Wy-decay-alpha-pos}
		\|W_y(t)\| \le C e^{-\gamma t} \|W_y(0)\| \,,\quad \text{ with }
  C:= \sqrt{\frac{1+\delta}{\left(1-\delta 
  - \frac{\sqrt{\alpha}}{2(\chi+\sqrt{\alpha})}\right)}}\,.
	\end{align}
 We conclude the first step of Theorem~\ref{thm:main2} by noting that
\begin{align*}
\|\bar f(t) - \bar f_\infty\|^2_{H^1} &= \|W(t)\|^2 + \|W_y(t)\|^2
\le \left(1+\frac{1}{\chi^2}\right) \|W_y(t)\|^2 \\
&\le \left(1+\frac{1}{\chi^2}\right) C^2 \|W_y(0)\|^2 e^{-2\gamma t} 
\le C_\alpha^2 \|\bar f_0 - \bar f_\infty\|^2_{H^1} e^{-2\gamma t}
\end{align*}
for all $t\ge 0$
with $C_\alpha:=\sqrt{\left(1+\frac{1}{\chi^2}\right)} C$.
\end{proof}

\begin{proof}[Proof of Theorem~\ref{thm:main2} - Step 2: convergence of the centre]
 Take $\eps_0=\eps_0(\chi,\alpha)>0$ small enough such that
\begin{equation}\label{cond-initial-alpha-pos}
\eps_0\le \frac{2\chi\sqrt{2(\chi+\sqrt{\alpha})}}{3C\lambda} 
\end{equation}
Then, using the Poincar\'e inequality \eqref{Poincare-1} from Proposition~\ref{prop:Poincare} and the decay estimate \eqref{Wy-decay-alpha-pos}, we have
for all $t\ge 0$
\begin{align*}
\|\Pi W_y(t)\|+2\chi\|\Pi W(t)\|
\le 3\|\Pi W_y(t)\|\le 3 \|W_y(t)\|
\le 3C\|W_y(0)\|\le 3C \eps_0 <\frac{\mu}{2}\,.
\end{align*}
Hence, we can apply estimate \eqref{eq:xdot-bound} from Proposition~\ref{prop:xdotbound} and the improved Poincar\'e inequality from Corollary~\ref{cor:Poincare-improved-est} together with Lemma~\ref{lem:cond-initial} and 
$\langle u-v\rangle_\lambda|_{t=0}=0$
to obtain
\begin{align*}
|\xdot(t)|&\le \frac{\|(\I-\Pi)W_y(t)\|+2\chi\|(\I-\Pi)W(t)\|}{\mu-\|\Pi W_y(t)\|-2\chi\|\Pi W(t)\|}\\
&\le \frac{\|(\I-\Pi)W_y(t)\|+2\chi\|(\I-\Pi)W(t)\|}{\mu/2}\\
&\le \frac{3\|(\I-\Pi)W_y(t)\|}{\mu/2}
\le \frac{6\|W_y(t)\|}{\mu}
\le \frac{6C\|W_y(0)\|}{\mu}e^{-\gamma t}\,.
\end{align*}
This exponential decay of $|\xdot|$ immediately implies that $\xdot$ is integrable between $0$ and $+\infty$, which concludes the existence of a limit $\x_\infty$.
\end{proof}

\appendix

\section{Rigorous proof of Proposition \ref{prop:H1dot}}\label{appendix: technical proof}
In order to prove Proposition~\ref{prop:H1dot}, we begin with a preliminary result.
\begin{lem}\label{lem:uglylim}
	For any $f,g\in C(\R;\R)$ not necessarily with continuous derivative at the origin, we have
	\begin{align*}
		\lim_{\eps\to 0} \int f_y^\eps (\partial_y K^\eps\ast(\sy g))\eta_\infty\dy
		= 2g(0)\llangle f_y\rrangle + \int \sign(y)f_y(y)g_y(y)\eta_\infty(y)\dy\,.
	\end{align*}
\end{lem}
\begin{proof}
	By direct computation,
	\begin{align*}
		&\int f_y^\eps (\partial_y K^\eps\ast(\sy g))\eta_\infty\dy
		= \int f_y^\eps ( K^\eps\ast(2\delta_0g+\sy g_y))\eta_\infty\dy\\
		&\quad=\frac{1}{\eps^2}\iiint K\left(\frac{y-z_1}{\eps}\right)f_y(z_1) K\left(\frac{y-z_2}{\eps}\right) (2\delta_0(z_2)g(z_2)+\sign(z_2) g_y(z_2))\eta_\infty(y)\dy\dz_1\dz_2\\
		&\quad=\frac{2g(0)}{\eps^2}\iint K\left(\frac{y-z_1}{\eps}\right)f_y(z_1) K\left(\frac{y}{\eps}\right) \eta_\infty(y)\dy\dz_1\\
		&\quad \quad+ \frac{1}{\eps^2}\iiint K\left(\frac{y-z_1}{\eps}\right)f_y(z_1) K\left(\frac{y-z_2}{\eps}\right) \sign(z_2) g_y(z_2)\eta_\infty(y)\dy\dz_1\dz_2\\
		&\quad=2g(0)\iint K\left(z_1\right)f_y(\eps(y-z_1)) K\left(y\right) \eta_\infty(\eps y)\dy\dz_1\\
		&\quad \quad+ \frac{1}{\eps^2}\iiint K\left(\frac{y-z_1}{\eps}\right)f_y(z_1) K\left(\frac{y-z_2}{\eps}\right) \sign(z_2) g_y(z_2)\eta_\infty(y)\dy\dz_1\dz_2\,.
	\end{align*}
	The result then follows by taking the limit $\eps\to0$.
\end{proof}

Using this preliminary result, we are now ready to prove Proposition~\ref{prop:H1dot}

\begin{proof}[Proof of Proposition~\ref{prop:H1dot}]
	Convolving system \eqref{kinetic10} with $K^\eps$ in the space variable, and then differentiating, we obtain
	\begin{subequations}
		\label{kinetic10-eps}
		\begin{align}
			\partial_t \partial_y u^\eps + (1-\xdot) \partial_{yy} u^\eps &= -2\xdot\chi \partial_y K^\eps\ast\left(\sy u\right) - \partial_y(u^\eps - v^\eps)\notag\\
			&\quad +\chi \partial_y K^\eps\ast\left(\sy (u-v)\right) -2\xdot\chi \partial_y K^\eps\ast\sy\,, 
			\label{kinetic10a-eps}\\
			\partial_t \partial_y v^\eps - (1+\xdot) \partial_{yy} v^\eps &= -2\xdot\chi\partial_y K^\eps\ast(\sy v) + \partial_y (u^\eps-v^\eps)\notag\\
			&\quad +\chi \partial_y K^\eps\ast\left(\sy (u-v)\right)  -2\xdot\chi \partial_y K^\eps\ast\sy\,.
			\label{kinetic10b-eps}
		\end{align}
	\end{subequations}
	We then compute directly for $W_y^\eps:=\partial_y K^\eps\ast W$,
	\begin{align*}
		\frac12\frac{\dd}{\dt}
		\|W_y^\eps\|^2
		&= \frac12\frac{\dd}{\dt} \int \left(|u_y^\eps|^2+|v_y^\eps|^2\right)\eta_\infty\dy
		= \int \left(u_y^\eps\partial_t u_y^\eps+v_y^\eps\partial_t v^\eps\right)\eta_\infty\dy\\
		&= -\frac{(1-\xdot)}{2} \int \partial_y |u_y^\eps|^2\eta_\infty\dy
		+\frac{(1+\xdot)}{2} \int \partial_y |v_y^\eps|^2\eta_\infty\dy
		- \int |u_y^\eps-v_y^\eps|^2\eta_\infty\dy\\
		&\quad +\chi\int \left[ (u_y^\eps+v_y^\eps)  \partial_y K^\eps\ast\left(\sy (u-v)\right)  \right]\eta_\infty\dy\\
		&\quad -2\chi\xdot\int \left[ (u_y^\eps+v_y^\eps)  \partial_y K^\eps\ast\sy \right]\eta_\infty\dy\\
		&\quad -2\chi\xdot \int \left[ u_y^\eps\partial_y K^\eps\ast\left(\sy u\right)+ v_y^\eps\partial_y K^\eps\ast\left(\sy v\right)\right]\eta_\infty\dy\\
		&= -\frac{(1-\xdot)}{2} \int \partial_y |u_y^\eps|^2\eta_\infty\dy
		+\frac{(1+\xdot)}{2} \int \partial_y |v_y^\eps|^2\eta_\infty\dy
		- \int |u_y^\eps-v_y^\eps|^2\eta_\infty\dy\\
		&\quad +\chi I_1^\eps(u,v)-2\chi\xdot I_2^\eps(u,v)  -2\chi\xdot I_3^\eps(u,v)\,.
	\end{align*}
	Note that the last three terms $I_1^\eps, I_2^\eps, I_3^\eps$ have a nice structure for which we can apply Lemma~\ref{lem:uglylim}. For $I_1^\eps$, we have $f=u+v$, $g=u-v$, and so obtain the limit
	\begin{align*}
		\lim_{\eps\to0} I_1^\eps(u,v)&=\lim_{\eps\to0}  \int \left[ (u_y^\eps+v_y^\eps)  \partial_y K^\eps\ast\left(\sy (u-v)\right)  \right]\eta_\infty\dy\\
		&= 2(u(0)-v(0))\llangle u_y+v_y \rrangle + \int \sign(y)(|u_y(y)|^2-|v_y(y)|^2)\eta_\infty(y)\dy\\
		&= 2(u(0)-v(0))\llangle u_y+v_y \rrangle + \frac{1}{2\chi}\int \partial_y(|u_y(y)|^2-|v_y(y)|^2)\eta_\infty(y)\dy\,.
	\end{align*}
	where we integrated by parts in the last line using $\partial_y \eta_\infty(y)=-2\chi\sy\eta_\infty(y)$.
	Similarly, with $f=u+v$, $g=1$, we have for $I_2^\eps$ that
	\begin{align*}
		\lim_{\eps\to0} I_2^\eps(u,v)&=\lim_{\eps\to0}  \int \left[ (u_y^\eps+v_y^\eps)  \partial_y K^\eps\ast\sy  \right]\eta_\infty\dy
		= 2\llangle u_y+v_y \rrangle \,.
	\end{align*}
	Finally, for $I_3^\eps$, we have either $f=g=u$ or $f=g=v$, and so
	\begin{align*}
		\lim_{\eps\to0} I_3^\eps(u,v)&=\lim_{\eps\to0}  \int\left[ u_y^\eps\partial_y K^\eps\ast\left(\sy u\right)+ v_y^\eps\partial_y K^\eps\ast\left(\sy v\right)\right]\eta_\infty\dy\\
		&= 2u(0)\llangle u_y\rrangle+ 2v(0)\llangle v_y\rrangle + \int \sign(y)(|u_y(y)|^2+|v_y(y)|^2)\eta_\infty(y)\dy\\
		&= 2u(0)\llangle u_y\rrangle+ 2v(0)\llangle v_y\rrangle + \frac{1}{2\chi}\int \partial_y(|u_y(y)|^2+|v_y(y)|^2)\eta_\infty(y)\dy\,.
	\end{align*}
	Putting the above expressions together, we conclude for the limit $\eps\to 0$,
	\begin{align*}
		\frac12\frac{\dd}{\dt} \|W_y\|^2
		&= \lim_{\eps\to0} \frac12\frac{\dd}{\dt} \|W_y^\eps\|^2\\
		&= -\frac{(1-\xdot)}{2} \int \partial_y |u_y|^2\eta_\infty\dy
		+\frac{(1+\xdot)}{2} \int \partial_y |v_y|^2\eta_\infty\dy
		- \int |u_y-v_y|^2\eta_\infty\dy\\
		&\quad + 2\chi (u(0)-v(0))\llangle u_y+v_y \rrangle + \frac{1}{2}\int \partial_y(|u_y(y)|^2-|v_y(y)|^2)\eta_\infty(y)\dy\\
		&\quad-4\chi\xdot \llangle u_y+v_y \rrangle 
		-4\chi\xdot \left( u(0)\llangle u_y\rrangle+ v(0)\llangle v_y\rrangle \right)\\
		&\quad
		- \xdot \int \partial_y(|u_y(y)|^2+|v_y(y)|^2)\eta_\infty(y)\dy\\
		&=  - \int |u_y-v_y|^2\eta_\infty\dy- \frac{\xdot}{2} \int \partial_y \left(|u_y|^2+|v_y|^2\right)\eta_\infty\dy\\
		&\quad + 2\chi (u(0)-v(0))\llangle u_y+v_y \rrangle -4\chi\xdot \left[ \llangle u_y+v_y \rrangle+    u(0)\llangle u_y\rrangle+ v(0)\llangle v_y\rrangle \right]\,.
	\end{align*}
	This concludes the proof.
\end{proof}


\section{Method of characteristics}\label{sec:char-method}
In this section our main goal is to control the nonlocal terms $\llangle u_y\rrangle$ and $\llangle v_y\rrangle$ appearing as part of the entropy decay in Proposition~\ref{prop: gronwall_hypo_general_nonlocal} for the case $\alpha=0$. This is achieved in Proposition~\ref{prop:badterms} which we will prove here. We begin by rewriting system~\eqref{kinetic10} for 
\begin{align*}
	\tilde u(y;t):= e^t\eta_\infty(y)u(y;t)\,,\qquad 
	\tilde v(y;t):= e^t\eta_\infty(y)v(y;t)\,,
\end{align*}
yielding the new system
\begin{subequations}
	\label{kinetic15}
	\begin{align}
		\label{kinetic15a}
		\partial_t\tilde u + (1-\xdot) \partial_y\tilde u &=
		\tilde v -\chi \sy\left(\tilde u + \tilde v + 2\xdot e^t\eta_\infty\right)\,, \\
		\label{kinetic15b}
		\partial_t \tilde v - (1+\xdot) \partial_y\tilde v &= 
		\tilde u + \chi\sy\left(\tilde u + \tilde v - 2\xdot e^t\eta_\infty\right)\,.
	\end{align}
\end{subequations}
Notice that 
$\partial_y u(y;t)=\partial_y\tilde u(y;t)e^{-t}\eta_\infty^{-1} + 2\chi\sy  u(y;t)$,
and therefore 
\begin{align*}
	\llangle u_y\rrangle = e^{-t}\llangle \tilde u_y \eta_\infty^{-1}\rrangle 
	+2\chi\llangle \sy u\rrangle=  e^{-t}\llangle \tilde u_y\rrangle
\end{align*}
thanks to the fact that $u(y;t)$ is continuous at $y=0$. Similarly, $ \llangle v_y\rrangle =  e^{-t}\llangle \tilde v_y\rrangle$. This allows us to work with $\llangle \tilde u_y\rrangle$ and $\llangle \tilde v_y\rrangle$ directly, which we compute explicitly using the method of characteristics.

\begin{lem}[Method of Characteristics]\label{lem:method-of-char}
	If $(\tilde u, \tilde v)$ solve \eqref{kinetic15}, then they have the following integral representations.
	\begin{align*}
		\tilde u(y;t) &= \tilde u(y-t+x(t)-x(0);0)
		+\int_0^t\tilde v(y-s+x(t)-x(t-s);t-s)\ds\\
		&\quad-\chi\int_0^t\sign(y-s+x(t)-x(t-s))\tilde g(y-s+x(t)-x(t-s);t-s)\ds\,,\\
		\tilde v(y;t) &= \tilde v(y+t+x(t)-x(0);0)
		+\int_0^t\tilde u(y+s+x(t)-x(t-s);t-s)\ds\\
		&\quad+\chi\int_0^t\sign(y+s+x(t)-x(t-s))\tilde h(y+s+x(t)-x(t-s);t-s)\ds\,,
	\end{align*}
	where
	\begin{align*}
		\tilde g (y;t):=\tilde u(y;t) + \tilde v(y;t) + 2\xdot(t) e^t\eta_\infty(y)\,,\quad 
		\tilde h (y;t):=\tilde u(y;t) + \tilde v(y;t) - 2\xdot(t) e^t\eta_\infty(y)\,.
	\end{align*}
\end{lem}

\begin{proof}
	We consider the particle trajectories 
	\begin{align*}
		\dot Y(s)=\pm 1 -\xdot (s)\,,
	\end{align*}
	which for any $s\ge 0$ has solution 
	\begin{align*}
		Y(s)=Y(0) \pm s +x(0) - x(s)\,.
	\end{align*}
	Integrating solutions to \eqref{kinetic15a} along $Y(s)$ using "$+$" in the trajectory dynamics, we have
	\begin{align*}
		\tilde u(Y(t);t)=\tilde u(Y(0);0)
		+\int_0^t\tilde v(Y(s);s)\ds -\chi\int_0^t\sign(Y(s))\tilde g(Y(s);s)\ds
	\end{align*}
	Fix $t\ge0$, $y\in\R$ and choosing the initial condition $Y(0)=y-t+x(t)-x(0)$, we have
	$$
	Y(s)=y+s-t+x(t) - x(s)\,,
	$$
	and so $Y(t)=y$. Substituting into the expression for $\tilde u (Y(t);t)$, and changing variables $s\mapsto t-s$, we obtain the first equation. The second equation for $\tilde v(y;t)$ follows in the same manner, using instead the "$-$" case in the trajectories.
\end{proof}

With these preliminary computations, we are now ready to present the proof of Proposition~\ref{prop:badterms}.
\begin{proof}[Proof of Proposition~\ref{prop:badterms}]
First, notice that the sign of $z_\pm(s,t)$ is a direct consequence of Taylor expanding $x(t)=x(t-s)+s\xdot(\xi)$ with $\xi\in (t-s,t)$ and using Assumption~\ref{ass:xdot-bound}. 
Next, we will show the bound on $\llangle u_y\rrangle$ in detail; the estimate for $\llangle v_y\rrangle$ follows in a similar manner. In order to control $\llangle \tilde u_y\rrangle$, we directly differentiate the integral representation in Lemma~\ref{lem:method-of-char} for $y\neq 0$,
	\begin{align}
		&\partial_y\tilde u(y;t)
		=\partial_y\tilde u(y-t+x(t)-x(0);0)
		+\int_0^t\partial_y\tilde v(y-s+x(t)-x(t-s);t-s)\ds\notag\\
		&\quad-\chi\int_0^t\sign(y-s+x(t)-x(t-s))\partial_y\tilde g(y-s+x(t)-x(t-s);t-s)\ds -2\chi I(y;t)\,,\label{badterms-est1}
	\end{align}
	where
	\begin{align*}
		I(y;t):= \int_0^t \delta_0(y-s+x(t)-x(t-s))\tilde g(y-s+x(t)-x(t-s);t-s)\ds
	\end{align*}
	characterises precisely the jump discontinuity of $\partial_y\tilde u$ at $y=0$. Changing variables to $\tau:= y-s+x(t)-x(t-s)$, we have
	\begin{align*}
		I(y;t) &=\int_y^{y-t+x(t)-x(0)} \delta_0(\tau)\tilde g(\tau;t-s(y;\tau,t))\left(\frac{1}{-1+\xdot(t-s(y;\tau,t))}\right) \dd \tau    \\
		&=\int_{y-t+x(t)-x(0)}^y \delta_0(\tau)\tilde g(\tau;t-s(y;\tau,t))\left(\frac{1}{1-\xdot(t-s(y;\tau,t))}\right) \dd \tau    \,.
	\end{align*}
	At $\tau=0$, we define $s^*(y;t):=s(y;0,t)$ solving $s^*=y+x(t)-x(t-s^*)$, if it exists. Firstly, note that $I(y;t)=0$ if $s^*(y;t)<0$ or $s^*(y;t)>t$. Consider the case $0\le s^*(y;t)\le t \le T$. From Assumption~\ref{ass:xdot-bound}, we deduce
 \begin{align*}
     x(t)-x(t-s^*(y;t))=\int_{t-s^*(y;t)}^t \xdot(\bar s) \dd\bar s \in [-ps^*(y;t),+ps^*(y;t)] 
 \end{align*}
 Hence 
	\begin{align*}
		\frac{y}{1+p}\le s^*(y;t)\le \frac{y}{1-p}\,,
		\quad \text{ which implies }\quad
		\lim_{y \to 0} s^*(y;t)=0\,,
	\end{align*}
	and so $\sign(s^*(y;t))=\sign(y)$. It follows that no such solution $s^*\ge 0$ exists if $y<0$, meaning $I(y;t)=0$ for all $y<0$, and so $I(0^-;t)=\lim_{y\to 0^-}I(y;t)=0$. For $y>0$ small enough, the Dirac delta has support in the integration range, and taking the limit $y\to 0^+$ we have
	\begin{align*}
		\lim_{y\to 0^+} I(y;t) = \lim_{y\to 0^+} 
		\frac{\tilde g(0;t-s^*(y;t))}{1-\xdot(t-s^*(y;t))}=\frac{\tilde g(0;t)}{1-\xdot(t)} = I(0^+;t)\,.
	\end{align*}
	Considering that all other terms in \eqref{badterms-est1} are continuous at $y=0$, we obtain
	\begin{align*}
		&\llangle \tilde u_y \rrangle = \frac12 \left(\tilde u_y(0^+;t)+\tilde u_y(0^-;t)\right)\\
		& =\partial_y\tilde u(-t+x(t)-x(0);0)
		+\int_0^t\partial_y\tilde v(-s+x(t)-x(t-s);t-s)\ds\notag\\
		&\quad-\chi\int_0^t\sign(-s+x(t)-x(t-s))\partial_y\tilde g(-s+x(t)-x(t-s);t-s)\ds -\chi I(0^+;t)\\
		& =\partial_y\tilde u(z_-(t,t);0)
		+\int_0^t\partial_y\tilde v(z_-(s,t);t-s)\ds
		\quad+\chi\int_0^t\partial_y\tilde g(z_-(s,t);t-s)\ds -\chi I(0^+;t)\,, 
	\end{align*}
	where we used that $z_-(s,t)=-s+x(t)-x(t-s)<0$. Unwrapping the $\backsim$ notation, we have
\begin{align*}
	\partial_y\tilde u(z_-(t,t);0) &=  \eta_\infty(z_-(t,t))
 \left[u_y(z_-(t,t);0)-2\chi\sign(z_-(t,t))u(z_-(t,t);0)\right]\\ 
	&= \eta_\infty(z_-(t,t)) \left[u_y(z_-(t,t);0)+2\chi u(z_-(t,t);0)\right]\,,\\
	\partial_y\tilde v(z_-(s,t);t-s)&= e^{t-s} \eta_\infty(z_-(s,t)) \left[v_y(z_-(s,t);t-s)-2\chi\sign(z_-(s,t))v(z_-(s,t);t-s)\right]\\
	&= e^{t-s} \eta_\infty(z_-(s,t)) \left[v_y( z_-(s,t);t-s)+2\chi v(z_-(s,t);t-s)\right]\,,\\
	\partial_y\tilde g(z_-(s,t);t-s)&= e^{t-s} \eta_\infty(z_-(s,t)) \left[(u_y+v_y)(z_-(s,t);t-s)\right.\\
 &\hspace{3cm}\left.-2\chi\sign(z_-(s,t))(u+v)(z_-(s,t);t-s)\right]\\
	&\quad -4\chi\xdot(t-s)e^{t-s} \eta_\infty(z_-(s,t))\sign(z_-(s,t))\\
	&= e^{t-s} \eta_\infty(z_-(s,t)) \left[(u_y+v_y)(z_-(s,t);t-s)+2\chi(u+v)(z_-(s,t);t-s)\right]\\
	&\quad +4\chi\xdot(t-s)e^{t-s} \eta_\infty(z_-(s,t))\,.
\end{align*}
	We conclude that
	\begin{align*}
		\llangle  u_y \rrangle &= e^{-t} \llangle \tilde u_y \rrangle\\
		& =e^{-t}\eta_\infty(z_-(t,t)) \left[u_y(z_-(t,t);0)+2\chi u(z_-(t,t);0)\right]\\
		&\quad+\int_0^t e^{-s} \eta_\infty(z_-(s,t)) \left[v_y(z_-(s,t);t-s)+2\chi v(z_-(s,t);t-s)\right]\ds\\
		&\quad+\chi\int_0^t
		e^{-s} \eta_\infty(z_-(s,t)) \left[(u_y+v_y)(z_-(s,t);t-s)+2\chi(u+v)(z_-(s,t);t-s)\right]\ds\\
		&\quad+4\chi^2\int_0^t\xdot(t-s)e^{-s} \eta_\infty(z_-(s,t))\ds
		-\chi e^{-t}I(0^+;t)\,.  
	\end{align*}
	We proceed similarly for $\llangle  v_y \rrangle$. Differentiating the integral representation of $\tilde v$ in Lemma~\ref{lem:method-of-char} for $y\neq 0$, we obtain
	\begin{align}
		&\partial_y\tilde v(y;t)
		=\partial_y\tilde v(y+t+x(t)-x(0);0)
		+\int_0^t\partial_y\tilde u(y+s+x(t)-x(t-s);t-s)\ds\notag\\
		&\quad+\chi\int_0^t\sign(y+s+x(t)-x(t-s))\partial_y\tilde h(y+s+x(t)-x(t-s);t-s)\ds +2\chi J(y;t)\,,\label{badterms-est2}
	\end{align}
	where
	\begin{align*}
		J(y;t):= \int_0^t \delta_0(y+s+x(t)-x(t-s))\tilde h(y+s+x(t)-x(t-s);t-s)\ds
	\end{align*}
	characterises the jump discontinuity of $\partial_y\tilde v$ at $y=0$. Changing variables to $\tau:= y+s+x(t)-x(t-s)$, we have
	\begin{align*}
		J(y;t) &=\int_y^{y+t+x(t)-x(0)} \delta_0(\tau)\tilde h(\tau;t-s(y;\tau,t))\left(\frac{1}{1+\xdot(t-s(y;\tau,t))}\right) \dd \tau  \,.
	\end{align*}
	At $\tau=0$, we define $s^*(y;t):=s(y;0,t)$ solving $s^*=-y-x(t)+x(t-s^*)$, if it exists.
 Again, noting that $J(y;t)=0$ if $s^*(y;t)<0$ or $s^*(y;t)>t$, we consider the case $0\le s^*(y;t)\le t \le T$. From Assumption~\ref{ass:xdot-bound}, we deduce
 \begin{align*}
     x(t)-x(t-s^*(y;t))=\int_{t-s^*(y;t)}^t \xdot(\bar s) \dd\bar s \in [-ps^*(y;t),+ps^*(y;t)] 
 \end{align*}
 Hence 
	\begin{align*}
		\frac{-y}{1+p}\le s^*(y;t)\le \frac{-y}{1-p}\,,
		\quad \text{ which implies }\quad
		\lim_{y \to 0} s^*(y;t)=0\,,
	\end{align*}
	and so $\sign(s^*(y;t))=-\sign(y)$. It follows that no such solution $s^*\ge 0$ exists if $y>0$, meaning $J(y;t)=0$ for all $y>0$, and so $J(0^+;t)=\lim_{y\to 0^+}J(y;t)=0$. For $y<0$ close enough to zero, the Dirac delta has support in the integration range, and taking the limit $y\to 0^-$ we have
	\begin{align*}
		\lim_{y\to 0^-} J(y;t) = \lim_{y\to 0^-} 
		\frac{\tilde h(0;t-s^*(y;t))}{1+\xdot(t-s^*(y;t))}=\frac{\tilde h(0;t)}{1+\xdot(t)} = J(0^-;t)\,.
	\end{align*}
	Considering that all other terms in \eqref{badterms-est2} are continuous at $y=0$, we obtain
	\begin{align*}
		&\llangle \tilde v_y \rrangle = \frac12 \left(\tilde v_y(0^+;t)+\tilde v_y(0^-;t)\right)\\
		& =\partial_y\tilde v(t+x(t)-x(0);0)
		+\int_0^t\partial_y\tilde u(s+x(t)-x(t-s);t-s)\ds\notag\\
		&\quad+\chi\int_0^t\sign(s+x(t)-x(t-s))\partial_y\tilde h(s+x(t)-x(t-s);t-s)\ds +\chi J(0^-;t)\\
		& =\partial_y\tilde v(z_+(t,t);0)
		+\int_0^t\partial_y\tilde u(z_+(s,t);t-s)\ds
		+\chi\int_0^t\partial_y\tilde h(z_+(s,t);t-s)\ds +\chi J(0^-;t)\,,
	\end{align*}
	where $z_+(s,t)=s+x(t)-x(t-s)>0$. Unwrapping the $\backsim$ notation, we have
	\begin{align*}
		\partial_y\tilde v(z_+(t,t);0)&= \eta_\infty(z_+(t,t)) \left[v_y(z_+(t,t);0)-2\chi\sign(z_+(t,t))v(z_+(t,t);0)\right]\,,\\ 
		&= \eta_\infty(z_+(t,t)) \left[v_y(z_+(t,t);0)-2\chi v(z_+(t,t);0)\right]\,,\\
		\partial_y\tilde u(z_+(s,t);t-s)&= e^{t-s} \eta_\infty(z_+(s,t)) \left[u_y(z_+(s,t);t-s)-2\chi\sign(z_+(s,t))u(z_+(s,t);t-s)\right]\,,\\
		&= e^{t-s} \eta_\infty(z_+(s,t)) \left[u_y(z_+(s,t);t-s)-2\chi u(z_+(s,t);t-s)\right]\,,\\
		\partial_y\tilde h(z_+(s,t);t-s)&= e^{t-s} \eta_\infty(z_+(s,t)) \left[(u_y+v_y)(z_+(s,t);t-s)-2\chi\sign(z_+(s,t))(u+v)(z_+(s,t);t-s)\right]\\
		&\quad +4\chi\xdot(t-s)e^{t-s} \eta_\infty(z_+(s,t))\sign(z_+(s,t))\\
		&= e^{t-s} \eta_\infty(z_+(s,t)) \left[(u_y+v_y)(z_+(s,t);t-s)-2\chi(u+v)(z_+(s,t);t-s)\right]\\
		&\quad +4\chi\xdot(t-s)e^{t-s} \eta_\infty(z_+(s,t))\,.
	\end{align*}
	We conclude that
	\begin{align*}
		\llangle  v_y \rrangle &= e^{-t} \llangle \tilde v_y \rrangle\\
		& = e^{-t}\eta_\infty(z_+(t,t)) \left[v_y(z_+(t,t);0)-2\chi v(z_+(t,t);0)\right]\\
		&\quad+\int_0^t e^{-s}\eta_\infty(z_+(s,t)) \left[u_y(z_+(s,t);t-s)-2\chi u(z_+(s,t);t-s)\right]\ds\\
		&\quad+\chi\int_0^t
		e^{-s} \eta_\infty(z_+(s,t)) \left[(u_y+v_y)(z_+(s,t);t-s)-2\chi(u+v)(z_+(s,t);t-s)\right]\ds\\
		&\quad+4\chi^2\int_0^t\xdot(t-s)e^{-s} \eta_\infty(z_+(s,t))\ds
		+\chi e^{-t} J(0^-;t)\,.  
	\end{align*}
\end{proof}

\subsection{From time to space}\label{sec:time-to-space}

 This section is devoted to an auxiliary lemma, that allows to turn integration in time into integration in space. It is used for the proof of Proposition~\ref{prop:badterms-est}.

\begin{lem}\label{lem:badterms-aux1}
	Fix $\mu>0$, and let Assumption~\ref{ass:xdot-bound} hold for some $T>0$ and $p\in(0,1)$. Consider a function $f:\R\times\R_{\ge0} \to \R$ with enough regularity for the expressions below to make sense. Then
	\begin{align*}
		\int_0^{T} \int_0^t e^{-\mu s}\left|f(z_\pm(s,t);t-s)\right|\ds\dt
		\le\frac{1}{1-p} \int_0^{T} \int_\R e^{-\frac{\mu}{1+p} |y|}\left|f(y;t)\right|\dy\dt\,.
	\end{align*}
\end{lem}

\begin{proof}
	Exchanging the order of integration, then changing variables to $\tau:=t-s$, followed by another exchange of order of integration, we have
	\begin{align*}
		&\int_0^{T} \int_0^t e^{-\mu s}\left|f(z_\pm(s,t);t-s)\right|\ds\dt
		=  \int_0^{T} \int_s^{T} e^{-\mu s}\left|f(\pm s+x(t)-x(t-s);t-s)\right|\dt \ds\\
		&\qquad =  \int_0^{T} \int_0^{{T}-s} e^{-\mu s}\left|f(\pm s+x(\tau+s)-x(\tau);\tau)\right|\dd \tau \ds\\ 
		&\qquad =  \int_0^{T} \int_0^{{T-\tau}} e^{-\mu s}\left|f(\pm s+x(\tau+s)-x(\tau);\tau)\right|\ds\dd \tau\,.
	\end{align*}
	Next, we introduce $y(s):=\pm s+x(\tau+s)-x(\tau)$. Then $\dy=(\pm 1+\xdot(\tau+s))\ds$. Writing $s=s(y;\tau)$, the last expression is equal to
	\begin{align*}
		\int_0^{T} \int_0^{\pm (T-\tau)+x(T)-x(\tau)} e^{-\mu s(y;\tau)}\left|f(y;\tau)\right|\frac{1}{(\pm 1+\xdot(\tau+s(y;\tau)))}\dy\dd \tau\,.
	\end{align*}
	We use again that 
 $$x(\tau+s(y;\tau) )-x(\tau)=\int_\tau^{\tau+s(y;\tau) }\xdot(\bar s) \dd\bar s\in [-ps(y;\tau) ,+ps(y;\tau) ]$$ thanks to Assumption~\ref{ass:xdot-bound} since 
 $[\tau,\tau+s(y;\tau) ]\subset [0,T]$.
 Looking at the "$+$" case first, we notice that
	\begin{align*}
		(1-p)s(y;\tau)\le y \le (1+p)s(y;\tau) 
		\quad \Longrightarrow \quad
		\frac{y}{1+p}\le s(y;\tau)\le \frac{y}{1-p}\,.
	\end{align*}
	and $0\le T-\tau+x(T)-x(\tau)<\infty$.
	Hence
	\begin{align*}
		&\int_0^{T} \int_0^{{T-\tau}+x({T})-x(\tau)} e^{-\mu s(y;\tau)}\left|f(y;\tau)\right|\frac{1}{(1+\xdot(\tau+s(y;\tau)))}\dy\dd \tau\\
		&\qquad\le \frac{1}{1-p} \int_0^{T} \int_0^\infty e^{-\frac{\mu y}{1+p}}\left|f(y;\tau)\right|\dy\dd \tau
		\le \frac{1}{1-p} \int_0^{T} \int_\R e^{-\frac{\mu}{1+p} |y|}\left|f(y;\tau)\right|\dy\dd \tau\,.
	\end{align*}
	Which concludes the proof for the "$+$" case. For the "$-$" case, we have similarly,
	\begin{align*}
		-(1+p)s(y;\tau)\le y \le -(1-p)s(y;\tau) 
		\quad \Longrightarrow \quad
		\frac{-y}{1+p}\le s(y;\tau)\le \frac{-y}{1-p}
	\end{align*}
	and $-\infty< -(T-\tau)+x(T)-x(\tau)\le 0$.
	Hence
	\begin{align*}
		&\int_0^{T} \int_0^{-(T-\tau)+x(T)-x(\tau)} e^{-\mu s(y;\tau)}\left|f(y;\tau)\right|\frac{1}{(-1+\xdot(\tau+s(y;\tau)))}\dy\dd \tau\\
		&\qquad\le \frac{1}{1-p} \int_0^{T} \int_{-\infty}^0 e^{\frac{\mu y}{1+p}}\left|f(y;\tau)\right|\dy\dd \tau
		\le \frac{1}{1-p} \int_0^{T} \int_\R e^{-\frac{\mu}{1+p} |y|}\left|f(y;\tau)\right|\dy\dd \tau\,.
	\end{align*}
	
\end{proof}

\section*{Acknowledgements}
This project has received funding from the European Research Council (ERC) under the European Union’s Horizon 2020 research and innovation programme (grant agreement No 865711).
GF has been partially financed by the Austrian Science Fund (FWF) projects \href{https://doi.org/10.55776/W1245}{10.55776/W1245} and \href{https://doi.org/10.55776/F65}{10.55776/F65}, and by Gruppo Nazionale per l’Analisi Matematica, la Probabilità e le loro Applicazioni (GNAMPA) No. E55F22000270001, and by the INdAM project N.E53C22001930001 “MMEAN-FIELDSS”. GF also wants to thank the Department of Applied Mathematics and Mathematical Physics of the Urgench State University, and Christian Schmeiser for his mentoring and support.
FH is supported by start-up funds at the California Institute of Technology. FH was also supported by the Deutsche Forschungsgemeinschaft (DFG, German Research Foundation) via project 390685813 - GZ 2047/1 - HCM.


\end{document}